\documentclass[12pt, a4paper]{amsart}
\usepackage{amssymb}
\usepackage{t1enc}
\usepackage[mathscr]{eucal}
\usepackage{xcolor}
\usepackage{tikz}
\usepackage{tikz-cd}
\usepackage{hyperref}
\usepackage[capitalize,noabbrev,nameinlink]{cleveref}
\numberwithin{equation}{section}
\usepackage[margin=1in]{geometry}
\usepackage{booktabs}
\usepackage{enumerate}
\usepackage{bbm}
\usepackage{bm}
\usepackage{xfrac}
\usepackage{comment}
\crefname{thm}{Theorem}{Theorems}
\crefname{lem}{Lemma}{Lemmas}
\crefname{cor}{Corollary}{Corollaries}
\crefname{prop}{Proposition}{Propositions}
\usepackage{listings}
\theoremstyle{plain}
\newtheorem{maintheorem}{Theorem}

\newtheorem{thm}{Theorem}[section]
\newtheorem{lem}[thm]{Lemma}
\newtheorem{cor}[thm]{Corollary}
\newtheorem{prop}[thm]{Proposition}

 \theoremstyle{definition}
\newtheorem{defn}[thm]{Definition}
\newtheorem{rem}[thm]{Remark}

\newtheorem{qst}{Question}

\let\H\undefined

\newcommand{\mb}[1]{\mathbb{#1}}
\newcommand{\mc}[1]{\mathcal{#1}}

\newcommand{\mr}[1]{\mathrm{#1}}

\newcommand{\Var}{\mathrm{Var}}
\newcommand{\End}{\operatorname{End}}
\newcommand{\SH}{\operatorname{SH}}
\newcommand{\GW}{\operatorname{GW}}
\newcommand{\Spec}{\operatorname{Spec}}
\newcommand{\Cr}{\mathrm{Cr}}
\newcommand{\1}{\mathbbm{1}}

\newcommand{\id}{\mathrm{id}}
\newcommand{\K}{K}
\newcommand{\Einfty}{\mathbb{E}_{\infty}}
\newcommand{\H}{\operatorname{H}}
\newcommand{\kbar}{\overline{k}}

\newcommand{\KQ}{KQ}
\newcommand{\Z}{\mathbb{Z}}
\newcommand{\sn}{\operatorname{sn}}

\newcommand{\Aut}{\operatorname{Aut}}

\newcommand{\mot}{\mathrm{mot}}

\newcommand{\Tr}{\operatorname{Tr}}
\newcommand{\ev}{\mathrm{ev}}
\newcommand{\et}{\acute{e}t}
\newcommand{\coev}{\mathrm{coev}}
\newcommand{\Hdg}{\mathrm{Hdg}}
\newcommand{\Gal}{\mathrm{Gal}}
\newcommand{\chikq}{\chi_1^{\KQ}}

\usepackage{mathtools}

\newcommand{\M}{\operatorname{M}}
\makeatletter
\@namedef{subjclassname@2020}{\textup{2020} Mathematics Subject Classification}

\makeatother

\begin{document}

\title{Motivic invariants of automorphisms and the plane Cremona group}

\author[Azouri]{Ran Azouri}
\address{Department of Mathematics, Faculty of Mathematics and Computer Science,  Weizmann Institute of Science} 
\email{ran.azouri@weizmann.ac.il}
\urladdr{https://ranazouri.com}

\author[Nanavaty]{Anubhav Nanavaty}
\address{Department of Mathematics \\ Cornell University} 
\email{nanavaty@cornell.edu}
\urladdr{https://anubhavn.github.io}

\author[Pajwani]{Jesse Pajwani}
\address{School of Mathematics, University of Bristol, Bristol, BS8
1TW, UK, and the Heilbronn Institute for Mathematical Research, Bristol, UK} 
\email{jesse.pajwani@bristol.ac.uk}
\urladdr{https://sites.google.com/view/jessepajwani/Jesse-Pajwani}

\subjclass[2020]{14E07, 14F42, 19B99, 19G38}
\keywords{Motivic homotopy theory, Cremona groups, $\K$-theory of varieties}

\begin{abstract}
Recent work of Campbell and Zakharevich studies the $\K$-theory spectrum of varieties. In previous work, the second author constructs a derived version of the $\mathbb{A}^1$-Euler characteristic, giving a map out of this spectrum. We analyse the morphism that the derived $\mathbb{A}^1$-Euler characteristic induces out of $\K_1(\Var_k)$. We obtain an explicit description of this map, which encodes information about automorphisms of varieties. We use Hermitian $\K$-theory to obtain a geometric description for part of this morphism in terms of Hodge cohomology. This allows us to obtain structural results about the group $\K_1(\Var_k)$. As an application, we use the $\Z/2\Z$-equivariant minimal model programme for surfaces to construct an arithmetic invariant of a large class of involutions in the plane Cremona group. We use this to construct two non-conjugate de Jonqui\`eres involutions with isomorphic positive genus fixed curves over any perfect field $k$ of characteristic not $2$ with at least 6 elements which admits a quadratic extension.

\end{abstract}

\maketitle
\tableofcontents

\section{Introduction}
Algebraic geometry often uses invariants associated to the category of varieties over a field $k$. A classical example of such an invariant when $k=\mathbb{C}$ is the \emph{compactly supported Euler characteristic}. For $X/\mathbb{C}$ a smooth projective variety, the Euler characteristic is defined to be the integer 
\[
\chi_0(X):=\sum_{i=0}^{\infty} (-1)^i\mathrm{dim}_{\mathbb{Q}} \H_{\rm{sing}}^i(X(\mathbb{C}); \mathbb{Q})\in\mb{Z},
\]
which extends to a ring homomorphism from the Grothendieck ring of varieties (Definition $\ref{defn:ringofvarieties}$) $\K_0(\Var_{\mathbb{C}}) \to \Z$.

Motivic homotopy theory allows us to apply techniques from homotopy theory to analyse varieties over general fields. The \emph{stable motivic homotopy category}, $\SH(k)$, is an algebro-geometric analogue of the stable homotopy category $\SH$ studied in algebraic topology, and there is a natural functor $\Sigma^{\infty}\colon \Var_k \to \SH(k)$, allowing us to study varieties by their image in $\SH(k)$. A recent trend in motivic homotopy theory is to use the rich arithmetic and geometric structure of this setting to construct arithmetic refinements of classical invariants coming from algebraic topology. Over a general perfect field $k$, Hoyois \cite{Hoy14} defines an arithmetic refinement of the Euler characteristic, known as the \emph{$\mathbb{A}^1$-Euler characteristic}. He builds on a categorical construction of Dold--Puppe \cite{DoldPuppe1999} and constructs this Euler characteristic as a categorical trace in the stable motivic homotopy category $\SH(k)$, taking values in the group of homotopy classes of endomorphisms of the motivic sphere spectrum, $\pi_0(\End_{\SH(k)}(\1_k))$. Work of Morel \cite[Theorem 6.4.1]{Morel2004} identifies $\pi_0(\End_{\SH(k)}(\1_k))=\GW(k)$, the Grothendieck--Witt ring of bilinear forms. Later, the $\mathbb{A}^1$-Euler characteristic was shown to extend to a ring homomorphism out of $\K_0(\Var_k)$ in \cite[Theorem 1.13]{AMBOWZ22}. The categorical definition does not lend itself to computations, however a theorem of Levine--Raksit \cite[Theorem 1.3]{LR20} gives a more computable definition. They show that if $X/k$ is a smooth projective variety, the $\mathbb{A}^1$-Euler characteristic $\chi_0(X)$ is given in $\GW(k)$ by taking a natural bilinear form furnished by the cup product pairing and coherent duality on the Hodge cohomology groups of $X$
\[
\chi_{0}^{\mb{A}^1}(X)=\Bigg[\bigoplus_{p,q}H^p(X,\Omega^q_{X/k})[q-p],\cup \Bigg]\in \GW(k).
\]
There is much arithmetic information contained in the $\mathbb{A}^1$-Euler characteristic, and the fact that it gives a homomorphism out of $\K_0(\Var_k)$ allows it to be manipulated easily. However, when passing from $\Var_k$ to $\K_0(\Var_k)$, one identifies which varieties can be `cut and pasted' together to form other ones, but we lose information about the different ways the pieces of the same variety fit together. Partly in an attempt to record such information, Zakharevich \cite{Zak17-assemblers} and later Campbell \cite{Campbell} constructed a $\K$-theory spectrum of varieties, $\K(\Var_k)$ (Definition $\ref{defn:KVark}$), an $\Einfty$-ring whose connected components are exactly $\K_0(\Var_k)$, and higher $\K$-groups of $\Var_k$ are defined as $\K_i(\Var_k) := \pi_i( \K(\Var_k))$. Importantly for this paper, \cite[Theorem A]{Annihilator} shows that the groups $\K_i(\Var_k)$ are intimately connected to birational geometry. The group $\K_1(\Var_k)$ has a known presentation due to \cite[Proposition 4.4]{KLMMS24} and \cite[Theorem B]{Ng}, and is closely related to automorphisms of varieties: any automorphism $f\colon X\to X$ defines an element which we denote $[X,f]\in \K_1(\Var_k)$ (Definition $\ref{def:SWK_1Aut}$). Despite having a presentation, little is known about the structure of the group $\K_1(\Var_k)$: for example it has been unknown whether $\K_1(\Var_k)$ is infinite for any field $k$.

The connection between the $\K$-theory spectrum of varieties and the $\mathbb{A}^1$-Euler characteristic was first shown by the second author \cite[Theorem 1.1]{Nan24} building on work of R\"ondigs \cite[Corollary 6.7]{Ron25}, where he proved that, in characteristic $0$ the $\mathbb{A}^1$-Euler characteristic $\chi_0$ lifts to a map of spectra
\[
\chi\colon \K(\Var_k) \to \End_{\SH(k)}(\1_k).
\]
In Definition $\ref{defn:derivedEulerChar}$, we extend the above definition to positive characteristic fields, so that if $k$ is a perfect field with exponential characteristic $e$, we obtain a map
\[
\chi\colon \K(\Var_k) \to \End_{\SH(k)[1/e]}(\1_k[1/e]).
\]

The main motivating question for the paper is the following.
\begin{qst}
    What arithmetic or geometric information is contained in the induced map $\chi_1 \colon \K_1(\Var_k) \to \pi_1\End_{\SH(k)[1/e]}(\1_k[1/e])$?
\end{qst}
In our first main theorem, we give a purely categorical description of this map. 
\begin{maintheorem}\label{maintheorem:Categorical}
    Let $X$ be a smooth, projective $k$-variety and $f\in\Aut(X)$. Then $\chi_1([X,f])$ is the element of $\pi_1(\End_{\SH(k)[1/e]}(\1_k[1/e]))$ induced by the $2$-commutative diagram (where, for notational ease, we write $X,f$ to denote the image of $\Sigma^{\infty} X_+, \Sigma^{\infty}f$ in $\SH(k)[1/e]$)
\begin{equation*}
			\begin{tikzcd}
				&|[alias=A]|  X \wedge X^\vee\arrow[rr,"\id\wedge\id"]\arrow[dd," f\wedge f^{\vee, -1}"]&& X \wedge X^\vee \arrow[rr,"\tau"]\arrow[ddll,Rightarrow]\arrow[dd,"f\wedge f^{\vee, -1}"] && X^\vee \wedge  X \arrow[dr,""{name=C,below,inner sep=1pt}]\arrow[dd,"f^{\vee, -1}\wedge f"']\arrow[ddll,Rightarrow]&\\
				\1_k[1/e]\arrow[ur]\arrow[dr,""{name=B,inner sep=1pt}] &&&&&& \1_k[1/e],\\
				& X \wedge X^\vee \arrow[rr,"\id\wedge\id"]&& X \wedge X^\vee\arrow[rr,"\tau"] && |[alias=D]| X^\vee \wedge X \arrow[ur]&
				\arrow[Rightarrow,from=A,to=B]
				\arrow[Rightarrow,from=C,to=D]
			\end{tikzcd}
		\end{equation*}
	where the $2$-morphisms in the above diagram are specified in Lemma $\ref{lem:trace}$.
\end{maintheorem}

In fact, the description above also applies to arbitrary varieties $p \colon X \to \Spec k$ if we replace $\Sigma^{\infty} X_+$ by the compactly supported motive $\M^c(X) := p_*p^! \1_k[1/e] $, and $\Sigma^{\infty} f_+$ by $\M^c(f)$, described by Equation~\ref{eq:M^c(f)}. While Theorem \ref{maintheorem:Categorical} provides a description of $\chi_1([X,f])$, it is not computable in practice. Although we have a description of the group $\pi_1\End_{\SH(k)[1/e]}(\1_k)[1/e]$ outside of characteristic $2,3$ due to \cite[(1.1)]{RSO19}, the above diagram is not easy to analyse using this description. To obtain computable information from $\chi_1([X,f])$, we compose $\chi_1$ with the unit map to the motivic spectrum representing Hermitian $\K$-theory $u\colon \pi_1(\End_{\SH(k)[1/e]}(\1_k)) \to \pi_{1,0}(\KQ[1/e])$. Because we are working with Hermitian $\K$-theory, from here, we require that $\mathrm{char}(k) \neq 2$. Writing $\chikq$ for the composition $u \circ \chi_1$, we can identify $\pi_{1,0}(\KQ[1/e])$ with $\pi_{1,0}(\KQ)$ (Corollary $\ref{cor:KQLocal}$), and we see that $\chikq([X,f])$ is an element of $\pi_{1,0}(\KQ)$. A theorem of Bass \cite[Corollary 4.7.7]{Bass} identifies $\pi_{1,0}(\KQ)$  with the group $\{\pm1\} \times k^{\times}/k^{\times 2}$. This allows us to prove the second main theorem of this paper, which can be thought of as a higher categorical version of the motivic Gauss--Bonnet formula \cite[Theorem 1.3]{LR20}. 

\begin{maintheorem}\label{maintheorem:Hermitian}
    Let $X$ be a smooth projective variety over $k$ and $f\colon X \to X$ be an automorphism. Let $f^*_{\Hdg}$ denote the pullback map on the Hodge cohomology of $X$. Let $\sn(f^*_{\Hdg})$ denote the spinor norm of $f^*_{\Hdg}$ (Definition $\ref{defn:spinornorm}$). Then 
    \[
  \chikq([X,f]) = ( \det(f^*_{\Hdg}), \sn(f^*_{\Hdg})) \in \{\pm 1\} \times k^{\times}/k^{\times 2}.
    \]
\end{maintheorem}

Unlike the categorical description from Theorem \ref{maintheorem:Categorical}, computing $( \det(f^*_{\Hdg}), \sn(f^*_{\Hdg}))$ is often relatively straightforward. This description allows us to compute $\chikq$ for many natural automorphisms of varieties. We obtain simple formulae for $\chikq$ of Bertini involutions on del Pezzo surfaces of degree $1$ (Corollary $\ref{cor:BertiniComputation}$), Geiser involutions on del Pezzo surfaces of degree $2$ (Lemma $\ref{lem:GeiserComputation}$), and involutions on conic bundle surfaces (Corollary $\ref{cor:conicbundlesn}$). These computations allow us to obtain new results about $\K_1(\Var_k)$. 

\begin{maintheorem}\label{maintheorem:K1Var}
    The group $\K_1(\Var_k)$ surjects onto $\{\pm1\} \times k^{\times}/k^{\times 2}$. In particular, if $k$ is finitely generated over $\mathbb{Q}$, then $\K_1(\Var_k)$ is not finitely generated.
\end{maintheorem}
Previous results about $\K_1(\Var_k)$ (see for example \cite[\S6.3]{CWZ19}) have shown that $\K_1(\Var_k)$ is often non-trivial and in certain cases surjects onto $\Z/2\Z$, or possibly is infinite after a localisation (\cite[Theorem 1.3]{BG21}), however Theorem \ref{maintheorem:K1Var} demonstrates that $\K_1(\Var_k)$ is infinite in a large class of natural cases. These computations also allow us to examine the $\Einfty$-ring structure on $\K(\Var_k)$ and examine the compatibility of $\chi$ with this structure (see Remark $\ref{rem:Einfty}$ for further discussion).\\
\textbf{Applications to the Plane Cremona Group}.  We give an application of $\chi_1$ to birational geometry by constructing a conjugacy invariant for involutions of the plane Cremona group $\mathrm{Cr}_2(k)$ (Definition $\ref{defn:PlaneCremona}$). The study of conjugacy classes of involutions was first initiated over the complex numbers by Bertini \cite{Bertini}, and completed by Bayle and Beauville \cite[Theorem 2.6, Proposition 2.7]{BayleBeauville}. A similar classification of conjugacy classes of involutions over other fields was explicitly asked for in \cite[\S9]{DolgachevIskovskikh}, and remains an open and an active area of research. For $k=\mathbb{R}$, this question has also been explored in \cite[Theorem 1.1]{Zimmermann2018}, and an almost complete classification has been obtained in \cite[Main Theorem]{CMYZ}. The approach in all cases has been to regularise birational involutions to $\mb{Z}/2\mb{Z}$-\emph{Mori fibre spaces} over $k$ (Definition \ref{defn:MFS}), and classify these up to $\mb{Z}/2\mb{Z}$-birational equivalence. We show that this approach fits naturally into the language of $\K_1(\Var_k)$, and examine our invariant $\chikq$ in the setting of birational geometry.
\begin{maintheorem}\label{maintheorem:birational}
 Let $f\in\Cr_2(k)$ be an involution such that $f\otimes\overline{k}\in \Cr_2(\overline{k})$ is not conjugate to an involution in $\mr{PGL}_3(\overline{k})$. Let $(X,f_X)$ be any $\mb{Z}/2\mb{Z}$-Mori fibre space regularisation of $f$. Then 
\[\chi_1^{\KQ}(f):=\chi_1^{\KQ}(X,f_X)\] is independent of the choice of $(X,f_X)$, and hence a conjugacy invariant of $f$. 
\end{maintheorem}
The proof of this theorem uses the $\mb{Z}/2\mb{Z}$-Sarkisov program for surfaces over a perfect field developed in \cite{DolgachevIskovskikh}  (Definition $\ref{def:linktypes}$). One sees $\mb{Z}/2\mb{Z}$-Sarkisov links give rise to relations in $\K_1(\Var_k)$, and we study how $\chi_1^{\KQ}$ behaves under these relations. 

A priori, an invariant valued in $k^\times/k^{\times2}$ may seem quite coarse, however, we give explicit computations which demonstrate that $\chikq$ records interesting arithmetic information tied to birational involutions, allowing us to distinguish their conjugacy classes. In Theorem \ref{thm:deJonquieres}, we construct an explicit family of degree $6$ de Jonqui\`eres involutions which are conjugate over $\kbar$, but such that the invariant $\chikq$ exhibits them as not conjugate over $k$.

In \cite[Proposition 2.7]{BayleBeauville}, it is shown that the conjugacy class of an involution in $\mathrm{Cr}_2(\mb{C})$ is completely determined by the isomorphism class of its normalised fixed curve. Over $k=\mb{R}$, \cite[Main Corollary]{CMYZ} constructs an uncountable family of non-conjugate involutions with isomorphic normalised fixed curves. For other fields $k$ however, it is an open question whether non-conjugate involutions with isomorphic fixed curves exist. Using Theorem \ref{maintheorem:birational}, and building on a construction from \cite[\S10]{CMYZ}, we construct such involutions whenever $k$ has at least one non-square.

\begin{maintheorem}\label{maintheorem:fixedcurve}
   Suppose that $k^{\times}/k^{\times 2}$ is non trivial and $k\neq \mb{F}_3, \mb{F}_5$. Then there exist two birational involutions $\phi, \phi' \in \Cr_2(k)[2]$ with isomorphic fixed curves of positive genus over $k$ that are not conjugate to each other. 
\end{maintheorem}
The above theorem is part of a larger story: invariants coming from motivic homotopy theory are particularly sensitive to the underlying arithmetic, and as such, they contribute to our understanding of birational geometry over non-algebraically closed fields.

\subsection{Outline of the paper}
Section $\ref{sec:prelims}$ of this paper is dedicated to recalling important notions for sections \ref{section:categorical} and \ref{section:Hermitian} of this paper, including the definition of $\K(\Var_k)$ and the construction of the spectrum level map $\chi\colon \K(\Var_k) \to \End_{\SH(k)}(\1_k)$. We also recall Ng's presentation of $\K_1(\Var_k)$ (\cite[Theorem B]{Ng}) to demonstrate the connection between automorphisms of varieties and $\K_1(\Var_k)$. The techniques used in this section are largely $\K$-theoretic and $\infty$-categorical.

In Section $\ref{section:categorical}$, we use the results of \cite{HSS17} to give a categorical description of the map $\chi_1\colon \K_1(\Var_k) \to \pi_1(\End_{\SH(k)[1/e]}(\1_k[1/e]))$. This section concludes with the proof of Theorem \ref{maintheorem:Categorical}. The techniques used in this section are again largely $\infty$-categorical.

Section $\ref{section:Hermitian}$ of this paper is dedicated to proving Theorem \ref{maintheorem:Hermitian} by using Hermitian $\K$-theory to understand the map constructed in Theorem \ref{maintheorem:Categorical}. We recall relevant notions of Hermitian $\K$-theory and its relation to $\SH(k)$, before using the motivic Gauss--Bonnet formula of Levine--Raksit, \cite[Theorem 1.3]{LR20} to obtain an explicit, geometric, description of the map constructed in Section $\ref{section:categorical}$. This section uses techniques from motivic homotopy theory.

Section $\ref{section:examples}$ is concerned with computing explicit examples. We show some elementary linear algebra results before using them to compute the value of the map $\chi_1$ for various canonical automorphisms of varieties. We also show transformation results about how $\chi_1$ behaves under blow-up operations. This section uses classical algebraic geometry techniques, especially those concerning del Pezzo surfaces.

In Section $\ref{section:K1Vark}$, we use the results of Section $\ref{section:examples}$ to prove Theorem \ref{maintheorem:K1Var} as well as other new results concerning the structure of $\K_1(\Var_k)$ and the interactions of the map $\chi$ with the $\Einfty$-ring structures on the source and target spectra. This section is largely $\K$-theoretic.

Section $\ref{section:birational}$ uses the results of Section $\ref{section:examples}$ in order to construct a new invariant for certain involutions in the plane Cremona group. We first recall relevant results from the equivariant minimal model programme which we reinterpret in the language of $\K_1(\Var_k)$. This allows us to show that equivariant birational equivalences
give rise to natural elements in $\K_1(\Var_k)$. We then use the results of Sections $\ref{section:Hermitian}$ and $\ref{section:examples}$ to give explicit conditions for when $\chi_1$ of these elements vanish. This section concludes by proving Theorem \ref{maintheorem:birational}, which provides us with a new, arithmetic invariant of involutions in the plane Cremona group. This section uses many results and techniques from birational geometry.

Finally in Section $\ref{sec:application}$, we use the results of Section $\ref{section:birational}$ to construct examples of non-conjugate elements in the plane Cremona group. We first adapt an example from \cite[Lemma 6.2]{CMYZ} to construct a family of non-conjugate, but geometrically conjugate, de Jonqui\`eres involutions. We end by proving Theorem \ref{maintheorem:fixedcurve}, which follows by adapting an example from \cite[\S10]{CMYZ} and applying Theorem \ref{maintheorem:birational}. This section uses classical algebraic geometry techniques, especially concerning conic bundle surfaces.
\subsection{Notation and conventions}
Throughout this article, we fix the following notation.
\begin{itemize}
    \item Let $k$ denote an arbitrary perfect field, which is of characteristic not $2$ from Section \ref{section:Hermitian} onwards.
    \item Let $e$ denote the exponential characteristic of $k$, so that $e=1$ if $\mathrm{char}(k)=0$ and $p$ if $\mathrm{char}(k) = p$.
    \item Let $k^\times$ denote the multiplicative group of $k$.
    \item Let $\kbar$ denote a fixed algebraic closure of $k$.
    \item Let $\Gal_k$ denote the absolute Galois group of $k$.
    \item Let $\mathrm{Var}_k$ denote the category of varieties over $k$. A \emph{variety} here will mean a reduced, separated scheme of finite type over $k$.
    \item For $X/k$ a smooth projective variety of dimension $n$, $\Omega^1_{X/k}$ will denote the \emph{cotangent bundle}, we will write $\Omega^i_{X/k} = \bigwedge^i \Omega^1_{X/k}$ and the \emph{canonical bundle} is defined to be $\omega_{X/k} := \Omega^n_{X/k}$.
    \item Let $\SH(X)$ denote the stable motivic homotopy category of $X$.
    \item For $f: X \to Y$ a morphism of varieties, let $f_*, f_!, f_{\#} \colon \SH(X) \to \SH(Y)$ and $f^*, f^!\colon \SH(Y) \to \SH(X)$ be the associated functors furnished by the six-functor formalism on $\SH(-)$ (see \cite[Theorem 1.1, Theorem 1.5]{Hoy17}).
    \item Let $\1_k$ denote the unit object of $\SH(k)$.
    \item Let $\wedge$ denote the smash product in $\SH(k)$.
    \item For $S,T \in \SH(k)$, write $\rm{Map}_{\SH(k)}(S,T)$ to mean the mapping spectrum. If $S=T$, we will write $\End_{\SH(k)}(S)$ to mean this mapping spectrum.
    \item Let $\KQ$ denote the Hermitian $\K$-theory spectrum in $\SH(k)$ constructed in \cite{PaninWalter2019} where it is written as $\mathbf{BO}$. It is often also referred to as $\mathrm{KO}$ such as in \cite{OrmsbyOstvaer}.
    \item Let $\GW(k)$ denote the \emph{Grothendieck--Witt ring} of $k$ (Definition $\ref{defn:GWk}$).
\end{itemize}

\subsection*{Acknowledgments}
The authors thank Ming Ng deeply for his continued support regarding understanding $\K_1(\Var_k)$. We thank Thor Wittich for helpful conversations about the structure of $\K_1$ of an assembler. We thank Susanna Zimmermann for helpful discussions regarding Sarkisov links and the minimal model programme. We thank Sam Streeter for discussions about the del Pezzo surfaces appearing in Section $\ref{section:examples}$, and Harry Shaw for pointing out an error in an earlier draft of Subsection $\ref{subsection:final}$. We thank Inna Zakharevich and Sidhanth Raman for numerous helpful discussions. We thank Thomas Brazelton and Frédéric Déglise for discussions on six-functor formalism. We thank the organisers of the 2024 Park City Mathematics Institute summer programme in motivic homotopy theory for bringing the authors together to begin this project, as well as the organisers of the 2026 SLMath special semester in motivic homotopy theory for bringing the authors together to get the project over the line. An early draft of sections 3 and 4 was written with Stephen McKean; both sections have entirely altered since.

The first author was supported by the ERC project `Foundation of motivic real K-theory' (ERC grant no. 949583). The third author was supported by the Heilbronn Institute for Mathematical Research.

The use of AI for this paper was limited to surveying background literature, checking the linear algebra computations at the beginning of Section $\ref{section:examples}$ and finding a suitable candidate for the sextic polynomial $f(z,t)$ in subsection \ref{subsection:example}. AI was not used to generate ideas for the constructions of this paper, prove results, or write any parts of the paper.
	\section{Preliminaries}\label{sec:prelims}
    \subsection{The $\K$-theory spectrum of varieties and derived motivic measures}
\begin{defn}\label{defn:ringofvarieties}
     The \emph{Grothendieck ring of varieties}, denoted $\K_0(\Var_k)$, is a ring where the generators are isomorphism classes $[X]$ where $X$ is a $k$-variety, subject to relations of the form $[X]=[Z]+[X\setminus Z]$ for any Zariski closed subvariety $Z\subset X$. The ring structure is induced by Cartesian product.
\end{defn}
A spectrum whose path components recovered $\K_0(\mr{Var}_k)$ was first constructed by Zakharevich in \cite[\S5.1]{Zak17-assemblers} using the framework of assemblers, and then Campbell in \cite[p.3]{Campbell} using SW-categories. While both frameworks produce homotopy equivalent spectra, we recall the framework of SW-categories as it is more suitable for this paper.
    \begin{defn}[{\cite[Definition 3.23]{Campbell}}] \label{def:CWcat}
        An \emph{SW-category} is a category $C$ equipped with three distinguished wide subcategories of cofibrations, complement maps, and weak equivalences. We write cofibrations with hooks $Z\hookrightarrow X$, and complements as $U\xrightarrow{\circ} X$. The category $C$ is further equipped with a collection of so-called ``subtraction sequences'' \[Z \hookrightarrow X \xleftarrow{o} U. \]
        This data is required to satisfy axioms spelled out in \cite[Definition 3.7, Definition 3.13, Definition 3.24]{Campbell}.
    \end{defn}
\begin{prop}[{\cite[Proposition 3.26]{Campbell}}]
    The category $\Var_k$ can be given the structure of an SW-category, where the cofibrations are closed immersions, complements are open immersions, and the weak equivalences are isomorphisms. The subtraction sequences are of the form $Z\hookrightarrow X \xleftarrow{\circ} U$ where the image of $Z$ is the set-theoretic complement of $U$. 
\end{prop}
\begin{defn}[{\cite[Definition 2.6]{CWZ19}}]\label{defn:KVark}
    If $\mathcal{C}$ is an SW-category, we define its \emph{$\K$-theory space} to be \[\K^t(\mathcal{C}):=\Omega|w\tilde{S}_\bullet\mathcal{C}|\]
    where the simplicial category $w\tilde{S}_\bullet\mathcal{C}$ is the \emph{modified $S_\bullet$-construction} of $\mathcal{C}$ (\cite[Definition 3.30]{Campbell}). There is an associated $K$-theory spectrum $K(\mathcal{C})$ which is defined by taking a (functorial) fibrant-cofibrant replacement in the stable model category of symmetric spectra of the spectrum whose $n$-th space is $|N_\bullet w\underbrace{\tilde{S}_\bullet\dots \tilde{S}_\bullet}_{n\ \mr{times}}\mathcal{C}|$. In particular, when $\mathcal{C}=\Var_k$, the associated \emph{$\K$-theory spectrum of varieties} is denoted $\K(\Var_k)$.
\end{defn}
For the rest of the paper, we restrict to the case $\mathcal{C}=\Var_k$. Similar to \cite{Lieb21}, we provide a framework to construct maps $\K(\mr{Var}_k)\to \K(\mathcal{A})$, where $\mathcal{A}$ is a nice $\infty$-category.
\begin{defn}\label{defn:Ktheoryspectrum}
    Let $\mathcal{C}$ be a pointed $\infty$-category admitting finite colimits. We define the \emph{$\K$-theory space of $\mathcal{C}$} to be
    \[
    \Omega^{\infty}\K(\mathcal{C}):=\Omega |S_\bullet^\infty(\mathcal{C})|,
    \]
    where the simplicial $\infty$-category $S_\bullet^\infty(\mathcal{C})$ is the $\infty$-categorical analogue of the $S_\bullet$ construction outlined in \cite[\S7.1]{BGT}. Again we have the associated spectrum $ \K(\mathcal{C})$.
\end{defn}

For the rest of the paper, we restrict to the case when $\mathcal{C}=\mathcal{A}^\omega$ is the $\infty$-category of compact objects of a presentable stable $\infty$-category.

The following definition is adapted from \cite[Section 6]{Kuijper}.
\begin{defn}
    Let $\mr{Span}_k^{\mr{sc}}$ denote the \emph{category of spans of scissors congruences of varieties}, where the objects are varieties over $k$, and a morphism $X\dashrightarrow Y$ is a span \[X\xleftarrow{\circ } U\xrightarrow{f} Y,\]
    where the left morphism is an open immersion and the right morphism is a closed embedding. Compositions are given by the usual pullback.
    Similarly, let $\mr{Span}_k$ denote the category of spans where the right leg of the span is allowed to be any proper morphism.

    Both $\mr{Span}_k^{\mr{sc}}$ and $\mr{Span}_k$ are 1-categories, since any isomorphism between pullbacks along open immersions is unique.
\end{defn}

The following definition is adapted from \cite[Definition 5.2]{Campbell}- one can make a more general definition, but we only require this level of generality for this paper.
    \begin{defn}\label{def:Wexact}
		Let $\mathcal{A}$ be a presentable stable $\infty$-category, and let $\mathcal{A}^\omega$ denote the subcategory of compact objects. An \emph{exact $\infty$-functor of spans} is an $\infty$-functor 
		\[
		F\colon \mr{Span}_k^{\mr{sc}}\longrightarrow\mathcal{A}^\omega
		\]
        where we identify $\mr{Span}_k^{\mr{sc}}$ with its nerve, such that
        \begin{itemize}
            \item $F$ sends $X\setminus U\dashrightarrow X\dashrightarrow U$ to a fibre sequence in $\mathcal{A}$ for any open $U\subset X$, where the left morphism is $X\setminus U\leftarrow X\setminus U\to X$ and the right morphism is $X\leftarrow U\to U$.
            \item $F(\emptyset)$ is a zero object of $\mathcal{A}$.
        \end{itemize}
\end{defn}

We introduce $\infty$-functors of spans for the following theorem, the proof of which is adapted from \cite[Proposition 5.3]{Campbell}.
\begin{thm}\label{thm:functorofspans}
    An exact $\infty$-functor of spans $F\colon\mr{Span}_k^{\mr{sc}}\to\mathcal{A}^\omega$ as in Definition $\ref{def:Wexact}$ defines a map of $\K$-theory spectra $\K(\Var_k) \to \K(\mathcal{A}^\omega)$ such that $[X] \in \K_0(\Var)$ gets mapped to the class $[F(X)]\in \K_0(\mathcal{A}^\omega)$.
\end{thm}
\begin{proof}
    We prove this theorem by showing that such an $\infty$-functor induces a morphism of simplicial categories \[N_\bullet w\tilde{S}_\bullet(\Var_k)\to N_\bullet S_\bullet^{\infty}(\mathcal{A}^\omega)_{\mr{iso}},\] where the left-hand side is the modified $S_\bullet$-construction (\cite[Definition 3.29]{Campbell}), and the right-hand side is the $S_\bullet$-construction of a pointed $\infty$-category with finite colimits (\cite[\S 7.1]{BGT}). Objects of $\tilde{S}_n(\mr{Var}_k)$ are functors $X\colon[n]^{\mr{op}}\times [n]\to \mr{Var}_k$ where the map $X_{i,j}\to X_{i,k}$ is a closed embedding for all $i$ and $j>k$, the map $X_{i,j}\to X_{\ell,j}$ is an open immersion for $i<\ell$, and the diagram
    \[X_{j-k,k}\to X_{i,k}\leftarrow X_{i,j}
    \] 
    is a subtraction sequence in $\mr{Var}_k$. This means that $X$ also furnishes a functor
    \[
    X\colon\mr{Ar}[n]\to \mr{Span}^{\mr{sc}}_k
    \]
    since open immersions $X_{i,j}\to X_{\ell,j}$ define canonical morphisms $ X_{\ell,j}\dashrightarrow X_{i,j}$ in $\mr{Span}_k^{\mr{sc}}$. Post-composing with the $\infty$-functor  $F\colon N\mr{Span}^{\mr{sc}}\to\mathcal{A}^\omega$ where $N$ denotes the nerve, we
    obtain another $\infty$-functor
    \[
    F\circ NX \colon N\mr{Ar}[n]\to \mathcal{A}^\omega.
    \]
    Since $X_{i,i}=\emptyset$, $(F\circ NX)_{0,0}$ is a zero object of $\mathcal{A}^\omega$, and since $F$ sends localisation sequences to fibre sequences (and therefore co-cartesian sequences in $\mathcal{A}$, since $\mathcal{A}$ is stable), we see that $F\circ NX\in\mr{Gap}([n],\mathcal{A}^\omega)$. This means we obtain an $\infty$-functor
    \[
    S_nF \colon N\tilde{S}_n(\mr{Var}_k)\to \mr{Gap}([n],\mathcal{A}^\omega).
    \]
    Moreover $S_nF$ commutes with the face and degeneracy maps on both sides, so can be promoted to the desired morphism of simplicial categories. This proof generalises to give a map between the iterations of the $S_\bullet$ constructions of both categories, which then induces the desired map of spectra $\K(F)\colon \K(\mr{Var}_k)\to \K(\mathcal{A}^\omega)$.
\end{proof}
\subsection{Relations in $\K_1$ Between Automorphisms}
Our work studies maps of spectra out of $\K(\Var_k)$ constructed using exact $\infty$-functors of spans. In particular, given such a map $\K(\Var_k)\to \K(\mathcal{A}^\omega)$, we focus on the induced map
\begin{equation}\label{eqn:K_1map}
\K_1(\Var_k)\to \K_1(\mathcal{A}^\omega).
\end{equation}
In order to do so, we consider presentations of these groups. Starting with the left-hand side of (\ref{eqn:K_1map}), there are two different presentations of $\K_1(\Var_k)$. One presentation comes from \cite[Theorem B]{Zak17-k1}, where the relations were shown to be a complete set of relations in \cite[Proposition 4.4]{KLMMS24}, and is obtained by viewing $\K(\Var_k)$ as the $\K$-theory spectrum of an assembler. However, we work with the presentation of $\K_1(\mr{Var}_k)$ given by Ng \cite[Theorem B]{Ng} obtained by viewing $\K(\mr{Var}_k)$ as the $\K$-theory spectrum of a pCGW-category, which subsumes the SW-category framework. By \cite[Theorem 1.3]{ZakhCampbell} these presentations are equivalent. We use Ng's presentation since it is easier to study the image of these generators in $\K_1(\mathcal{A}^\omega)$ sent by a map of the form (\ref{eqn:K_1map}). 
\begin{thm}[{\cite[Theorem B]{Ng}}]\label{thm:K1Presentation}
    The abelian group $\K_1(\mr{Var}_k)$ is generated by double exact squares, which are pairs
    \[
    \left(
    \begin{tikzcd}[sep=small]
    \emptyset\arrow[r]\arrow[d]&U\arrow[d,"g_1"]\\
    Z\arrow[r,"f_1"]&X
        \end{tikzcd}; \begin{tikzcd}[sep=small]
    \emptyset\arrow[r]\arrow[d]&U\arrow[d,"g_1'"]\\
    Z\arrow[r,"f_1'"]&X
        \end{tikzcd}
  \right)
    \]
    where $U,Z,X$ are varieties, the horizontal arrows are closed immersions, the vertical arrows are open immersions and $\mr{im}(f_1)$ (resp. $\mr{im}(f_1')$) is the complement of $\mr{im}(g_1)$ (resp. $\mr{im}(g_1')$). These generators satisfy three relations.
    \begin{enumerate}
        \item 
        \[
        \left(
    \begin{tikzcd}[sep=small]
    \emptyset\arrow[r]\arrow[d]&X\arrow[d,"\mr{id}_X"]\\
    \emptyset\arrow[r]&X
        \end{tikzcd};\ \begin{tikzcd}[sep=small]
    \emptyset\arrow[r]\arrow[d]&X\arrow[d,"\mr{id}_X"]\\
    \emptyset\arrow[r]&X
        \end{tikzcd}
  \right)=0,
        \]
    \item \[
        \left(
    \begin{tikzcd}[sep=small]
    \emptyset\arrow[r]\arrow[d]&\emptyset \arrow[d]\\
    X\arrow[r,"\mr{id}_X"]&X
        \end{tikzcd};\ \begin{tikzcd}[sep=small]
    \emptyset\arrow[r]\arrow[d]& \arrow[d]\emptyset\\
    X\arrow[r,"\mr{id}_X"]&X
        \end{tikzcd}
  \right)=0,
        \]
    \item Given compositions of closed embeddings $f_2\colon Z\xrightarrow{f_0}X\xrightarrow{f_1}V$ and $g_2\colon Z\xrightarrow{g_0}X\xrightarrow{g_1}V$ satisfying the conditions of \cite[Proposition 4.1]{Ng}, we have that the sum
\[
        \left(
    \begin{tikzcd}[sep=small]
    \emptyset\arrow[r]\arrow[d]&X\setminus Z \arrow[d,"g_0"]\\
    Z\arrow[r,"f_0"]&X
        \end{tikzcd};\ \begin{tikzcd}[sep=small]
    \emptyset\arrow[r]\arrow[d]&X\setminus Z \arrow[d,"g_0'"]\\
    Z\arrow[r,"f_0'"]&X
        \end{tikzcd}
  \right)+\left(
     \begin{tikzcd}[sep=small]
    \emptyset\arrow[r]\arrow[d]&V\setminus X
    \arrow[d,"g_1"]\\
    X\arrow[r,"f_1"]&V
        \end{tikzcd};\ \begin{tikzcd}[sep=small]
     \emptyset\arrow[r]\arrow[d]&V\setminus X
    \arrow[d,"g_1'"]\\
    X\arrow[r,"f_1'"]&V
        \end{tikzcd}
  \right)
  \]
  is equal to the double exact square 
  \[\left(
    \begin{tikzcd}[sep=small]
    \emptyset\arrow[r]\arrow[d]& V\setminus Z\arrow[d,"g_2"]\\
    Z\arrow[r,"f_2"]&V
        \end{tikzcd};\ \begin{tikzcd}[sep=small]
    \emptyset\arrow[r]\arrow[d]& V\setminus Z\arrow[d,"g_2'"]\\
    Z\arrow[r,"f_2'"]&V
        \end{tikzcd}
  \right).
        \]
    
    \end{enumerate}
\end{thm}
\begin{rem}
    Theorem $\ref{thm:K1Presentation}$ gives a complete set of generators and relations for $\K_1(\Var_k)$. For most of this paper, we only work with certain distinguished elements of $\K_1(\Var_k)$, which have a simpler form and arise by looking at automorphisms of varieties.
\end{rem}

\begin{defn}\label{def:SWK_1Aut}
    Let $X$ be a variety over $k$, and let $f\colon X\to X$ be an automorphism of $X$. Define the element $[X,f]$ of $\K_1(\Var_k)$ to be the doubly exact square
    \[
    [X,f]\colon =\left(
    \begin{tikzcd}[sep=small]
    \emptyset\arrow[r]\arrow[d]&X\arrow[d,"\mr{id}_X"]\\
    \emptyset\arrow[r]&X
        \end{tikzcd};\ \begin{tikzcd}[sep=tiny]
    \emptyset\arrow[r]\arrow[d]&X\arrow[d,"f"]\\
    \emptyset\arrow[r]&X
        \end{tikzcd}
  \right).
    \]
\end{defn}
The following results can also be proved using the assembler formalism from \cite{Zak17-assemblers}.
\begin{lem}\cite[Corollary 4.11(ii)]{Ng}\label{lemma:AutEmbedding}
    Let $X$ be a variety over $k$. Then the assignment $f \mapsto [X,f]$ defines a group homomorphism $\mr{det}_X\colon \Aut(X) \to \K_1(\Var_k)$, which factors through $\Aut(X)^{ab}$.
\end{lem}
\begin{rem}
    The morphism $\det_X$ is written $\det_X$ since it behaves like a universal determinant, see the discussion in the introduction of \cite{Zakharevich2022}. This is the higher homotopical analogue of how the assignment of a variety $X$ to its class in $\K_0(\Var_k)$ behaves like a universal motivic measure (see \cite{Bit04} for a discussion).
\end{rem}

\begin{cor}\label{cor:equalautos}
    Given a diagram of isomorphisms in $\Var_k$
    \[
    \begin{tikzcd}        X\arrow[r,"g"]\arrow[d,"f_X"]&Y\arrow[d,"f_Y"]\\
        X\arrow[r,"g"]&Y,
    \end{tikzcd}
    \]
    we have the equality $[X,f_X]=[Y,f_Y]$ in $\K_1(\Var_k)$.
\end{cor}
\begin{proof}
    Consider the optimal $3 \times 3$ diagram as in \cite[Definition 4.6]{Ng}
    \[
   \left( \begin{tikzcd}
        X\arrow[r,"g"]\arrow[d,"\mr{id}_X"]&Y\arrow[d,"\mr{id}_Y"]&\emptyset\arrow[d]\arrow[l]\\
        X\arrow[r,"g"]&Y& \emptyset\arrow[l]\\        \emptyset\arrow[u]\arrow[r,equal]&\emptyset\arrow[u]&\emptyset\arrow[l,equal]\arrow[u]
    \end{tikzcd}; \begin{tikzcd}
        X\arrow[r,"g"]\arrow[d,"f_X"]&Y\arrow[d,"f_Y"]&\emptyset\arrow[d]\arrow[l]\\
        X\arrow[r,"g"]&Y& \emptyset\arrow[l]\\        \emptyset\arrow[u]\arrow[r,equal]&\emptyset\arrow[u]&\emptyset\arrow[l,equal]\arrow[u]
    \end{tikzcd}\right).
    \]
By Theorem 4.9 of \cite{Ng}, relation (2) gives us equality between the sum
\[
\left(
    \begin{tikzcd}[sep=small]
    \emptyset\arrow[r]\arrow[d]&\emptyset\arrow[d]\\
    X\arrow[r,"g"]&Y
        \end{tikzcd}; \begin{tikzcd}[sep=small]
    \emptyset\arrow[r]\arrow[d]&\emptyset\arrow[d]\\
    X\arrow[r,"g"]&Y
        \end{tikzcd}
  \right)-\left(
    \begin{tikzcd}[sep=small]
    \emptyset\arrow[r]\arrow[d]&\emptyset\arrow[d]\\
    X\arrow[r,"g"]&Y
        \end{tikzcd}; \begin{tikzcd}[sep=small]
    \emptyset\arrow[r]\arrow[d]&\emptyset\arrow[d]\\
    X\arrow[r,"g"]&Y
        \end{tikzcd}
  \right)+\left(
    \begin{tikzcd}[sep=small]
    \emptyset\arrow[r]\arrow[d]&\emptyset\arrow[d]\\
    \emptyset\arrow[r]&\emptyset
        \end{tikzcd}; \begin{tikzcd}[sep=small]
    \emptyset\arrow[r]\arrow[d]&\emptyset\arrow[d]\\
    \emptyset\arrow[r]&\emptyset
        \end{tikzcd}
  \right)
\]
where the first two double exact squares are equivalent, and the  sum
\[
    [X,f_X]-[Y,f_Y]+\left(\begin{tikzcd}[sep=small]
    \emptyset\arrow[r]\arrow[d]&\emptyset\arrow[d]\\
    \emptyset\arrow[r]&\emptyset
        \end{tikzcd}; \begin{tikzcd}[sep=small]
    \emptyset\arrow[r]\arrow[d]&\emptyset\arrow[d]\\
    \emptyset\arrow[r]&\emptyset
        \end{tikzcd}\right).
    \]
Relation (1) tells us that $\left(
    \begin{tikzcd}[sep=small]
    \emptyset\arrow[r]\arrow[d]&\emptyset\arrow[d]\\
    \emptyset\arrow[r]&\emptyset
        \end{tikzcd}; \begin{tikzcd}[sep=small]
    \emptyset\arrow[r]\arrow[d]&\emptyset\arrow[d]\\
    \emptyset\arrow[r]&\emptyset
        \end{tikzcd}
  \right)=0$, giving us that $[X,f_X]=[Y,f_Y]$.
\end{proof}

An important relation in $\K_1(\Var_k)$ is the following cut-and-paste relation. Let $X$ be a variety over $k$, $f\colon X\to X$ be an automorphism, and $Z\subset X$ be a closed subvariety such that $f$ restricts to an automorphism $f|_{Z}\colon Z\to Z$.
\begin{prop}\label{prop:Autocutanpaste}
     In $\K_1(\Var_k)$, we have
    \[
    [X,f]=[X\setminus Z,f|_{X\setminus Z}]+[ Z,f|_{Z}].
    \]
\end{prop}
\begin{proof}
    This proof is identical to Corollary \ref{cor:equalautos}, replacing the optimal $3 \times 3$ diagram in the proof with the optimal $3 \times 3$ diagram
    \[
   \left( \begin{tikzcd}
        Z\arrow[r]\arrow[d,"\mr{id}_Z"]&X\arrow[d,"\mr{id}_X"]&X\setminus Z\arrow[d,"\mr{id}_{X\setminus Z}"]\arrow[l]\\
        Z\arrow[r]&X& X\setminus Z\arrow[l]\\        \emptyset\arrow[u]\arrow[r,equal]&\emptyset\arrow[u]&\emptyset\arrow[l,equal]\arrow[u]
    \end{tikzcd}; \begin{tikzcd}
        Z\arrow[r]\arrow[d,"f|_Z"]&X\arrow[d,"f"]&X\setminus Z\arrow[l]\arrow[d,"f|_{X\setminus Z}"]\\
        Z\arrow[r]&X& X\setminus Z\arrow[l]\\        \emptyset\arrow[r,equal]\arrow[u]&\emptyset\arrow[u]&\emptyset\arrow[l,equal]\arrow[u]
    \end{tikzcd}\right). \qedhere
    \]
\end{proof}
Now let $X$ be a variety and $Z$ a closed subvariety. Let $f\colon X \to X$ be an automorphism such that $f(Z) = Z$. Consider the blow-up $\mathrm{Bl}_Z(X)$ and write $\pi\colon \mathrm{Bl}_ZX \to X$ and let $E_Z = (\pi^{-1}(Z))_{red}$ be the (reduced) exceptional divisor. Since $f(Z) = Z$, $f$ lifts to an automorphism $\tilde{f}$ of $\mathrm{Bl}_Z(X)$ such that $\tilde{f}(E_Z) = E_Z$.

\begin{cor}\label{cor:BlowUpRelation}
There is an equality in $\K_1(\Var_k)$
\[
[X,f] - [Z, f|_Z] = [\mathrm{Bl}_Z(X), \tilde{f}] - [E_Z, \tilde{f}|_{E_Z}].
\]
\end{cor}
\begin{proof}
 Proposition \ref{prop:Autocutanpaste} implies
\begin{align*}
[X,f]&=[Z,f_{Z}]+[X\setminus Z,f|_{X\setminus Z}],\\
[\mr{Bl}_{Z}(X),\tilde{f}]&=[E_{Z},\tilde{f}|_{E_Z}]+[\mr{Bl}_{Z}(X)\setminus E_{Z},\tilde{f}_{\mr{Bl}_{Z}(X)\setminus E_{Z}}].
\end{align*}
Furthermore, we have a commutative diagram of isomorphisms
\[
\begin{tikzcd}
    \mr{Bl}_{Z}(X)\setminus E_{Z}\arrow[r,"\pi"]\arrow[d,"\tilde{f}"] &X\setminus Z\arrow[d,"f"]\\
    \mr{Bl}_{Z}(X)\setminus E_{Z}\arrow[r,"\pi"] &X\setminus Z.
\end{tikzcd}
\]
By Corollary \ref{cor:equalautos}, $[\mr{Bl}_{Z}(X)\setminus E_{Z},\tilde{f}_{\mr{Bl}_{Z}(X)\setminus E_{Z}}]=[X\setminus Z,f_{X\setminus Z}]$, which gives the result.
\end{proof}
\begin{rem}
    In \cite[Theorem 3.1]{Bit04}, it is shown that, in characteristic $0$, $\K_0(\Var_k)$ is generated by the classes of smooth projective varieties, subject only to a blow-up relation. It remains to be seen, and is beyond the scope of this paper, whether an analogous presentation for $\K_1(\Var_k)$ exists in characteristic $0$. Corollary $\ref{cor:BlowUpRelation}$ demonstrates the first relation that would be needed for such a presentation, but it is likely one would need other relations that are not determined by automorphisms.
\end{rem}
\begin{rem}\label{rem:equivariantrov}
    One curiosity of the above relations is the following. Following \cite[\S7]{Bit04}, let $\K^{\Z/n\Z}_0(\Var_k)$ denote the Grothendieck ring of $\Z/n\Z$-varieties. There is a homomorphism of the underlying abelian groups
    \[
    \K^{\Z/n\Z}_0(\Var_k)\to \K_1(\Var_k)
    \]
    constructed as follows. Let $X$ be a variety over $k$ and let $f\colon \Z/n\Z \to \Aut(X)$ be the homomorphism defining a $\Z/n\Z$ action. We obtain the above homomorphism by sending $(X,f)$ to $[X, f(1)]$.
\end{rem}

We now discuss $\K_1(\mathcal{C})$, where $\mathcal{C}$ is a pointed stable $\infty$-category with finite limits. One can associate to $\mathcal{C}$ the $\K$-theory spectrum of a Waldhausen category (which is equivalent to the $\infty$-categorical $S_\bullet$ construction in \cite{BGT}). By \cite[Theorem 2.2]{MT08}, the generators and relations for $\K_1(\mathcal{C})$ are known. We introduce the simplest of these generators, namely the ones corresponding to homotopy classes of automorphisms. 
\begin{defn}\label{def:WK1Aut}
    Let $\mathcal{A}$ be a presentable stable $\infty$-category, and $X$ be an object of $\mathcal{A}^\omega $. Let $\mathcal{A}_{\mr{iso}}^\omega$ denote its core. By definition of the $S_\bullet$ construction we have that $\mathcal{A}^\omega$ and $\mr{Gap}([1],\mathcal{A}^\omega)$ are equivalent, which produces a map of topological spaces
    \[
    \Sigma|\mathcal{A}_{\mr{iso}}^\omega|\to |S_\bullet^{\infty}(\mathcal{A}^\omega)_{\mr{iso}}|.
    \]
  The adjoint map is $\iota \colon |\mathcal{A}^\omega_{\mr{iso}}|\to \Omega|S_\bullet(\mathcal{A}^\omega)_{\mr{iso}}|$. By definition of the $S_\bullet$ construction, $\pi_1(|\mathcal{A}^\omega_{\mr{iso}}|, [X]) \simeq \pi_0(\mr{Aut}(X))$, and $\pi_1(\Omega|S_\bullet(\mathcal{A}^\omega)_{\mr{iso}}|, [X]) \simeq \pi_1 (\Omega^\infty\K(\mathcal{A}^\omega),[X]) \simeq \K_1(\mathcal{A}^\omega)$, where $\mr{Aut}(X)$ is the core of the mapping space (spectrum) of $X$.  Therefore  $\pi_1(\iota)$ defines a map $\mr{det}_X\colon\pi_0(\mr{Aut}(X))\to \K_1(\mc{A}^\omega)$. We write the image of $\phi\in \mr{Aut}(X)$ as $[X,\phi]\in \K_1(\mc{A}^\omega)$.
\end{defn}
We use ``det'' both when speaking about $\K_1(\mr{Var})$ and $\K_1(\mr{A}^\omega)$, which are defined differently. The following proposition shows that this abuse of notation is justified. 
\begin{prop}\label{prop:maponk1}
    Let $F\colon \mr{Span}_k^{\mr{sc}}\to\mc{A}^\omega$ be an exact $\infty$-functor of spans and write the induced map on $\pi_1$ as $\K_1(F)\colon \K_1(\mc{C})\to\K_1(\mc{A})$. Consider $X\in\mr{Var}_k$ and an automorphism $f\colon X\to X$. Then
    \[
    \K_1(F)([X,f])=[F(X),F(f)].
    \]
\end{prop}
\begin{proof}
    Let $w\tilde{S}_\bullet(\mc{C})$ denote the simplicial category defined by the modified $\tilde{S}_\bullet$ construction given in \cite{Campbell}, and let $s\tilde{S}_\bullet(\mathcal{C})$ be the associated simplicial set. The same argument as in Definition \ref{def:WK1Aut} shows that we have a map of spaces 
    \[
    \Sigma |\mathcal{C}^\times|\to |w\tilde{S}_\bullet(\mc{C})|.
    \] 
    An automorphism $f\colon X\to X$ defines a 1-simplex in $N_\bullet\mathcal{C}^\times$ and the image of the induced map $\Delta^1\times\Delta^1\to N_\bullet w\tilde{S}_\bullet(\mc{C})$ gives the corresponding element $\mr{det}_X(f)\in \pi_2\Big(|w\tilde{S}_\bullet(\mc{C})|\Big)$. We claim this element is the same element as defined by $[X,f]$ from Ng's presentation in Definition \ref{def:SWK_1Aut}, adapting a strategy in \cite{MT1type}. Once this is proven, the rest follows by an identical argument to \cite[Proposition 2.20]{CWZ19}. 
    
    As $\pi_1|\Omega s\tilde{S}_\bullet\mathcal{C}|=\pi_2|s\tilde{S}_\bullet\mathcal{C}|$, Ng's double exact squares in Definition \ref{def:SWK_1Aut} build a $2$-sphere in $s\mathcal{S}_\bullet\mathcal{C}$ by gluing the 2 simplices $(0\to X\xrightarrow{f} X)$ and $(0\to X\xrightarrow{\mr{id}}X)$ along their common boundary.
    Viewing $s\tilde{S}_\bullet\mathcal{C}$ as a bisimplicial set with the `Nerve' degree constant, we have an inclusion of bisimplicial sets $s\tilde{S}_\bullet\mathcal{C}\hookrightarrow N_\bullet(w\tilde{S}_\bullet(\mathcal{C}))$ which turns Ng's 2-simplices into (0,2) simplices. Consider the (1,2)-simplex in $N_\bullet(w\tilde{S}_\bullet(\mathcal{C}))$
    \[
    \begin{tikzcd}
        X\arrow[r,"\mr{id}"]\arrow[d,"\mr{id}"]&X\arrow[d,"f"]&\emptyset\arrow[l]\arrow[d,equal]\\
        X\arrow[r,"f"]&X&\emptyset\arrow[l].
    \end{tikzcd}
    \]
    This (1,2)-simplex is isomorphic to a copy of $\Delta^1\times\Delta^2$. The induced map on the boundary $\partial(\Delta^1\times\Delta^2)\to N_\bullet(w\tilde{S}_\bullet(\mathcal{C}))$ is homotopically trivial (this mirrors relation (7) in Definition 1.4 of \cite{MT1type}). The boundary of this simplex is the gluing of the two simplices of the rows along a common boundary given by Ng's 2-spheres, along with the three 2-simplices: $X\xrightarrow{\mr{id}}X, X\xrightarrow{f}X, \emptyset\to\emptyset$ (these are  2-simplices because the objects in the morphism are 1-simplices). The only non-degenerate 2-simplex is $X\xrightarrow{f}X$, which is precisely the image of $\Delta^1\times\Delta^1\to N_\bullet w\tilde{S}_\bullet(\mathcal{C})$ defined by $f$. So, we have the relation in $\pi_2|w\tilde{S}_\bullet(\mathcal{C})|$
    \[
    \mr{det}_X(f)\pm [X,f]=0,
    \]
    where the sign will be the same for all $X$ and $f$ as it depends on the conventions for the boundary maps. We choose the sign that gives us $\mr{det}_X(f)=[X,f]$, as required.
\end{proof}

    \subsection{The stable motivic homotopy category}
This paper studies a spectrum level lift of the compactly supported $\mathbb{A}^1$-Euler characteristic, a motivic measure constructed using $\mb{A}^1$-homotopy theory. We will construct this spectrum level lift by using the previous section to build an exact $\infty$-functor of spans landing in a suitable subcategory of the stable motivic homotopy category, $\SH(k)$. In this subsection, we recall some basic definitions regarding $\SH(k)$ and the $\mathbb{A}^1$-Euler characteristic.

\begin{defn}
    Write $\SH(k)$ to mean the stable $\infty$-category of motivic spectra over $\Spec(k)$ as constructed in \cite{Voevodsky}. Write $\1_k$ to mean the sphere spectrum in $\SH(k)$ (which is sometimes written as $\mathbb{S}_k$ in other literature). For $X$ a variety over $k$, let $\Sigma^{\infty} X_+$ denote its infinite suspension spectrum in $\SH(k)$.
\end{defn}
\begin{defn}
     Let $e$ denote the exponential characteristic of  $k$. Let $\mr{Sp}$ be the stable $\infty$-category of spectra and $\mr{Sp}[1/e]$ the stable $\infty$-category of $e$-inverted spectra. Write \[\SH(k)[1/e]:=\SH(k)\otimes \mr{Sp}[1/e]\]
     where the right-hand side denotes the Lurie tensor product. Write $\1_k[1/e]$ to mean the $p$-inverted sphere spectrum in $\SH(k)[1/e]$.
\end{defn}
If $\mr{char}(k)=0$ then $\SH(k)[1/e]\simeq \SH(k)$. The \emph{smash product} $\wedge$ furnishes $\SH(k)$ with the structure of a stable symmetric monoidal $\infty$-category with unit $\1_k$ (see \cite[Remark 5.1.15]{Morel2004}). Since $\mr{Sp}[1/e]$ is given by a smashing localization (since $\mr{Sp}[1/e]$ is an idempotent $\Einfty$ algebra) of $\SH(k)$, it is a stable symmetric monoidal $\infty$-category, with unit $\1_k[1/e]$.

\begin{defn}
    Let $\mathcal{A}$ be a presentable, stable symmetric monoidal $\infty$-category. An object $A\in \mathcal{A}$ is \emph{strongly dualisable} if for every $B\in \mathcal{A}$, the canonical map
    \[
    \mr{Hom}(A,\1_\mathcal{A})\otimes B\to \mr{Hom}(A,B)
    \]
    is an equivalence. In the category $\mathcal{A}= \SH(k)$ \cite[Theorem 6.9]{Hoy17} implies that if a $k$-scheme $X$ is smooth and projective, then $\Sigma^{\infty} X_+ \in \SH(k)$ is strongly dualisable and its dual is the Thom spectrum of its stable normal bundle (see \cite[Theorem 6.9]{Hoy17} and \cite[Construction 2.12]{LR20} for a full exposition). For $C$ a strongly dualisable object of $\SH(k)$, write $C^\vee$ for its dual. 
\end{defn}

    Many constructions in this paper come from the endomorphism spectrum $\End_{\SH(k)}(\1_k)$. We view $\End_{\SH(k)}(\1_k)$ as an $S^1$-spectrum (by viewing it as the $\Sigma^{\infty}\mathbb{G}_m$-delooping of the bi-spectrum $\1_k$) and note that $\pi_1(\End_{\SH(k)}(\1_k)) = [ \Sigma \1_k, \1_k] = \pi_{1,0}(\1_k)$ (see the discussion in \cite[Corollary 6.7]{Ron25} for more details). In \cite[Theorem 6.4.1]{Morel2004}, Morel shows that $\pi_0(\End_{\SH(k)}(\1_k))$ is isomorphic to the Grothendieck--Witt ring of $k$ (Definition $\ref{defn:GWk}$). We mostly avoid using this isomorphism, since it can be misleading when talking about higher homotopy groups. See the discussion at the beginning of Section \ref{section:Hermitian} for further details.

\begin{defn}\label{defn:categoricalchi0}
Following \cite{Hoy14}, we define the \emph{$\mb{A}^1$-Euler characteristic} of a smooth projective variety $X$, written $\chi^{\mb{A}^1}(X) \in \pi_0\End_{\SH(k)}(\1_k)$, to be the categorical trace of the identity of $\Sigma^\infty X_+$, in the sense of \cite[Theorem 1.3]{DoldPuppe1999}. This is the class of the trace
\[\1_k\to\Sigma^\infty X_+\wedge(\Sigma^\infty X_+)^\vee\xrightarrow{\tau}(\Sigma^\infty X_+)^\vee\wedge\Sigma^\infty X_+\to\1_k.\]

For any variety $p \colon X \to k$ over a perfect field $k$ of exponential characteristic $e$, define the \emph{compactly supported $\mb{A}^1$-Euler characteristic} $\chi_0$ by replacing $\Sigma^{\infty} X_+$ in the above definition with the image of the \emph{compactly supported motive} $h^c(X) := p_! \1_X$ in $\SH(k)[1/e]$. This furnishes a ring homomorphism
\[ \chi_0 \colon \K_0(\Var_k)  \to \pi_0\End_{\SH(k)[1/e]}(\1_k[1/e])\]

(see \cite[\S~5.1]{LPLS24}, \cite[Definition~2.7 and Remark~2.10]{Azo}). The above ring homommorphism satisfies that for smooth projective varieties $X$, $\chi_0([X])$ is the same as $\chi^{\mb{A}^1}(X)$ after inverting $e$. For ease of notation, we write $\chi_0(X)$ to mean $\chi_0([X])$. 
\end{defn}
    For ease of notation, we will not use the ``compactly supported'' prefix, though readers should note that this is an analogue of the compactly supported Euler characteristic, not the standard one. The above definition of $\chi^{\mb{A}^1}$ is the original definition given in \cite{Hoy14}. It is often also referred to as the ``quadratic Euler characteristic'', and occasionally as the ``motivic Euler characteristic''. In Section \ref{section:Hermitian}, we will discuss a more geometric, and easily computable, interpretation of $\chi_0$ coming from \cite[Theorem 1.3]{LR20}.
   
The trace in Definition \ref{defn:categoricalchi0} was first lifted to a map of spectra 
\begin{equation}\label{map:derived}\K(\Var_k)\to \mr{End}_{\SH(k)}(\1_k)
\end{equation}
by the second author in Theorem 1.1 of \cite{Nan24}, using model categories. His work also relies on work of R\"ondigs \cite[Theorem 6.6]{Ron25}, which in turn relies on a construction by Hoyois--Scherotzke--Sibilla \cite[Theorem 6.5, Remark 6.6]{HSS17} when the characteristic of $k$ is $0$. Little is known about the map $\ref{map:derived}$ aside from its induced map on $\pi_0$. This paper is concerned with generalising the spectrum level construciton to positive characteristic and analysing the map induced by $\chi$ on the first homotopy group
\[
\chi_1\colon \K_1(\Var_k) \to \pi_1(\End_{\SH(k)[1/e]}(\1_k[1/e])). 
\]

\section{The derived $\mathbb{A}^1$-Euler characteristic, six-functors and traces}\label{section:categorical}
Recall that $k$ is a perfect field of exponential characteristic $e$. In this section we give a construction of the derived $\mathbb{A}^1$-Euler characteristic
	\[
	\chi \colon \K(\Var_k)\longrightarrow
	\End_{\SH(k)[1/e]}(\1_k[1/e]).
	\]
   Our construction is more general and more streamlined than in \cite{Nan24}, and allows us to analyse the induced map on homotopy groups more easily. To construct our $\chi$, we first construct a map of spectra
	\[
	\mot\colon	\K(\Var_k)\longrightarrow	\K(\SH(k)[1/e]^\omega) ,
	\]
	 by using a descent principle for compactly supported extension, outlined in \cite{Kuijper}. On $\K_0$, the map $\mot$ sends the class of a variety $p \colon X \to k$ to $p_* p^! \1_k$, localised by inverting the characteristic. Inverting the characteristic ensures that the target is a rigid stable $\infty$-category, which allows us to apply the higher categorical trace of \cite[Theorem 6.5, Remark 6.6]{HSS17} (see also \cite[Theorem 6.6]{Ron25} for the application of the methods of Hoyois--Scherotzke--Sibilla in this setting). We then discuss the induced map on $\pi_1$, which is the central invariant of this paper
		\[\chi_1:=\pi_1\chi\colon \K_1(\Var_k)\to\pi_1\End_{\SH(k)[1/e]}(\1_k[1/e]).\]
        For $X$ a variety and $f\colon X \to X$ an automorphism, we construct $\chi_1([X,f])$ as an explicit element in the motivic homotopy group $\pi_{1,0} (\1_k[1/e])$ by describing it in terms of the $2$-categorical trace, following \cite[pp.~108--109]{HSS17}.

	\subsection{The compactly supported motivic realisation map}

If $X$ is a smooth projective variety, then by \cite[Proposition C.12]{Hoy14}, we have $\Sigma^\infty X_+\in \SH(k)^\omega$. The following proposition relies on the cohomological descent theorem of \cite[Theorem 1.2]{Kuijper}.
\begin{prop}\label{prop:Mc}
 There is an $\infty$-functor
    \[
    h^c \colon \mr{Span}^{op}_k\to \SH(k)
    \]
    such that for any variety $p_X \colon X \to \Spec k$, we have $h^c(X) \simeq p_{X,!} \1_X \in \SH(k)$.
 For an isomorphism $f\colon X\xrightarrow{\sim} Y$, we have that $h^c(f)\colon h^c(Y)\to h^c(X)$ is equivalent to
\begin{equation} \label{eq:h^c(f)}
    p_{Y!}\1_Y \xrightarrow{p_{Y!}(\eta_f)}  p_{Y!}f_*f^*\1_Y \simeq p_{X!}\1_X,
\end{equation}
where $\eta_f$ is the unit of the adjunction $(f^*,f_*)$.
 Moreover, $h^c$ sends  $\emptyset$ to a zero object, and  a sequence of the form $X\setminus U\dashrightarrow X\dashrightarrow U$ to the fibre sequence $h^c(U) \to h^c(X) \to h^c(X \setminus U)$.
\end{prop}
\begin{proof}
    Let $\mr{Comp}_k$ denote the nerve of the full subcategory consisting of proper varieties in the 1-category $\Var_k$. In \cite[Corollary 6.25]{Hoy17} there is an $\infty$-functor defined on the nerve of the 1-category of schemes over $k$ which we can restrict to $\mr{Comp}_k$ and write it as
    \[
   h_0 \colon \mr{Comp}^{op}_k\to \SH(k), (p \colon X \to \Spec k ) \mapsto p_*p^*\1_k.
    \]
  In loc. cit. Hoyois proves that $h_0$ is a sheaf for the cdh topology, so it satisfies descent for abstract blowup squares. This is equivalent to the following. Suppose we have a cartesian square in $\mathrm{Comp}_k$
    \[\begin{tikzcd}
        X\arrow[r]\arrow[d]&Y\arrow[d,"p"]\\
        Z\arrow[r,"i"]&W
    \end{tikzcd}\]
    where $i$ is a closed immersion, and $p$ is a proper morphism such that $p$ induces an isomorphism $p^{-1}(W\setminus Z)\to W\setminus Z$. Then $h_0$ sends this to a cartesian, and hence cocartesian square in $\SH(k)$. In addition $h_0(\emptyset) \simeq 0$ by definition. As a result, \cite[Theorem 1.2]{Kuijper} tells us that $h_0$ extends to an $\infty$-functor (unique up to homotopy) which we write 
    \[    h^c\colon\mr{Span}^{op}_k\to \SH(k).
    \]
    By construction, this functor commutes with the embedding $\mr{Comp}_k\to \mr{Span}_k$ up to homotopy. Moreover, $h^c$ sends $X\setminus U\dashrightarrow X\dashrightarrow U$ to a fibre sequence in $\SH(k)$, and since $\Sigma^\infty(\emptyset)$ is a zero object in $\SH(k)$, we have that $h^c(\emptyset)$ is also a zero object.
    Let $p_X \colon X \to \Spec k$ be any variety, and choose a compactification $j \colon X \to \bar{X}$ (so that $\bar{X}$ is proper) with closed complement $i \colon Z \to \bar{X}$. We have the fibre sequence
    \[ h^c(X) \to h^c(\bar{X}) \to h^c(Z) .\]
    Since the maps $p_ {\bar{X}} \colon \bar{X} \to \Spec k$, $p_{Z} \colon Z \to \Spec k$ are proper, we have isomorphisms
    $h^c(\bar{X}) \simeq p_{\bar{X}, !} \1_{\bar{X}} $,   $h^c(Z) \simeq p_{Z, !} \1_Z $, and $h^c(X) \simeq \mathrm{fib}(p_{\bar{X}, !} \1_{\bar{X}}  \to p_{Z, !} \1_Z )$, where by definition of $h_0$ the last map is induced by the inclusion $i$. By the localisation property of the six functor formalism \cite[Theorem 6.18(4)]{Hoy17} we have the fibre sequence
    \[ p_{{X}, !} \1_{X} \to p_{\bar{X}, !} \1_{\bar{X}} \to p_{Z, !} \1_Z ,\]
    which gives us an isomorphism $h^c(X) \simeq p_{X, !} \1_X$.
    
    Now let $f\colon X\to Y$ be an isomorphism. Choose a compactification $j_Y\colon Y\hookrightarrow \bar{Y}$, and define $j_X=j_Y\circ f$, which gives $\bar{Y}$ also as a compactification of $X$. The identity map on $\bar{Y}$ extends $f$ and induces the identity on $Z = \bar{Y} \setminus Y$.
    By functoriality, $h^c(f)$ is the resulting map between the fibres of the two localisation sequences. After applying the identifications $h^c(X) \simeq p_{X, !} \1_X$, $h^c(Y) \simeq p_{Y, !} \1_Y$, the map $h^c(f)$ is obtained by applying proper pushforward along $q \colon \bar{Y} \to \Spec k$ to
 \begin{equation} \label{eq:h^c(f)proof}
    j_{Y!}\1_Y  \xrightarrow{j_{Y!}(\eta_f)}  j_{Y!}f_*f^*\1_Y  \simeq j_{X!}\1_X.
    \end{equation}
    Indeed, since $j_X=j_Y\circ f$, the counit for $(j_{X!},j_X^!)$ is the composite of $j_{Y!}(\varepsilon_f)$ with the counit for $j_Y$, where $\varepsilon_f\colon f_!f^!\1_Y\to\1_Y$ is the counit for $f$.
    Since $f$ is an isomorphism, the equivalences $f_! \simeq f_*$, and $f^! \simeq f^*$, identify $\varepsilon_f$ with the inverse of $\eta_f$. Consequently, the map \eqref{eq:h^c(f)proof} above is compatible with the maps to $\1_{\bar Y}$. Applying $q_!$ gives Equation~\ref{eq:h^c(f)}, since $q\circ j_X=p_X$ and $q\circ j_Y=p_Y$.
\end{proof}
\begin{prop}\label{prop:compdual}
    Every object of $(\SH(k)[1/e])^\omega$ is strongly dualisable. 
\end{prop}
\begin{proof}
    By \cite[Proposition 6.4]{Hoy17}, we have that $\SH(k)^\omega$ is generated by $T^{-n}\wedge \Sigma_+^\infty X$ where $X$ is a smooth affine $k$ scheme and $T\simeq S^1\wedge \mb{G}_m$. Since $T$ is invertible, it is strongly dualisable. Since the localisation $L:\SH(k)\to\SH(k)[1/e]$ produces an adjunction with right adjoint given by inclusion $i:\SH(k)[1/e]\to \SH(k)$, we see that $\SH(k)[1/e]^\omega$ is also generated by objects of the form $L(T^{-n}[1/e])\wedge L(\Sigma_+^\infty X)$.
    
    Strong dualisability can be checked in the homotopy category of $(\SH(k)[1/e])^\omega$. The object $L(T)$ remains invertible, so it is still strongly dualisable in $(\SH(k)[1/e])^\omega$. When $\mr{char}(k)=0$, strong dualisability of $\Sigma_+^\infty X$ is the statement of \cite[Theorem 4.11]{Ron25} and when $\mr{char}(k)=p$ it is the statement of \cite[Corollary B.2]{Riouappendix}. As the smash product preserves dualisability, this gives the result. 
\end{proof}
 Write $h^c_{1/e}\colon \mr{Span}^{\mr{op}}_k\to \SH(k)[1/e]$ to denote the composition of $h^c$ with the exact localisation functor $\SH(k)\to\SH(k)[1/e]$. 
 
\begin{prop}\label{prop:mot1}
   The essential image of $h^c_{1/e}$ is in $(\SH(k)[1/e])^\omega$. 
\end{prop}
    \begin{proof}
       Firstly, by the discussion after \cite[Definition 1.3.6]{BD17}, the object $p_! \1_X$ is constructible and compact in $\SH(k)$ (viewed as a triangulated category), which remains true if we view $\SH(k)$ as an $\infty$-category. The localisation functor $L:\SH(k)\to\SH(k)[1/e]$ preserves compact objects, so every object in the image of $h^c_{1/e}$ is compact, as required.
       \end{proof}

       In order to use Proposition \ref{prop:maponk1} to obtain a map on $K$-theory spectra, we need a functor with similar properties to those of the functor of Proposition~\ref{prop:Mc}, but with opposite covariance.  
       
       \begin{defn}
       Define the functor $\M^c := (h^c_{1/e})^{\vee} \colon \mr{Span}_k \to \SH(k)[1/e]$. This is well defined, since Proposition $\ref{prop:mot1}$ shows that $h^c_{1/e}$ has essential image landing in compact objects of $\SH(k)[1/e]$, which are dualisable by Proposition \ref{prop:compdual}.
       \end{defn}
       \begin{lem}\label{lem:McProperties}
           If $p\colon  X \to \Spec(k)$ and $f\colon X \to X$ is an automorphism, then we have an isomorphism $\M^c(X) \simeq p_*p^! \1_k[1/e]$, and $\M^c(f)$ is equivalent to  \begin{equation} \label{eq:M^c(f)}
    p_*p^!(\1_k[1/e])  \simeq p_*f_*f^!p^!(\1_k[1/e])  \xrightarrow{p_*(\varepsilon_f)}  p_*p^!(\1_k[1/e]),
    \end{equation}
    where the first equivalence is by $p\circ f=p$, and $\varepsilon_f$ is the counit of $(f_!, f^!)$.

        Moreover if $X$ is smooth and proper, then $\M^c(X) = \Sigma^{\infty} X_+[1/e]$ and $\M^c(f) = \Sigma^{\infty} f_+[1/e]$.    
       \end{lem}
\begin{proof}
The statement about $\M^c(X)$ follows since
        \[
			\M^c(X) = L(p_!\1_X)^\vee = \underline{\operatorname{Hom}}_{\SH(k)[1/e]}\bigl(L(p_!\1_X),\1_k[1/e]\bigr) \simeq p_*p^!\1_k[1/e].
		\]
For the statement about $\M^c(f)$, note that for an isomorphism $f \colon X \to Y$, with $p_X, p_Y$ the structure maps on $X,Y$, since $p_X = p_Y \circ f$, we have the composition
            \[
        p_{X,*}p_X^!(\1_k[1/e])  \simeq p_{Y,*}f_*f^!p_Y^!(\1_k[1/e])  \xrightarrow{p_{Y*}(\varepsilon_f)}  p_{Y*}p_Y^!(\1_k[1/e]).
            \]
            Now, by the duality equivalence $\underline{\operatorname{Hom}}_Y(f_* \1_X, -) \simeq f_*f^!(-) $ the above composition is the dual of $h^c_{1/e}(f)$, and the unit of the adjunction $(f^*,f_*)$ corresponds under duality to  the counit of $(f_!, f^!)$. Taking $X=Y$ we obtain that $\M^c(f)$ is the composition \eqref{eq:M^c(f)}. 
            
            Finally, for $X$ smooth and proper, $\Sigma_+^\infty X$ is dualisable in $\SH(k)$ by the smooth pullback property (\cite[Theorem 6.18(2)]{Hoy17}) and so $\M^c(X) \simeq L(p_{\#} \1_X) \simeq L(\Sigma^{\infty}X_+)$. This identification takes $\M^c(f)$ to $(\Sigma^\infty f_+)[1/e]$. 
\end{proof}
       
       \begin{prop}\label{prop:mot2}
           The composition of functors $\mr{Span}^{\mr{sc}}_k\to \mr{Span}_k \xrightarrow{\M^c} \SH(k)[1/e]$ is an exact $\infty$-functor of spans, which in turn induces a map of spectra
    \[
    \mr{mot}\colon \K(\Var_k)\to \K\Big((\SH(k)[1/e])^\omega\Big).
    \]
    The induced map on $\pi_0$ sends $[X]$ to $[\M^c(X)]$. On $\pi_1$ it sends $[X, f]$ to $[\M^c(X), \M^c(f)]$.
       \end{prop}
       \begin{proof}
       The first part follows by dualising Proposition~\ref{prop:mot1}: dualising induces an exact contravariant equivalence on compact objects, so $\M^c$ is a (covariant) exact $\infty$-functor of spans, and Theorem~\ref{thm:functorofspans} gives the desired map of $K$-theory spectra. The assertion on the induced map on $\pi_1$ now follows from Proposition~\ref{prop:maponk1}.
       \end{proof}

    \begin{defn}
        Since every object in $\SH(k)[1/e]^{\omega}$ is dualisable, we obtain a map, the \emph{trace of the identity map}
        \[
        \mathrm{tr}_{\mr{id}}\colon \K_0(\SH(k)[1/e]^{\omega}) \to \pi_0(\End_{\SH(k)[1/e]}(\1_k[1/e]))
        \]
        given on the level of objects by sending an object $\mathcal{E}$ to the morphism
        \[
        \1_k[1/e] \to \mathcal{E} \wedge \mathcal{E}^\vee \xrightarrow{\tau} \mathcal{E}^\vee \wedge \mathcal{E} \to \1_k[1/e].
        \]
    \end{defn}

        \begin{thm}[{\cite[Remark 6.6]{HSS17}}]\label{prop:tr}
            The trace of the identity map lifts to a map of spectra (which we abuse notation and also call $\mathrm{tr}_{\mr{id}}$)
            \[
            \mathrm{tr}_{\mr{id}}\colon \K\Big((\SH(k)[1/e])^\omega\Big) \to \mr{End}_{\SH(k)[1/e]}(\1_{k}[1/e]).
            \]
            Passing to infinite loop spaces, we have a map
            \[|\SH(k)[1/e]^\omega_{\mr{iso}}|\xrightarrow{\iota }\Omega^\infty K\Big((\SH(k)[1/e])^\omega\Big)\xrightarrow{\Omega^\infty\mathrm{tr}_{\mr{id}}} \Omega^\infty\mr{End}_{\SH(k)[1/e]}(\1_{k}[1/e]),\] 
            where the first map is what induces the determinant map on $\pi_1$ from Definition \ref{def:WK1Aut}, and the composition is homotopy equivalent to the trace functor in \cite[Definition 3.9]{HSS17}.
        \end{thm}

    To prove Theorem \ref{maintheorem:Categorical}, it will be helpful to examine the proof of Theorem $\ref{prop:tr}$, which a simple case of the previous theorem in loc. cit. (i.e. \cite[Theorem 6.5]{HSS17}).
        
        \begin{defn}\label{defn:derivedEulerChar}  Define the \emph{derived $\mathbb{A}^1$-Euler characteristic} to be the map given by composing the maps in Theorem~\ref{prop:tr} and Proposition~\ref{prop:mot2}
				\[
				\chi : = \mathrm{tr}_{\mr{id}} \circ \mot \colon \K(\Var_k) \to \mathrm{End}_{\SH(k)[1/e]}(\1_k[1/e]). \]
            Let $\chi_1$ be the map induced by $\chi$ on the first homotopy groups, so that
	\[
	\chi_1 \colon \K_1(\Var_k) \to \pi_1 \End_{\SH(k)[1/e]} (\1_k[1/e]).
	\]
            \end{defn}

            \begin{rem}
    For $f\colon X \to X$ an automorphism of a smooth projective variety, there is another natural associated element in the homotopy groups of motivic spheres: we could use the categorical trace, $\mathrm{tr}_{\mathbb{A}^1}(f)$, defined to be the endomorphism
    \[
    \1_k \to \Sigma^{\infty}X_+ \wedge (\Sigma^{\infty}X_+)^\vee \xrightarrow{\Sigma^{\infty}f \wedge \id} \Sigma^{\infty}X_+ \wedge (\Sigma^{\infty}X_+)^\vee\xrightarrow{\tau}  (\Sigma^{\infty}X_+)^\vee \wedge \Sigma^{\infty}X_+ \rightarrow \1_k,
    \]
    where we omit $[1/e]$ for notational simplicity. Here, we see that $\mathrm{tr}_{\mathbb{A}^1}(\id_X) = \chi_0(X)$. This trace lives in $\pi_0(\End_{\SH(k)}(\1_k))$, so gives a different invariant to $\chi_1([X,f])$. The trace of non-identity automorphisms has been used in recent work such as \cite{BHSVW}. It is beyond the scope of this paper to analyse the comparisons between $\mathrm{tr}_{\mb{A}^1}(f)$ and $\chi_1([X,f])$.
\end{rem}

            \begin{lem}
			The map $\chi$ induces on connected components is the $\mathbb{A}^1$-Euler characteristic 
            \[\pi_0 \chi = \chi_0\colon  \K_0(\Var_k) \to  \pi_0(\mathrm{End}_{\SH(k)[1/e]}(\1_k[1/e])).
            \]
            \end{lem}
            \begin{proof}
        This is immediate by construction. Indeed, for a variety $p \colon X \to k$, we have
        \[
        \pi_0\chi ([X]) = \mathrm{tr}_{\mathbb{A}^1} (\id_{\M^c(X)}) = \mathrm{tr}_{\mathbb{A}^1} ( \id_{(p_! \1_X)[1/e]}) = \chi_0 (X).
        \]
      In the above equation, the second equality follows from the invariance of the trace of the identity under duality, and the third equality follows automatically by Definition~\ref{defn:categoricalchi0}.
\end{proof}

    \begin{rem}
        The spectra $\K(\Var_k)$ and $\mathrm{End}_{\SH(k)}(\1_k)$ both have additional structure: they are $\mb{E}_\infty$-ring spectra. The construction of $\chi$ does not tell us anything about whether $\chi$ respects the $\mb{E}_\infty$-ring structures, see Remark \ref{rem:Einfty} for more discussion.\end{rem}

\begin{cor}\label{cor:k(f)}
	Let $X$ be any $k$-variety and let $f\in\Aut(X)$, then
	\[
	\chi_1([X,f])	=	\pi_1 \mathrm{tr}_{\mr{id}}\bigl([\M^c(X),\M^c(f)]\bigr).
	\]
\end{cor}
\begin{proof}
   This follows immediately from Proposition $\ref{prop:mot2}$.
\end{proof}
	
	\begin{cor}\label{cor:chi1-localisation}
		Let $f\in\Aut(X)$ and let $Z\subset X$ be an $f$-invariant closed
		subvariety. Let $U$ denote $X\setminus Z$. Then
		\[	
        [\M^c(X),\M^c(f)]=[\M^c(Z),\M^c(f|_Z)] + [\M^c(U),\M^c(f|_U)]
		\]
    in $\K_1((\SH(k)[1/e])^\omega)$. In particular, $\chi_1([X,f]) = \chi_1([Z,f|_Z])+\chi_1([U,f|_U])$.
	\end{cor}
	
	\begin{proof}
	The identities follow from Proposition~\ref{prop:mot2} and Proposition~\ref{prop:Autocutanpaste}.
	\end{proof}
	
	\subsection{The 2-categorical trace}
	
	An upshot of Corollary \ref{cor:k(f)} is that $\chi_1([X,f])$ depends only on the autoequivalence $\M^c(f)$ of $\M^c(X)$. With this observation in hand, we are now set to describe the composite $\chi_1\circ\det_X$ in terms of a 2-categorical trace. So far, we have described the map $\mot_1 \circ \mathrm{det}_X \colon\Aut(X)\to\K_1\big(\SH(k)[1/e])^\omega\big)$, so it only remains to describe the trace map $\K_1\big((\SH(k)[1/e])^\omega\big)\to\pi_1\End_{\SH(k)[1/e]}(\1_k[1/e])$.
	
	\begin{lem}\label{lem:trace}
		Let $\mc{E} \in \SH(k)[1/e]$ be strongly dualisable. Let $F\colon \mc{E}\to\mc{E}$ be an automorphism. The 2-trace $\pi_1 \mathrm{tr}_{\mr{id}} \colon \K_1\big((\SH(k)[1/e])^\omega\big)\to\pi_1\End_{\SH(k)[1/e]}(\1_k[1/e])$ sends the element $[\mc{E},F]$ of $\K_1\big((\SH(k)[1/e])^\omega\big)$ to the class of the diagram
		\begin{equation}\label{eq:trace diagram}
			\begin{tikzcd}
				&|[alias=A]| \mc{E}\wedge\mc{E}^\vee\arrow[rr,"\id\wedge\id"]\arrow[dd,"F\wedge(F^{-1})^\vee"]&&\mc{E}\wedge\mc{E}^\vee\arrow[rr,"\tau"]\arrow[ddll,Rightarrow]\arrow[dd,"F\wedge(F^{-1})^\vee"] && \mc{E}^\vee\wedge\mc{E}\arrow[dr,""{name=C,below,inner sep=1pt}]\arrow[dd,"(F^{-1})^\vee\wedge F"']\arrow[ddll,Rightarrow]&\\
				\1_k\arrow[ur]\arrow[dr,""{name=B,inner sep=1pt}] &&&&&& \1_k,\\
				&\mc{E}\wedge\mc{E}^\vee\arrow[rr,"\id\wedge\id"]&&\mc{E}\wedge\mc{E}^\vee\arrow[rr,"\tau"] && |[alias=D]| \mc{E}^\vee\wedge\mc{E}\arrow[ur]&
				\arrow[Rightarrow,from=A,to=B]
				\arrow[Rightarrow,from=C,to=D]
			\end{tikzcd}
		\end{equation}
		where the four 2-morphisms are
		\begin{align*}
			(F\wedge(F^{-1})^\vee)\circ\coev_\mc{E}&\to\coev_\mc{E},\\
			\alpha\wedge\id\colon (F\wedge(F^{-1})^\vee)\circ\id\wedge\id&\cong \id\wedge\id\circ(F\wedge(F^{-1})^\vee),\\
			(F^{-1})^\vee\wedge F\circ\tau &\simeq\tau\circ(F\wedge(F^{-1})^\vee),\\
			\ev_\mc{E}&\to\ev_\mc{E}\circ((F^{-1})^\vee\wedge F).
		\end{align*}
	\end{lem}

	\begin{proof}
		In their proof of Theorem~\ref{prop:tr}, the authors of \cite{HSS17} describe $\pi_1(\mathrm{tr}_{\mr{id}}\circ\iota)$ via \cref{eq:trace diagram} \cite[Diagram (2.1)]{HSS17}. We will explain the notation necessary to deduce \cref{eq:trace diagram} from \cite[Diagram (2.1)]{HSS17}.
		
		Since $\SH(k)[1/e]$ is symmetric monoidal under the smash product $\wedge$, there is a natural equivalence $\tau\colon \mc{A}\wedge\mc{B}\to\mc{B}\wedge\mc{A}$ for any $\mc{A},\mc{B}\in\SH(k)[1/e]$. \cref{eq:trace diagram} is the trace of the 1-morphism $(F,\alpha)\colon (\mc{E},\id)\to(\mc{E},\id)$, where $\alpha\colon F\circ\id_\mc{E}\to\id_\mc{E}\circ F$ is the identity 2-morphism. One can check that the right dual of $F$ is $F^{-1}$. Using this notation and interposing $\tau$, we thus obtain \cref{eq:trace diagram} from \cite[Diagram (2.1)]{HSS17}.
	\end{proof}

    Lemma $\ref{lem:trace}$ is the key ingredient to proving Theorem \ref{maintheorem:Categorical}.

    \begin{proof}[Proof of Theorem \ref{maintheorem:Categorical}]
    Recall that we need to show that $\chi_1([X,f])$ is given by the element of $\pi_1(\End_{\SH(k)[1/e]}(\1_k[1/e]))$ described by the following diagram
    \begin{equation*}
				\begin{tikzcd}
				&|[alias=A]|  X \wedge X^\vee\arrow[rr,"\id\wedge\id"]\arrow[dd," f\wedge f^{\vee, -1}"]&& X \wedge X^\vee \arrow[rr,"\tau"]\arrow[ddll,Rightarrow]\arrow[dd,"f\wedge f^{\vee, -1}"] && X^\vee \wedge  X \arrow[dr,""{name=C,below,inner sep=1pt}]\arrow[dd,"f^{\vee, -1}\wedge f"']\arrow[ddll,Rightarrow]&\\
				\1_k[1/e]\arrow[ur]\arrow[dr,""{name=B,inner sep=1pt}] &&&&&& \1_k[1/e],\\
				& X \wedge X^\vee \arrow[rr,"\id\wedge\id"]&& X \wedge X^\vee\arrow[rr,"\tau"] && |[alias=D]| X^\vee \wedge X \arrow[ur]&
				\arrow[Rightarrow,from=A,to=B]
				\arrow[Rightarrow,from=C,to=D]
			\end{tikzcd}
		\end{equation*}
    where we write $X,f$ to mean $\Sigma^{\infty} X_+[1/e]$ and $\Sigma^{\infty}f_+[1/e]$. 
    The result then follows instantly by combining Lemma $\ref{lem:McProperties}$, Corollary $\ref{cor:k(f)}$ and Lemma $\ref{lem:trace}$.
    \end{proof}

Theorem \ref{maintheorem:Categorical} gives us a purely categorical description of the group homomorphism $\chi_1\colon  \K_1(\Var_k) \to \pi_{1} \mathrm{End}_{\SH(k)[1/e]}(\1_k[1/e])$. While this categorical description has nice formal properties, it is not clear from this how to meaningfully compute $\chi_1([X,f])$.

\section{The image in Hermitian $\K$-theory}\label{section:Hermitian}
For the rest of the paper, assume that $k$ is a perfect field of characteristic $\neq 2$. In this section we use Hermitian $\K$-theory to extract geometric data from the map $\chi_1 \colon \K_1(\Var_k) \to \pi_1(\End_{\SH(k)}(\1_k))$ obtained in Theorem \ref{maintheorem:Categorical}. We discuss motivation for using Hermitian $\K$-theory in this way, then use Theorem \ref{maintheorem:Categorical} to express the image of $\chi_1([X,f_X])$ in Hermitian $\K$-theory in terms of computable geometric invariants.

For this section, write $\KQ$ to be the object of $\SH(k)$ representing Hermitian $\K$-theory (see \cite[Theorem 1.3]{PaninWalter2019} for a construction, where this spectrum is denoted $\mathbf{BO}$). By construction (see \cite[Definition 7.1, Theorem 1.3]{PaninWalter2019} \cite[Definition 2.10]{Schlichting}), we have that $\rm{Map}_{\SH(k)}(\1_k, \KQ)$ is the spectrum given by Waldhausen's $\mathrm{S}_\bullet$-construction on the category of non-degenerate symmetric bilinear forms over $k$. 

\subsection{Background motivation and $\GW(k)$}
The motivation for using Hermitian $\K$-theory to study $\chi_1$ comes from the $\mathbb{A}^1$-Euler characteristic $\chi_0: \K_0(\Var_k) \to \pi_{0} \mathrm{End}_{\SH(k)}(\1_k)$. We first note that the target of this map has a more classical formulation.

	\begin{defn}\label{defn:GWk}
    The \emph{Grothendieck--Witt ring} of $k$, $\GW(k)$, is the Grothendieck group completion of isometry classes of non-degenerate symmetric bilinear forms on finite dimensional $k$-vector spaces. For $V$ a finite dimensional vector space and $B\colon V \times V \to k$ a non-degenerate symmetric bilinear form, write $[(V,B)] \in \GW(k)$ for the associated class.
    
    By \cite[\S2, Theorem 4.1]{Lam05}, $\GW(k)$ is generated by elements $\langle a \rangle$ for $a\in k^\times$, which are classes of one-dimensional forms $(x,y) \mapsto axy$, subject to certain relations. Write $\mathbb{H} := \langle -1 \rangle + \langle 1 \rangle$, which corresponds to the bilinear form on the $2$-dimensional vector space $\mathrm{span}\{x,y\}$ where $\mathbb{H}(x,x) = \mathbb{H}(y,y) = 0$ and $\mathbb{H}(x,y)=1$. By definition of Hermitian $\K$-theory, there is a canonical isomorphism $\pi_{0,0}(\KQ) \cong \GW(k)$.
	\end{defn}

\begin{thm}[Theorem 6.4.1 \cite{Morel2004}, Corollary 6.2 of \cite{Ananyevskiy2016}]
The unit map of Hermitian $\K$-theory, $u:\1_k \to \KQ$ induces an isomorphism $u\colon \pi_0( \mathrm{End}_{\SH(k)}(\1_k)) \cong \GW(k)$.
\end{thm}
Since $u$ is a canonical morphism, and since \cite[Theorem 6.4.1]{Morel2004} is independent of the connection with Hermitian $\K$-theory, most literature identifies $\pi_0\mathrm{End}_{\SH(k)}(\1_k)$ with $\GW(k)$. We avoid doing this in this paper since $u$ does not induce canonical isomorphisms on higher homotopy groups (see Theorem $\ref{thm:RSO}$).

The unit map is extremely useful for computing $\chi_0$. While it is possible to use the categorical description of $\chi_0$ from Definition $\ref{defn:categoricalchi0}$ to do some computations (Hoyois computes $\chi_0( \Spec(L))=[\mathrm{Tr}_{L/k}]$ in \cite[Theorem 1.9]{Hoy14}, and Hoyois and Levine compute $\chi_0(X)$ for $X$ a cellular variety in \cite[Proposition 1.5]{LevineAspects}), this formulation is not amenable to computing the $\mathbb{A}^1$-Euler characteristic of more complicated varieties. In practice, we use a description of $\chi_0$ in terms of the Hodge cohomology, coming from \cite[Theorem 1.3]{LR20}.

\begin{defn}\cite[Construction 1.2]{LR20}
    Let $X$ be a smooth projective variety over $k$. Let $n=\dim{X}$. Define the \emph{Hodge complex} to be the perfect complex of $k$-vector spaces
    \[
    \Hdg(X) :=\bigoplus_{i,j=0}^n \H^i(X,\Omega^j_{X/k})[j-i].
    \]
    Equip $\Hdg(X)$ with the symmetric non-degenerate bilinear form induced by the cup product and trace map
\[
B_{\Hdg} \colon \H^i(X,\Omega^j_{X/k})\otimes \H^{n-i}(X,\Omega^{n-j}_{X/k})\xrightarrow{\smile}\H^n(X,\Omega^n_{X/k})\xrightarrow{\Tr}k.
\]
 Write $\chi^{\Hdg}(X) \in \GW(k)$ to mean the class corresponding to the complex $\Hdg(X)$ equipped with the bilinear form $B_{\Hdg}$. 
\end{defn}

\begin{thm}[Theorem 1.3 of \cite{LR20}]\label{thm:GaussBonnet}
For $X/k$ a smooth projective variety 
\[u \circ \chi_0(X) = \chi^{\Hdg}(X) \in \GW(k).\]
\end{thm}
\begin{rem}
    Theorem $\ref{thm:GaussBonnet}$ is often stated without the homomorphism $u$, since $u$ is a canonical isomorphism. We choose to include writing $u$ since the proof uses the comparison with Hermitian $\K$-theory (see \cite[Theorem 8.6]{LR20}).
\end{rem}

This allows us to compute $\chi_0(X)$ for more complicated varieties such as hypersurfaces \cite[Theorem 3.12]{LPLS24} and intersections of two same degree hypersurfaces, \cite[Theorem 1]{Viergever}. This also allows us to relate parts of $\chi_0(X)$ to other arithmetic invariants.
\begin{defn}\label{defn:rankdisc}
        There are two homomorphisms out of $\GW(k)$ which are of particular interest to us. For $B$ a non-degenerate symmetric bilinear form on a vector space $V$, the \emph{rank} of $B$ is defined to be $\mathrm{dim}_k(V)$. This furnishes a ring homomorphism 
        \begin{align*}
    \mathrm{rank}: \GW(k) &\to \Z\\
    [(V,B)] &\mapsto \mathrm{dim}_k(V).
        \end{align*}
        Similarly, the assignment $[(V,B)] \mapsto \mathrm{det}(B) \in k^\times/k^{\times 2}$ furnishes a homomorphism of abelian groups, called the \emph{discriminant} \begin{align*}
            \mathrm{disc}\colon \GW(k) &\to k^\times/k^{\times 2} \\
            [(V,B)] &\mapsto \mathrm{det}(B).
        \end{align*}
    \end{defn}
For $X$ a smooth projective variety Theorem $\ref{thm:GaussBonnet}$ allows us to compare $\chi_0(X)$ to more traditional invariants. For example, \cite[Remark 1.3.1]{LevineAspects} demonstrates that $\mathrm{rank}(\chi_0(X))$ isthe compactly supported Euler characteristic of the topological space $X(\mathbb{C})$. Moreover, \cite[Theorem 2.29]{YauZaslow} shows that $\mathrm{disc}(\chi_0(X))$ can be computed by looking at the Galois representation $\mathrm{det}_\ell(X) := \bigotimes_{i=0}^{\infty} \mathrm{det}(H^i_{\et}(X_{\kbar}, \mathbb{Q}_\ell))^{\otimes (-1)^i}$, which can in turn be calculated by looking at the point counting $\zeta$-functions of the reduction of the variety to finite fields.

\subsection{The first homotopy group of $\KQ$}

We recall the results of \cite{RSO19} to study the morphism $\pi_{1,0}(\1_k) \to \pi_{1,0}(\KQ)$ induced by the unit map $\1_k \to \KQ$ in $\SH(k)$.

\begin{thm}[{R\"ondigs--Spitzweck--\O stv\ae r, \cite[Theorem 5.5]{RSO19}}]\label{thm:RSO}
    When the characteristic of $k$ is not $2$ or $3$ the unit map of Hermitian $\K$-theory $u\colon \1_k \to \KQ$ induces a surjection $\pi_{1,0}(\1_k) \twoheadrightarrow \pi_{1,0}(\KQ)$ with kernel isomorphic to $\K^M_2(k)/(24)$.
\end{thm}
\begin{rem}
    Assuming $\mathrm{char}(k) \neq 2$ is crucial for many results of the paper, due to the reliance on Hermitian $\K$-theory. In characteristic $3$, we still obtain the unit map $\1_k \to \pi_{1,0}(\KQ)$, but we are unable to use Theorem $\ref{thm:RSO}$ to ensure whether it is surjective. We can instead show this directly, however, in Corollary $\ref{lem:surjectivityforC}$.
\end{rem}

\begin{cor}\label{cor:k=CIso}
    When $k=\mathbb{C}$, the unit map is an isomorphism.
\end{cor}
\begin{proof}
    When $k=\mathbb{C}$, $\K^M_2(\mathbb{C})$  is a divisible group, so $\K^M_2(\mathbb{C})/(24)$ is trivial and the unit map in Theorem $\ref{thm:RSO}$ is an isomorphism.
\end{proof}
\begin{rem}
    Note the difference between Theorem $\ref{thm:RSO}$ with the $\pi_{0,0}$-case where the unit map $\pi_{0,0}(\1_k) \to \pi_{0,0}(\KQ)$ is an isomorphism. In particular, while we use the unit map to Hermitian $\K$-theory to study $\chi_1$ in this paper, we (in general) lose information. It is beyond the scope of this paper to study the information contained in $\K^M_2(k)/(24)$.
\end{rem}

\begin{defn}\label{def:chi1KQ}
    Let $\KQ[1/e]$ denote the image of $\KQ$ in $\SH(k)[1/e]$, and write $u_{1/e}$ to mean the localised unit map $\1_k[1/e]\to \KQ[1/e]$. Define the map \[\chikq\colon \K_1(\Var_k) \to \pi_{1,0}(\KQ[1/e])\] to be the composition $u_{1/e} \circ \chi_1$. 
\end{defn}
In Corollary $\ref{cor:KQLocal}$, we show that $\pi_{1,0}(\KQ)\cong \pi_{1,0}(\KQ[1/e])$, so this localisation does not affect computations. The rest of this section is dedicated to giving an explicit description of $\chi_1^{\KQ}([X,f_X])$ when $X$ is a smooth projective variety and $f_X\colon X \to X$ is an automorphism. It will be first helpful to have an explicit description of $\pi_{1,0}(\KQ)$. For this, we introduce an invariant of automorphisms of bilinear forms, the spinor norm, which has been well studied in other literature. See \cite[\S 9.3]{Sch85} for a textbook treatment.

\begin{defn}\label{defn:spinornorm}
    Let $V$ be a $k$-vector space, and $B\colon V \times V \to k$ be a non-degenerate symmetric bilinear form. A map $\sigma \colon V\to V$ is an \emph{automorphism of $(V,B)$} if $\sigma$ is an isomorphism and $B( \sigma(v_1), \sigma(v_2)) = B(v_1, v_2)$ for all $v_1,v_2\in V$. Write $\sigma\colon(V,B)\to(V,B)$.

    Given $v\in V$, we say that $\sigma$ is a \emph{reflection by $v$}\footnote{There is a potential point of confusion here: reflection \emph{by} $v$ means that $\sigma$ reflects $V$ across the hyperplane orthogonal to $v$, not that $\sigma$ reflects $V$ across the span of $v$.} if for every $w\in V$, we have
    \[\sigma(w)=w-2\frac{B(v,w)}{B(v,v)}v.\]
    If $\sigma$ is a reflection by $v$, the \emph{spinor norm} of $\sigma$ is defined as $\sn(\sigma)=B(v,v)$. Note that if $\sigma$ is a reflection by $v$, then $\sigma$ is a reflection by $\lambda v$ for $\lambda\in k^\times$, so for $\sn(\sigma)$ to be well defined, it takes values in  $k^\times/k^{\times 2}$. Reflections generate $\Aut(V,B)$ by the Cartan--Dieudonné theorem (see \cite[Chapter 5, \S 1]{Lam05}), and we can extend $\sn$ to a group homomorphism
    \[\sn:\Aut((V,B))\to k^\times/k^{\times 2}.\]
\end{defn}

\begin{thm}[{\cite[Corollary 4.7.7]{Bass}}]\label{thm:bass}
    Every element of $\pi_{1,0}\KQ$ can be written as a symbol $[(V,B), \sigma]$ where $V$ is a finite dimensional vector space over $k$, $B\colon V \times V \to k$ is a non-degenerate symmetric bilinear form and $\sigma \colon (V,B) \to (V,B)$ is an automorphism. Moreover, there is an isomorphism 
    \begin{align*}
        \KQ_1(k) &\cong \{\pm 1\} \times k^\times/k^{\times 2}\\
        [(V,B), \sigma] &\mapsto (\det(\sigma), \sn(\sigma)).
    \end{align*}
\end{thm}

In particular, Theorem $\ref{thm:bass}$ means that the morphism $\chikq$ can be considered to be a group homomorphism $\K_1(\Var_k) \to \{\pm 1\} \times k^{\times}/k^{\times 2}$. Even before we understand how to compute $\chikq$, we obtain the following corollary.
\begin{cor}\label{cor:KQLocal}
    The localisation $L:\SH(k)\to \SH(k)[1/e]$ induces the canonical isomorphism $\pi_{1,0}(\KQ)=  \pi_{1,0}(\KQ[1/e]) = \{\pm 1 \} \times k^{\times}/k^{\times 2}$. 
\end{cor}
\begin{proof}
The adjunction properties of the localisation $\SH(k)\to \SH(k)[1/e]$ give the isomorphism $\pi_{1,0}(\KQ[1/e])\cong \pi_{1,0}(\KQ)[1/e]$, but since $\pi_{1,0}(\KQ)$ is $2$-torsion by Theorem \ref{thm:bass}, $e$ is already invertible in $\pi_{1,0}(\KQ)$.
\end{proof}
In light of Corollary $\ref{cor:KQLocal}$, we view $\chikq$ as a map $\K_1(\Var_k) \to \pi_{1,0}(\KQ)$.
\begin{cor}\label{cor:SquaresVanish}
    Let $f\colon X \to X$ be an automorphism of a smooth projective variety, and suppose that there exists $g\colon X \to X$ such that $f = g \circ g$. Then $\chikq([X,f]) = 0$. 
\end{cor}
\begin{proof}
    This follows immediately since by Theorem $\ref{thm:bass}$, $\pi_{1,0}(\KQ)$ is $2$-torsion, and 
    \[
    \chikq([X,f]) = 2 \cdot \chikq([X,g]) = 0. \qedhere\]
\end{proof}

\subsection{The Hodge complex and $\KQ_1$}
It therefore remains to compute $\chikq([X,f])$. We do this by showing that $f$ induces an automorphism of $\chi^{\Hdg}(X)$ and compare this to construction in Theorem \ref{maintheorem:Categorical}.

\begin{defn}\label{defn:chihdg}
Let $f \colon X  \to Y$ be a morphism of $k$-varieties and consider the differential maps $df_j \colon f^* \Omega^j_{Y/k} \to \Omega^j_{X/k}$.
The \emph{pullback in Hodge cohomology, $f^*_{\Hdg}$} is defined on cohomology groups by the composition
\[
f^*_{\Hdg} \colon \H^i(Y,\Omega^j_{Y/k}) \xrightarrow{f^*} \H^i(X,f^*\Omega^j_{Y/k}) \xrightarrow{\H^i(df_j)} \H^i(X,\Omega^j_{X/k}).
\]
Define the map $f^*_{\Hdg}$ to be the induced morphism of complexes 
\[f^*_{\Hdg} \colon \Hdg(Y) \to \Hdg(X).\]
\end{defn}

\begin{lem}\label{lem:fHdgCompat}
    Let $f \colon X  \to Y$ be a morphism of $k$-varieties. Then
    \begin{enumerate}
    \item
    The morphism $f^*_{\Hdg}$ commutes with cup products.
    \item 
    For $X$ and $Y$ smooth and proper over $k$ of dimension $n$, and $f$ an isomorphism, we have that $f^*_{\Hdg}$ commutes with $\Tr$.
    \item 
    If $X=Y$ and $f \in \Aut(X)$, then $f^*_{\Hdg}$ induces an automorphism of symmetric bilinear forms on $(\Hdg(X), B_{\Hdg})$.
    \end{enumerate}
\end{lem}

	\begin{proof}
		The maps $df_j \colon f^*\Omega^j_{Y/k}\longrightarrow \Omega^j_{X/k}$ are given by exterior powers of $df_1$, so the diagram
		\[
		\begin{tikzcd}
			f^*\Omega^p_{Y/k}\otimes f^*\Omega^q_{Y/k}
			\ar[r,"\wedge"]
			\ar[d,"df_p\otimes df_q"']
			&
			f^*\Omega^{p+q}_{Y/k}
			\ar[d,"df_{p+q}"]
			\\
			\Omega^p_{X/k}\otimes\Omega^q_{X/k}
			\ar[r,"\wedge"]
			&
			\Omega^{p+q}_{X/k}
		\end{tikzcd}
		\]
		commutes. By naturality of the cup product with respect to morphisms of the
		coefficient sheaves, for $\eta\in \H^i(X,f^*\Omega^j_{Y/k})$ and $\xi\in \H^{i'}(X,f^*\Omega^{j'}_{Y/k})$, we have
		\[
		\H^{i+i'}(df_{j+j'})(\eta\smile\xi)
		=
		\H^{i}(df_j)(\eta)\smile \H^{i'}(df_{j'})(\xi).
		\]
		This allows us to prove assertion (1): if $\alpha \in \H^i(Y, \Omega^j_{Y/k})$ and $\beta \in \H^{i'}(Y, \Omega^{j'}_{Y/k})$, then
		\[
		\begin{aligned}
			f^*_{\Hdg}(\alpha\smile\beta)
			&=
			\H^{i+i'}(df_{j+j'})\bigl(f^*(\alpha\smile\beta)\bigr)
			\\
			&=
			\H^{i+i'}(df_{j+j'})
			\bigl(f^*\alpha\smile f^*\beta\bigr)
			\\
			&=
			\H^i(df_j)(f^*\alpha)\smile \H^{i'}(df_{j'})(f^*\beta)
			\\
			&=
			f^*_{\Hdg}(\alpha)\smile f^*_{\Hdg}(\beta).
		\end{aligned}
		\]

For assertion (2), since $f\colon X\to Y$ is a $k$-isomorphism, the structure morphisms satisfy $p_X=p_Y\circ f$. Applying
base-change compatibility of the Grothendieck duality trace (see e.g.\cite[Lemma~48.7.1]{stacks}) yields compatibility of the trace maps under base-change $f^*p_Y^!k\simeq p_X^!k$. Since $X$ and $Y$ are smooth of dimension $n$, we have
\[
p_Y^!k\simeq\omega_{Y/k}[n],
\qquad
p_X^!k\simeq\omega_{X/k}[n].
\]
Under these identifications, the base-change isomorphism
becomes 
\[
df_n[n]\colon f^*\omega_{Y/k}[n]\xrightarrow{\sim}\omega_{X/k}[n].
\]
Consequently we obtain assertion (2), i.e.,
\[
\Tr_X\circ f^*_{\Hdg}
=
\Tr_X\circ \H^n(df_n)\circ f^*
=
\Tr_Y.
\]
		
		For assertion (3) suppose that $X=Y$ and
		$f\in\Aut_k(X)$. Since $f$ is an isomorphism,
		$(f^{-1})^*_{\Hdg}$ is inverse to $f^*_{\Hdg}$, so
		$f^*_{\Hdg}$ is an automorphism of $\Hdg(X)$. The result follows by parts (1) and (2), since $B_{\Hdg}$ is given by composition of cup product and trace.
		\end{proof}
        
\subsection{The comparison theorem}
We come to the main theorem of this section.
\begin{thm}\label{thm:weil precise}
     Let $X$ be a smooth projective variety over $k$ and let $f\colon X \to X$ be an automorphism. Then \[\chikq([X,f]) = [(\Hdg(X), f^*_\Hdg)] \in\KQ_1(k).\]
\end{thm}
\begin{proof}
    Let us first work in $\SH(k)$ without inverting the characteristic, which we will extend to the characteristic inverted case later. We examine the constructions in \cite[Theorem~1.3]{LR20}, in which the top and bottom composites of the diagram
    \begin{equation}\label{eq:pushforward by u}
    \begin{tikzcd}
       & \Sigma^\infty X_+\wedge(\Sigma^\infty X_+)^\vee\arrow[r,"\tau"]\arrow[dd,"f\wedge (f^{-1})^\vee",near start] & (\Sigma^\infty X_+)^\vee\wedge\Sigma^\infty X_+\arrow[dr]\arrow[dd,"(f^{-1})^\vee\wedge f"',near end] &&\\
       \1_k\arrow[ur]\arrow[dr] &&& \1_k\arrow[r,"u"] & \KQ\\
       & \Sigma^\infty X_+\wedge(\Sigma^\infty X_+)^\vee\arrow[r,"\tau"] & (\Sigma^\infty X_+)^\vee\wedge\Sigma^\infty X_+\arrow[ur] &&
    \end{tikzcd}
    \end{equation}
    both represent $u\circ\chi_0(X)$ (where we write $f$ as an abuse of notation for $\Sigma^{\infty} f$). It remains to calculate the isomorphism $u\circ\chi_0(X)\to u\circ\chi_0(X)$ induced by $f\wedge(f^{-1})^\vee$ and the 2-morphisms implicit in \cref{eq:pushforward by u}.

    Let $n=\dim{X}$. The proof of \cite[Theorem 1.3]{LR20} proceeds via \cite[Theorem 1.5]{LR20}, in which Levine and Raksit prove that $u\circ\chi_0(X)=\pi_*e(T_{X/k})$, where $\pi\colon X\to\Spec{k}$ is the structure map, $T_{X/k}\to X$ is the tangent bundle, and $e(T_{X/k})\in\KQ^{2n,n}(X;\omega_{X/k})$ is the Euler class. To define this Euler class, let $s\colon X\to T_{X/k}$ denote the zero section, and let $u_*\1_X\in\KQ(X)$ denote the unit object of Hermitian $\K$-theory. Then the Euler class $e(T_{X/k})\colon =s^*s_*u_*\1_X$.
    Recall that classes in $\KQ^{2n,n}(X;\mc{L})$ can be represented by perfect complexes $C^\bullet$ of vector bundles on $X$, together with a pairing $C^\bullet\otimes C^\bullet\to\mc{L}[n]$. 
The class $e(T_{X/k}) \in\KQ^{2n,n}(X;\omega_{X/k})$ is represented by the Hermitian complex of coherent sheaves on $X$ given by $\bigoplus_{j=0}^n\Omega^j_{X/k}[j]$ equipped with the pairing given by
 \[\Omega^j_{X/k}[j]\otimes \Omega^{n-j}_{X/k}[n-j]\to \omega_{X/k}[n].\]

    The pullback $f^*e(T_{X/k})\in\KQ^{2n,n}(X;f^*\omega_{X/k})$ is represented by the complex of coherent sheaves on $X$ given by $\bigoplus_{j=0}^nf^*\Omega^j_{X/k}[j]$ equipped with the pairing given by the product
    \[f^*\Omega^j_{X/k}[j]\otimes f^*\Omega^{n-j}_{X/k}[n-j]\to f^*\omega_{X/k}[n].\]

	Since $f$ is an automorphism, the differential map $df_\bullet\colon f^*\Omega^j_{X/k} \rightarrow \Omega^j_{X/k}$ is an isomorphism of exterior power algebras, and so the following diagram commutes
			 \begin{equation}\label{eq:euler classes}
			\begin{tikzcd}
				f^*\Omega^j_{X/k}[j]\otimes
				f^*\Omega^{n-j}_{X/k}[n-j]
				\ar[r,"\wedge"]
				\ar[d,"\cong", "df_j\otimes df_{n-j}"']
				&
				f^*\omega_{X/k}[n]
				\ar[d,"df_n", "\cong" ']
				\\
				\Omega^j_{X/k}[j]\otimes
				\Omega^{n-j}_{X/k}[n-j]
				\ar[r,"\wedge"]
				&
				\omega_{X/k}[n].
			\end{tikzcd}
			\end{equation}
			We now describe the automorphism induced after applying $\pi_*$. Since
$\pi\circ f=\pi$, the adjunction unit $e(T_{X/k}) \to f_*f^*e(T_{X/k})$ induces a map on the Hermitian complexes represented by the Euler classes
\[
\pi_*e(T_{X/k})
\longrightarrow
\pi_*f_*f^*e(T_{X/k})
=
\pi_*f^*e(T_{X/k}) = \pi_*e(f^*T_{X/k}).
\]
On the summand $\H^i(X,\Omega^j_{X/k})$ of $\pi_*e(T_{X/k})$, this is the pullback
\[
\H^i(X,\Omega^j_{X/k})
\xrightarrow{f^*}
\H^i(X,f^*\Omega^j_{X/k}).
\]
Composing with the isomorphism induced by $df_\bullet$ gives
\[
\H^i(X,\Omega^j_{X/k})
\xrightarrow{f^*}
\H^i(X,f^*\Omega^j_{X/k})
\xrightarrow{\H^i(df_j)}
\H^i(X,\Omega^j_{X/k}),
\]
which is precisely $f^*_{\Hdg}$, by definition. By Lemma $\ref{lem:fHdgCompat}$, $f^*_{\Hdg}$ is an automorphism of the Hermitian perfect
complex $(\Hdg(X),B_{\Hdg})$ representing $u\circ\chi_0(X)$.

Tracing back the $2$-morphisms in Equation~\eqref{eq:pushforward by u}, the
automorphism of $u\circ\chi_0(X)$ induced by the middle vertical map is
exactly the composite
\[
\pi_*e(T_{X/k})
\longrightarrow
\pi_*f^*e(T_{X/k})
\longrightarrow
\pi_*e(T_{X/k})
\]
described above. Thus it is represented on the Hodge complex by $f^*_{\Hdg}$. 

It only remains to show the result when inverting the characteristic. The localisation functor $L:\SH(k)\to\SH(k)[1/e]$ sends the diagram \ref{eq:pushforward by u} to the diagram in \ref{eq:trace diagram} with $\mathcal{E}=L(\Sigma_+^\infty X), F=L(\Sigma_+^\infty f_+)$, composed with the map $\1_k[1/e]\to \KQ[1/e]$, which is exactly $\chikq([X,f])$. Since $L$ induces the isomorphism $\pi_1(\KQ)\to \pi_1(\KQ[1/e])$, we may identify $\chikq([X,f])= [\Hdg(X), f^*_{\Hdg}] \in\pi_{1,0}(\KQ)$, as required.
\end{proof}
This allows us to prove Theorem \ref{maintheorem:Hermitian}.

\begin{proof}[Proof of Theorem \ref{maintheorem:Hermitian}]
Recall that we need to show that $\chikq([X,f]) = (\det(f^*_{\Hdg}), \sn(f^*_{\Hdg}))$. This follows instantly by combining Theorem $\ref{thm:weil precise}$ with Theorem $\ref{thm:bass}$.
\end{proof}

This allows us to compute $\chikq([X,f])$ for $f\colon  X \to X$ an automorphism of a smooth projective variety $X$ by reducing it to a question of linear algebra. For example, we have the following result.
\begin{cor}
    When $k=\mathbb{C}$, $\chi_1([X,f]) = \mathrm{det}(f^*_{\Hdg}) \in \{\pm 1\} \cong \pi_1(\End_{\SH(k)}(\1_k))$.
\end{cor}
\begin{proof}
    This follows immediately from Theorem \ref{maintheorem:Hermitian} and Corollary $\ref{cor:k=CIso}$.
\end{proof}

\section{Example computations of $\chikq$}\label{section:examples}
In this section, we compute $\chikq([X,f])$ for some natural examples of pairs $[X,f]$. The pairs of varieties $X$ with automorphisms $f \colon X \to X$ we compute allow us to deduce results about the structure of $\K_1(\Var_k)$ in Section $\ref{section:K1Vark}$, and are required for the birational geometry results in Section $\ref{section:birational}$. We largely restrict to computing $\chikq([X,f])$ where $f$ is $2$-torsion, in part since Corollary $\ref{cor:SquaresVanish}$ implies that $\chikq([X,f])=0$ if $f$ has odd order.

\subsection{Linear algebra prerequisites}
In this subsection, we prove basic linear algebra results for computing spinor norms of certain automorphisms of vector spaces with bilinear forms, which serve as computational tools to compute $\chikq$ in later examples.

\begin{lem}\label{lem:product}
    Let $(V,B)$ and $(V',B')$ be two vector spaces with bilinear forms. Suppose $f\colon  (V,B) \to (V,B)$ is an automorphism of $(V,B)$ which is not the identity. Let $\epsilon = 0$ if $\mathrm{det}(f)=1$ and $1$ if $\mathrm{det}(f) = -1$. 
    
    Then $f \otimes \mathrm{id}$ is an automorphism of $(V \otimes V', B \otimes B')$ and 
\begin{align*}
\mathrm{det}(f \otimes \mathrm{id}) &= \mathrm{det}(f)^{\mathrm{dim}(V')}\\
\mathrm{sn}(f \otimes \mathrm{id}) &= \mathrm{sn}(f) ^{\mathrm{dim}(V')}\mathrm{disc}(V', B')^{\epsilon}.
\end{align*}
\end{lem}
\begin{proof}
    For a vector $w$ in a vector space with a bilinear form $(W,B_W)$, write $r_w\colon  W \to W$ for the reflection in $w$. First, suppose that $f=r_v$ is a reflection in $v$ for $v \in V$, so that $\mathrm{det}(f) = -1$ and $\mathrm{sn}(f) = B(v,v)$. Pick $e_1, \ldots, e_n$ an orthogonal basis for $V'$. An easy computation shows $f \otimes \mathrm{Id} = r_{v \otimes e_n} \circ r_{v \otimes e_{n-1}} \circ \ldots r_{v \otimes e_1}$. In particular, 
    \[
    \mathrm{sn}(f \otimes \mathrm{id}) = \prod B(v,v)B'(e_i, e_i) = \mathrm{sn}(f)^{\mathrm{dim}(V')} \cdot \mathrm{disc}(V',B'),\]
    which gives the result when $f$ is a single reflection.  The general result follows by applying the Cartan--Dieudonné theorem (\cite[Chapter 5, \S1]{Lam05}) to write $f$ as a product of reflections, then iteratively applying the above calculation for when $f$ is a single reflection (where we note that $\mathrm{det}(f)$ is $-1$ if we can write $f$ as an odd number of reflections and $1$ if we can write $f$ as an even number of reflections).
\end{proof}

\begin{lem}\label{negative}
    Let $\alpha \in k^\times$ and $\langle \alpha \rangle$ denote the $1$-dimensional vector space $k$ equipped with the non-degenerate symmetric bilinear form $(x,y) \mapsto \alpha xy$. Let $\sigma\colon  \langle \alpha \rangle \to \langle \alpha \rangle$ be the morphism $x \mapsto -x$. Then $\sn(\sigma) = \alpha$.
\end{lem}
\begin{proof}
    It is easy to check this is a reflection in the element $1$, and $\langle\alpha\rangle(1,1)=\alpha$.
\end{proof}
\begin{lem}\label{involutiondecomp}
    Let $V$ be a vector space equipped with a non-degenerate symmetric bilinear form $B$. Let $\sigma\colon  V \to V$ be an involution of $(V,B)$ so that $\sigma \circ \sigma = \mathrm{Id}$. Let $V^+, V^-$ be the eigenspaces for $\sigma$.  Then $(V^+)^{\perp} = V^-$, so $B$ decomposes as $B|_{V^+} \oplus B|_{V^-}$.
\end{lem}
\begin{proof}
    Note that $V = V^+ \oplus V^-$ since $\sigma$ is diagonalisable. It only remains to show that these eigenspaces are orthogonal. Since $\sigma\colon  V \to V$ is an automorphism of $(V,B)$, if we take $v \in V^+$ and $w \in V^-$, we see that $B(v,w) = B(\sigma(v), \sigma(w)) = B(v,-w) = -B(v,w)$, which gives $B(v,w) = 0$ as required. 
\end{proof}
\begin{cor}\label{SpinorInvolution}
    The spinor norm $\sn(\sigma) = \mathrm{disc}(B|_{V^-}) = \mathrm{disc}(B) \cdot \mathrm{disc}(B|_{V^+})$.
\end{cor}
\begin{proof}
    The first equality follows by diagonalising $B|_{V^-}$ and repeatedly applying Lemma $\ref{negative}$. The second equality follows since the discriminant is a homomorphism of abelian groups $\GW(k) \to k^\times/k^{\times 2}$.
\end{proof}
\begin{cor}\label{lem:sn of H}
    Let $\mathbb{H}$ denote the hyperbolic form, i.e., $\mathbb{H}$ is a $2$ dimensional vector space with basis $x,y$ equipped with a symmetric bilinear form $B$ such that $B(x,x) = B(y,y) = 0$ and $B(x,y) = 1$. Let $\sigma\colon  \mathbb{H} \to \mathbb{H}$ be the involution defined by $\sigma(x) = y$, $\sigma(y)=x$. 
    
    Then $\mathrm{sn}(\sigma) = -2$ and $\mathrm{det}(\sigma)=-1$. 
\end{cor}
\begin{proof}
    It is easy to check that the $(-1)$-eigenspace of $\sigma$ is spanned by $x-y$ and moreover we see $B(x-y, x-y) = -2$, which gives the result by Corollary $\ref{SpinorInvolution}$.
\end{proof}

\begin{lem}\label{lem:snofGalois}
    Let $L/k$ be a finite field extension, and let $\sigma \colon  L \to L$ be an involution of $L$ as a $k$-algebra, so that $\sigma$ is also an involution of the trace bilinear form on $L$. Let $L^+$ denote the fixed field of $\sigma$. Suppose that $\sigma$ is not the identity. Then
   \[\mathrm{sn}(\sigma) =  2^{[L^+:k]} \mathrm{disc}(L^+/k)\mathrm{disc}(L/k).\]

   Else, if $\sigma$ is the identity, $\mathrm{sn}(\sigma)=1$.
\end{lem}
\begin{proof}
    The case where $\sigma$ is the identity is trivial, so suppose not. By Corollary $\ref{SpinorInvolution}$, the spinor norm of $\sigma$ is the discriminant of the trace form restricted to $L^-$, so is given by $\mathrm{disc}(L/k) \mathrm{disc}( \mathrm{Tr}_{L/k}|_{L^+})$. We can compute
    \[
    \mathrm{Tr}_{L/k}|_{L^+}(x^2) = \mathrm{Tr}_{L/k}(x^2) = [L:L^+] \mathrm{Tr}_{L^+/k}(x^2).
    \]
    Since $\sigma$ is not the identity, $[L:L^+] = 2$, so $\mathrm{Tr}_{L/k}|_{L^+} = \langle 2 \rangle \cdot \mathrm{Tr}_{L^+/k}$. The discriminant is therefore given by $2^{[L^+:k]} \mathrm{disc}(L^+/k)$. This shows $\mathrm{sn}(\sigma) = 2^{[L^+:k]} \mathrm{disc}(L^+/k)\mathrm{disc}(L/k)$.
\end{proof}
\begin{cor}\label{cor:sn of quadratic}
    Let $\alpha \in k^\times$ be a non-square and $L=k(\sqrt{\alpha})$ be a quadratic extension of $k$ with $\sigma\colon  \Spec(L) \to \Spec(L)$ be the involution given by the non trivial element of $\Gal(L/k)$. Then $\chikq( [\Spec(L), \sigma]) = (-1, 2\alpha) \in \{\pm 1\} \times k^\times/k^{\times 2}$.
\end{cor}
\begin{proof}
    This follows instantly by applying Theorem \ref{maintheorem:Hermitian} to Lemma $\ref{lem:snofGalois}$, noting that the bilinear form on $\chi^{\Hdg}(\Spec(L))$ is exactly the trace form by \cite[Theorem 1.9]{Hoy14}.
\end{proof}
\begin{cor}\label{lem:surjectivityforC}
    The morphism $\chikq\colon  \K_1(\Var_k) \to \{\pm 1\} \times k^{\times} / k^{\times 2}$ is surjective.
\end{cor}
\begin{proof}
    Elements of the form $(-1, 2\alpha)$ generate the right hand group and lie in the image of $\chikq$ by Corollary $\ref{cor:sn of quadratic}$, which gives the result.
\end{proof}
We postpone examining the consequences of this corollary to Section $\ref{section:K1Vark}$.

\subsection{The invariant of canonical surface involutions}
In this subsection, we use the computations from the previous subsection to compute $\chikq([(S,f])$ for some distinguished pairs of rational surfaces equipped with an involution. Working with rational surfaces greatly simplifies computations due to the following lemma.
\begin{lem}\label{lem:rationalsurface}
Let $S$ be a smooth projective rational surface. Then the Hodge cohomology groups $\H^i(S, \Omega^j_{S/k})$ satisfy 
\[
\H^i(S, \Omega^j_{S/k}) = \begin{cases} k & \text{ if } i=j=0 \text{ or } i=j=2 \\
0 &\text{ if } i \neq j. \end{cases}
\]
Moreover, if $f\colon S \to S$ is an automorphism, then $f$ fixes $\H^0(S, \Omega^0_{S/k})$ and $\H^2(S, \Omega^2_{S/k})$.
\end{lem}
\begin{proof}
    First we see $\H^0(S, \Omega^0_{S/k}) = k$ since $S$ is projective. By Serre duality  (\cite[\S3, Corollary 7.13]{Hartshorne}) the same is true for $\H^2(S, \Omega^2_{S/k})$. 

    For the remaining cases, note that $\H^0(\mathbb{P}^2, \Omega^1_{\mathbb{P}^2/k}) = \H^0(\mathbb{P}^2, \Omega^2_{\mathbb{P}^2/k})=0$. Therefore, since $S$ is rational and the Hodge cohomology groups $\H^0(S, \Omega^i_{S/k})$ are birational invariants (\cite[Theorem 3.22]{Shafarevich}), we have that $\H^0(S, \Omega^i_{S/k}) = 0$ for $i=1,2$. The result follows by applying Serre duality to obtain the required vanishing on the other cohomology groups.

    To prove an automorphism fixes $\H^0(S, \Omega^0_{S/k})$ and $\H^2(S, \Omega^2_{S/k})$, by Serre duality it is enough to show the result for $\H^0(S, \Omega^0_{S/k})$. The result follows since $\H^0(S, \Omega^0_{S/k})$ is the group of global sections of $\mathcal{O}_{S/k}$, which are all constant.
\end{proof}
\begin{cor}\label{cor:H11Rational}
    Let $S$ be a smooth projective rational surface and let $f\colon S \to S$ be an automorphism of $S$. Then $\chikq([S, f]) = ( \mathrm{det}( f^*_{\Hdg}|_{\H^1(S, \Omega^1_{S/k})}),\sn (f^*_{\Hdg}|_{\H^1(S, \Omega^1_{S/k})})) $.
\end{cor}

We also fix the following notation.
\begin{defn}
    For $\mathcal{L}$ a line bundle on $S$, let $c_1(\mathcal{L}) \in \H^1(S, \Omega^1_{S/k})$ denote the first Chern class of $\mathcal{L}$. For  the canonical bundle $\omega_{S/k}$, we write $K_S := c_1(\omega_{S/k})$. 
\end{defn}

The three main involutions on surfaces we study are Bertini involutions on del Pezzo surfaces of degree $1$, Geiser involutions on del Pezzo surfaces of degree $2$ and involutions of rational conic bundle surfaces (of which, de Jonqui\`eres involutions are a subclass). Before getting to these surfaces, we begin with the following two examples.

\begin{lem}\label{lem:P1}
    Let $f\colon  \mathbb{P}^1_k \to \mathbb{P}^1_k$ be an automorphism of $\mathbb{P}^1_k$. Then $\chikq([ \mathbb{P}^1_k, f]) = (1,1)$.
\end{lem}
\begin{proof}
    Note that $f^*_{\Hdg}$ fixes $\H^0(\mathbb{P}^1_k, \Omega^0_{\mathbb{P}^1_k/k})$, since this group is the group of sections of the trivial bundle, which are all constant. By Serre duality (\cite[\S3 Corollary 7.13]{Hartshorne}), this means $f^*_{\Hdg}$ also fixes $\H^1(\mathbb{P}^1_k, \Omega^1_{\mathbb{P}^1_k/k})$. All other Hodge cohomology groups of $\mathbb{P}^1_k$ are trivial, so $f^*_{\Hdg}=\mathrm{Id}$, and the result follows trivially.
    \end{proof}

\begin{lem}\label{lem:surject on det}
    Let $X=\mb{P}^1_k\times\mb{P}^1_k$, and let $f\colon X\to X$ denote the map $f(x,y)=(y,x)$. Then $\chikq(f)=(-1,-2)$.
\end{lem}
\begin{proof}
    We use Theorem \ref{maintheorem:Hermitian} to reduce this to a computation of the action of $f^*_{\Hdg}$. Corollary $\ref{cor:H11Rational}$ reduces this to studying $\H^1(X, \Omega^1_{X/k})$. The K\"unneth formula  gives
    \[\H^1(X,\Omega^1_{X/k})=
            \left(\H^1(\mathbb{P}^1, \Omega^1_{\mathbb{P}^1/k}) \otimes_k \H^0(\mathbb{P}^1, \Omega^0_{\mathbb{P}^1/k})\right) \oplus 
    \left(\H^0(\mathbb{P}^1, \Omega^0_{\mathbb{P}^1/k}) \otimes_k \H^1(\mathbb{P}^1, \Omega^1_{\mathbb{P}^1/k})\right).\]

   Let $B_{1,1}$ denote the induced bilinear form given by restricting $B_{\Hdg}$ to $\H^1(X, \Omega^1_{X/k})$. The description from the K\"unneth formula shows $B_{1,1} = \mathbb{H}$ writing $x_{01}$ and $x_{10}$ to be basis elements of $\left(\H^1(\mathbb{P}^1, \Omega^1_{\mathbb{P}^1/k}) \otimes_k \H^0(\mathbb{P}^1, \Omega^0_{\mathbb{P}^1/k})\right)$ and $\left(\H^0(\mathbb{P}^1, \Omega^0_{\mathbb{P}^1/k}) \otimes_k \H^1(\mathbb{P}^1, \Omega^1_{\mathbb{P}^1/k})\right)$ respectively, we can pick $x_{01}$ and $x_{10}$ such that $B_{1,1}(x_{01}, x_{01}) = B_{1,1}(x_{10},x_{10}) = 0$ and $B_{1,1}(x_{01}, x_{10}) = 1$, and moreover $f^*_{\Hdg}|_{B_{1,1}}$ is given by the matrix $\left(\begin{smallmatrix}
        0 & 1\\ 1 & 0
    \end{smallmatrix}\right)$.
    
    The computation $\det(f^*_\Hdg|_{\H^1(X, \Omega^1)})=-1$ is straightforward, and applying Corollary $\ref{lem:sn of H}$ gives us $\sn(f^*_\Hdg|_{\H^1(X, \Omega^1)})=-2$.
\end{proof}

We now turn our attention to the Bertini involution on del Pezzo surfaces of degree $1$.

\begin{defn}
    Let $S$ be a del Pezzo surface of degree $1$. Since the canonical bundle $\omega_{S/k}$ is anti-ample, $\omega^{\otimes -2}_{S/k}$ defines a double cover $S \to Q$ where $Q \subseteq \mathbb{P}^3$ is a quadric cone. The \emph{Bertini involution} $B_S\colon  S \to S$ is the deck transformation of this double cover.

    When $S$ is given by the blow up of $\mathbb{P}^2$ at $8$ geometric points, the involution $B_S$ has a classical description as a birational transformation of the plane induced by the pencil of cubics passing through those $8$ geometric points, which gives us a description of $B_S$ in the plane Cremona group $\mathrm{Cr}_2(k)$ (see Definition \ref{defn:PlaneCremona}).
\end{defn}

The goal now is to compute $\chikq( [S, B_S])$ by using Theorem \ref{maintheorem:Hermitian}. Using Corollary $\ref{cor:H11Rational}$, it is enough to study the cohomology group $\H^1(S, \Omega^1_{S/k})$. Write $\H^1(S, \Omega^1_{S/k})^-$ to mean the $(-1)$-eigenspace of $B_S^*$ acting on $\H^1(S, \Omega^1_{S/k})$, and $\H^1(S, \Omega^1_{S/k})^+$ to mean the $1$-eigenspace.
\begin{lem}\label{lem:BertiniComp}
    The space $\H^1(S, \Omega^1_{S/k})^+$ is spanned by $K_S$. 
\end{lem}
\begin{proof}
    It is clear that $K_S$ lies in the $1$-eigenspace, since $B_S^*$ is an automorphism so preserves the canonical bundle. It is therefore enough to show that the $(-1)$-eigenspace of the action of $B_S^*$ on $\H^1(S, \Omega^1)$ is $8$-dimensional. Since $B_S^*$ is compatible with base change, it is enough to show this when $k$ is algebraically closed. Since $S$ is smooth and projective, the Hodge--de Rham spectral sequence degenerates on the first page and since $\H^0(S, \Omega^2) = \H^2(S, \Omega^0) = 0$ (see Lemma $\ref{lem:rationalsurface}$), we see that $\H^2_{dR}(S/k) \cong \H^1(S, \Omega^1)$. By applying a comparison isomorphisms between de Rham and singular cohomology, it is enough to show the equivalent statement for $\H^2_{\mr{sing}}(S(\mathbb{C}), \mathbb{Q})$. As in \cite[Proposition 8.8.2]{Dolgachev}, the result follows by the Lefshetz trace formula.
\end{proof}
\begin{cor}\label{cor:bertinidisc}
    The discriminant of $[\H^1(S, \Omega^1)^-]$ is exactly $-\mathrm{disc}(\H^1(S, \Omega^1_{S/k}))$.
\end{cor}
\begin{proof}
    Since $S$ is a del Pezzo surface of degree $1$, $B_{\Hdg}(K_S, K_S) = -1$. This implies that $[\H^1(S, \Omega^1)^+] = \langle -1 \rangle \in \GW(k)$, so the result follows by Corollary $\ref{SpinorInvolution}$.
\end{proof}
\begin{cor}\label{cor:BertiniComputation}
    For $S$ a del Pezzo surface of degree $1$ and $B_S$ its Bertini involution
    \[\chikq([S, B_S]) = (1, \mathrm{disc}(\chi_0(S))).
    \]
\end{cor}
\begin{proof}
     This follows by combining Corollary $\ref{cor:H11Rational}$, Corollary $\ref{SpinorInvolution}$ and Corollary $\ref{cor:bertinidisc}$.
\end{proof}

\begin{rem}\label{rem:computingBertini}
    For a given del Pezzo surface of degree $1$, it is often relatively straightforward to compute $\chi_0(S)$. When $S \cong \mathrm{Bl}_Z\mathbb{P}^2$ where $Z$ is a dimension $0$ subvariety of $\mathbb{P}^2$, we may write $Z = \coprod_{i=1}^r \Spec(L_i)$ for $L_i/k$ some finite extensions where $\sum_{i=1}^r [L_i:k] = 8$. Using the fact that $\chi_0\colon \K_0(\Var_k) \to \GW(k)$ is a ring homomorphism, the blow-up relation in $\K_0(\Var_k)$, and the fact that $\chi_0(\mathbb{P}^2) = \mathbb{H} + \langle 1 \rangle$ (see \cite[Proposition 1.4(5)]{LevineAspects}) we see that $\chi_0(S) = \mathbb{H} + \langle 1 \rangle + \sum_{i=1}^r [\mathrm{Tr}_{L_i/k}]$ and so \[\mathrm{disc}(S) = -\prod_{i=1}^r \mathrm{disc}(L_i/k).\] 
    Alternatively, by \cite[Chapter 3, Theorem 3.5.1]{Kollar}, a del Pezzo surface of degree $1$ is a sextic hypersurface in weighted projective space $\mathbb{P}(1,1,2,3)$, so we may use \cite[Theorem 4.5]{LPLS24} to compute $\chi_0(S)$. 
\end{rem}

We now turn our attention to a different class of involutions: Geiser involutions.

\begin{defn}
    Let $S$ be a del Pezzo surface of degree $2$. The anti-canonical map induced by $\omega^{-1}_S$ gives a surjective, degree $2$ map $\pi\colon  S \to \mathbb{P}^2$, branched over a plane quartic curve. Since $S \to \mathbb{P}^2$ is a double cover, it furnishes an involution $G_S\colon  S \to S$ the \emph{Geiser involution} associated to $S$. When $S$ is rational, the Geiser involution also admits a classical description in the plane Cremona group.
\end{defn}
\begin{lem}\label{lem:GeiserComputation}
    We can compute
    \[
    \chikq( [S, G_S]) = (-1, 2 \cdot \mathrm{disc}(\chi_0(S))) \in \{\pm1\} \times k^{\times}/k^{\times 2}.
    \]
\end{lem}

\begin{proof}
    Corollary $\ref{SpinorInvolution}$ and Corollary $\ref{cor:H11Rational}$ reduces the statement to computing the $1$-eigenspace of the $G_S$ action on $\H^1(S, \Omega^1_{S/k})$.
    
    Let $\mathcal{L}$ be a line bundle on $S$. If $\mathcal{L}=\pi^*\mathcal{L}'$ for some line bundle $\mathcal{L}'$ on $\mathbb{P}^2$, then $\mathcal{L}$ is in the $1$-eigenspace of $G_S^*$. Since $S \to \mathbb{P}^2$ is a double cover, this is an if and only if.
    
    Note that $c_1(\mathcal{O}_{\mathbb{P}^2}(1))$ is the hyperplane class and so $\H^1(\mathbb{P}^2, \Omega^1) \to \H^1(S, \Omega^1)$ is injective. Moreover, $K_S = \pi^*\mathcal{O}_{\mathbb{P}^2(1)}$ by construction of the double cover. Therefore the $1$-eigenspace $\H^1(S, \Omega^1_{S/k})^+$ is $1$-dimensional and generated by the Chern class of canonical bundle $K_S$.

    Since $S$ is a del Pezzo surface of degree $2$, $K_S \cdot K_S = -2$ so $[\H^1(S, \Omega^1_{S/k})^+] = \langle -2 \rangle$. The result follows from here by an identical argument to Corollary $\ref{cor:BertiniComputation}$.
\end{proof}
\begin{rem}
    As in Remark $\ref{rem:computingBertini}$, there are many different methods of computing $\chi_0(S)$. If $S$ is the blow up of $\mathbb{P}^2$ at $Z = \coprod_{i=1}^r \Spec(L_i)$ where $\sum_{i=1}^r [L_i:k]$=7, then we see that $\mathrm{disc}(\chi_0(S)) = - \prod_{i=1}^r \mathrm{disc}(L_i/k)$. Alternatively, \cite[Chapter III, Theorem 3.5.2]{Kollar} exhibit $S$ as a quartic hypersurface in weighted projective space $\mathbb{P}(1,1,1,2)$, so \cite[Theorem 4.5]{LPLS24} computes $\chi_0(S)$. Finally, we could view $S$ as a ramified double cover of $\mathbb{P}^2$ and use \cite[Corollary 2.29]{Broring} to compute $\chi_0(S)$.
\end{rem}

Finally, we turn our attention to another class of involution. Let $\pi\colon  S \to \mathbb{P}^1$ be a rational conic bundle. Suppose that $f_S\colon  S \to S$ is an involution such that
\begin{enumerate}
    \item There exists an involution $f_{\mathbb{P}^1}\colon  \mathbb{P}^1 \to \mathbb{P}^1$ with $\pi \circ f_S = f_{\mathbb{P}^1} \circ \pi$.
    \item The Picard group $\mathrm{Pic}(S_{\kbar})^{f_S^*}$ has rank $2$.
\end{enumerate}

 Such surfaces appear later as $\Z/2\Z$-Mori fibre spaces (see Definition $\ref{defn:MFS}$). Moreover, conic bundle surfaces as this appear as models for involutions known as \emph{de Jonqui\'eres involutions}, see \cite[\S6]{CMYZ} for a treatment over $\mathbb{R}$.

\begin{lem}\label{lem:SNConicBundle}
    The $1$-eigenspace of the $f^*_{\Hdg}$ action, $\H^1(S, \Omega^1_{S/k})^+$, is generated by $K_S$ and $c_1(\pi^*(\mathcal{O}_{\mathbb{P}^1}(1)))$. Moreover, $\mathrm{disc} [ \H^1(S, \Omega^1_{S/k})^+] = 1$.
\end{lem}

\begin{proof}
    We first claim that $\H^1(S, \Omega^1_{S/k})^+$ is $2$-dimensional. This follows by a standard base-field change argument (analogous to \cite[Theorem 2.14]{YauZaslow}), which we recall here for convenience of the reader. Firstly, this dimension is stable under base change, i.e., if $L/k$ is a field extension then $\mathrm{dim}_k \H^1(S, \Omega^1_{S/k})^+ = \mathrm{dim}_L \H^1(S_L, \Omega^1_{S_L/L})^+$. Since $S$ is a variety, it is finite type over $k$ and similarly $f_S\colon S \to S$ is also finite type over $k$. Therefore the pair $(S,i)$ is defined over a finitely generated field so assume $k$ is finitely generated over $\mathbb{Q}$. Embed $k$ into $\kbar$, and since this dimension is stable under base change, it is enough to show $\mathrm{dim}_{\kbar} \H^1(S_{\kbar}, \Omega^1_{S/\kbar})^{f_S^*} = 2$. Since $S$ is a rational surface, $\H^0(S, \Omega^2_{S/\kbar}) = 0 = \H^2(S, \Omega^0_{S/\kbar})$. By \cite[Exercise 5.1.8]{Hartshorne}, the first Chern class gives an injection $c_1\colon \mr{Pic}(S)\otimes_{\mb{Z}}\kbar\to H^1(S,\Omega^1_{S/\kbar})$, and a straightforward Riemann Roch argument shows that the dimensions are equal giving an isomorphism. So, the result holds by the assumption that $\mathrm{Pic}(S_{\kbar})^{f_S^*}$ has rank $2$.

    It only remains to study the bilinear form on $\H^1(S, \Omega^1_{S/k})^+$. Note that both $K_S$ and $c_1(\pi^*(\mathcal{O}_{\mathbb{P}^1}(1)))$ are contained in the $1$-eigenspace of the involution, since clearly the pullback $f_S^*K_S \cong K_S$ and similarly $f_S^*\pi^*(\mathcal{O}_{\mathbb{P}^1}(1))=\pi^*(\mathcal{O}_{\mathbb{P}^1}(1))$. The result holds as long as $K_S$ and $c_1(\pi^*(\mathcal{O}_{\mathbb{P}^1}(1)))$ are linearly independent. We show this by computing the bilinear form on $\mathrm{span}\{K_S, c_1(\mathcal{O}_{\mathbb{P}^1}(1)\}$, which also allows us to recover the result on the discriminant. By a slight abuse of notation, let $K_S$ be the divisor associated to the canonical bundle, and let $F$ be the divisor associated to $\pi^*(\mathcal{O}_{\mathbb{P}^1}(1))$. The bilinear form on $\H^1(S, \Omega^1_S)$ is given by the extension of the intersection pairing, so we compute the intersection form on $\langle K_S, F \rangle$.

    Note that $F \cdot F = 0$, since the conic bundle fibres do not intersect each other. Moreover $K_S \cdot F = \mathrm{deg}_F(\mathcal{O}_F(1))$ by the adjunction formula. Since $F$ is a conic, $\mathrm{deg}_F(\mathcal{O}_F(1))=-2$. This shows that $K_S$ and $F$ must be linearly independent, since $F$ is non-zero, and $K_S \cdot F$ is not a multiple of $F \cdot F$. Moreover, since $\pi\colon  S \to \mathbb{P}^1$ is a conic bundle, we see that the self intersection $K_S \cdot K_S = 8-\#\{\text{singular fibres}\}=N_S$. We therefore see that the bilinear form on $\H^1(S, \Omega^1_S)^+$ is given as a matrix by $\begin{pmatrix} 0 & -2 \\ -2 & N_S\end{pmatrix}$. In particular, $\mathrm{disc}(\H^1(S, \Omega^1_S)^+) = 1$.
\end{proof}
\begin{cor}\label{cor:conicbundlesn}
    We can compute $
    \chikq([S,f_S]) = \left((-1)^{\mathrm{rank}(\chi_0(S))},\mathrm{disc}(\chi_0(S)) \right)$.
\end{cor}
\begin{proof}
    Theorem \ref{maintheorem:Hermitian} shows that we can compute this as the determinant and the spinor norm of $f_S^*$ acting on the Hodge cohomology. Corollary $\ref{cor:H11Rational}$ reduces this to a computation of the action of $f_S^*$ on $\H^1(S, \Omega^1_{S/k})$ and Lemma $\ref{SpinorInvolution}$ and Lemma $\ref{lem:SNConicBundle}$ implies that we have $\mathrm{det}(f^*_S) = (-1)^{\mathrm{dim}( \H^1(S, \Omega^1_{S/k})^-)}$ and $\mathrm{sn}(f^*_S) = \mathrm{disc}( \H^1(S, \Omega^1_{S/k})^-)$.
    
    Since $\H^1(S, \Omega^1_{S/k})^+$ is $2$ dimensional, \[(-1)^{\mathrm{dim}( \H^1(S, \Omega^1_{S/k})^-)} = (-1)^{h^{1,1}(S)-2} = (-1)^{h^{1,1}(S)+2}=(-1)^{\mathrm{rank}(\chi_0(S))},\]
    which gives the determinant part of the calculation. For the spinor norm, we use Corollary $\ref{SpinorInvolution}$ to obtain that $\mathrm{sn}([S,f_S]) = \mathrm{disc}( \H^1(S, \Omega^1_{S/k})^-)$. Lemma $\ref{lem:SNConicBundle}$ gives us an equality
    \[
    \mathrm{disc}(\H^1(S, \Omega^1_{S/k})^-)) = \mathrm{disc}(\H^1(S, \Omega^1_{S/k})^+) \cdot \mathrm{disc}(\H^1(S, \Omega^1_{S/k}))
    \]
    so since $\mathrm{disc}(\H^1(S, \Omega^1_{S/k})^+) = 1$, we have that $\mathrm{disc}(\H^1(S, \Omega^1_{S/k})^-)) = \mathrm{disc}(\H^1(S, \Omega^1_{S/k}))$. Finally, since $\chi_0(S) = \mathbb{H} + [\H^1(S, \Omega^1_{S/k})]$, we have $\mathrm{disc}(\chi_0(S)) = -\mathrm{disc}([\H^1(S, \Omega^1_{S/k})])$.
\end{proof}

\subsection{The behaviour of $\chikq$ under surface blow-ups}
In this subsection, we study the behaviour of the $\chikq$ invariant out of $\K_1(\Var_k)$ under blow-ups. Let $X$ be a smooth projective surface with an automorphism $f_X\colon  X \to X$. Let $Z \subseteq X$ be a closed dimension $0$ subvariety with $f_X(Z) = Z$. Let $\mr{Bl}_Z(X)$ be the blow-up of $X$ along $Z$, and let $\tilde{f_X}$ be the induced automorphism on $\mr{Bl}_Z(X)$. Let $E$ be the exceptional divisor. In this section, we describe how $\chikq([X, f_X])$ relates to $\chikq([Z, f_X|_Z])$ and $\chikq([\mathrm{Bl}_ZX, \tilde{f_X}])$.

\begin{lem}\label{lem:exdiv}
    For the exceptional divisor $E$, we have $E \cong \mathbb{P}^1_k \times_k Z$.
\end{lem}
\begin{proof}
    This is immediate since $Z$ is $0$ dimensional and $E$ is a $\mathbb{P}^1$ bundle over $Z$.
\end{proof}
\begin{rem}
    Lemma $\ref{lem:exdiv}$ may not be true as varieties with a $\Z/2\Z$-action: in particular, even if $f_X|_Z = \mathrm{Id}_Z$, $\tilde{f_X}|_E$ may not necessarily be the identity on $E$.
\end{rem}
\begin{cor}\label{lem:BlowUpdetsn}
Suppose that $f_X|_Z$ is not the identity on $Z$. Let $\chikq( [Z, f_X|_Z]) = (\epsilon, d)$ where $\epsilon \in \{\pm 1\}$ and $d \in k^\times/k^{\times2}$. Then
    \[
\chikq( [\mr{Bl}_Z(X), \tilde{f_X}]) = \chikq([X, f_X]) + ( \epsilon, -d).
    \]
\end{cor}
\begin{proof}
    Using Lemma $\ref{cor:BlowUpRelation}$ and the fact that $\chikq( [Z, f_X|_Z]) = (\epsilon, d)$ it is enough to prove
    \[
\chikq([E, \tilde{f_X}|_E]) = (1,-1).
    \]
    Note that $\tilde{f_X}|_E\colon  E \to E$ is a morphism of $\mathbb{P}^1 \times_k \Spec(A)$ for $A/k$ the finite étale algebra such that $Z \cong \Spec(A)$. Let $\sigma\colon  A \to A$ be the automorphism of $A$ induced by $f_X\colon  Z \to Z$.
    
By construction $\H^0(E, \Omega^0) = A$, and $(\tilde{f_X})^*_{\Hdg}|_{\H^0(E, \Omega^0)} = \sigma$.

Serre duality (\cite[\S3, Corollary 7.13]{Hartshorne}) then implies that $\H^1(E, \Omega^1_{E/k}) \cong A$ and $(\tilde{f_X})^*_{\Hdg}|_{\H^1(E, \Omega^1)} = \sigma$ as well.

Consider $A$ as a vector space over $k$, equipped with its trace bilinear form. Let $\mathbb{H}$ be the hyperbolic form. We see that $\chi^{\Hdg}(E) = A \otimes \mathbb{H}$, as vector spaces with bilinear forms. Under this identification, we see that $(\tilde{f_X})^*_{\Hdg} = \sigma \otimes \mathrm{Id}_{\mathbb{H}}$. Since $\mathrm{dim}(\mathbb{H}) = 2$ and $\mathrm{disc}(\mathbb{H}) = -1$, Lemma $\ref{lem:product}$ gives us $\chikq([E, \tilde{f_X}|_E]) = (1, -1)$ as required.
\end{proof}

Corollary $\ref{lem:BlowUpdetsn}$ shows that, in general, $\chikq$ is far from being a birational invariant of automorphisms since it will in general change when blowing up. However, the following corollary details a situation where the blow-up does not change $\chikq$.

\begin{cor}\label{cor:blowupfixedpoint}
Suppose that $f_X|_Z = \mathrm{Id}_Z$. Then $\chikq( [X,f_X]) = \chikq( [ \mr{Bl}_Z(X), \tilde{f_X}])$.
\end{cor}

\begin{proof}
First note that $\chikq([Z, f_X|_Z]) = (1,1)$ trivially, since $f_X|_Z = \mathrm{Id}_Z$.  The result then follows by an identical argument to Corollary $\ref{lem:BlowUpdetsn}$, however, we see that $\sigma = \mathrm{Id}_A$ so $(\tilde{f_X})^*_{\Hdg}=\mathrm{Id}$, and therefore $\chikq( [E, \tilde{f_X}|_E]) = (1,1)$. 
\end{proof}

\section{Consequences for $\K_1(\Var_k)$}\label{section:K1Vark}

In this short section, we deduce new results about the structure of $\K_1(\Var_k)$. Very little is known about the structure of $\K_1(\Var_k)$. It has not previously been known whether $\K_1(\Var_k)$ can be infinite. In \cite[\S6.3]{CWZ19}, Campbell--Wolfson--Zakharevich construct a surjective homomorphism $\K_1(\Var_k) \to \Z/2\Z$ for $k$ a global or local field with a place of residue cardinality congruent to $3$ mod $4$. This morphism vanishes on the inclusion $\Z/2\Z \cong \K_1(\mathrm{Finite Sets}) \to \K_1(\Var_k)$, which tells us that for such fields, $\K_1(\Var_k)$ has at least $4$ elements. When $k$ is a number field, a theorem of Braunling--Groechenig \cite[Theorem 1.3(1)]{BG21} shows that $\K_1(\mathcal{O}_k)$ is a direct summand of $\pi_1(\K(\Var_k)[A^{-1}])$ for some abelian variety $A$, however, this theorem requires working in the localisation. 

Despite this, the results of the previous section immediately allow us to prove Theorem \ref{maintheorem:K1Var}, which states that if $k$ is a finitely generated field of characteristic $0$ then $\K_1(\Var_k)$ is not finitely generated.

\begin{proof}[Proof of Theorem \ref{maintheorem:K1Var}]
    The map $\chikq\colon  \K_1(\Var_k) \to \{\pm 1\} \times k^\times/k^{\times2}$ is surjective by Corollary $\ref{lem:surjectivityforC}$. If $k$ is a number field, then $k^\times/k^{\times 2}$ is not finitely generated, which gives the result. For the general case, let $K$ denote the algebraic closure of $\mathbb{Q}$ in $k$, which is a number field since $k$ is finitely generated. The injection $K^\times/K^{\times 2} \hookrightarrow k^\times/k^{\times 2}$ shows that $k^\times/k^{\times 2}$ is also not finitely generated as required.
\end{proof}
We in fact prove a slightly stronger result. Consider the category of varieties of dimension $0$, $\Var_k^{\leq 0}$. This category is equivalent to the opposite category of finite étale $k$-algebras via $\Spec(-)$, which is also equivalent to the category of finite $\Gal_k$-sets. This is an assembler (\cite[Page 5]{Annihilator}) so has an associated $\K$-theory spectrum $\K(\Var_k^{\leq 0})$. 
\begin{cor}
    Suppose that $k$ is a finitely generated field of characteristic $0$. Then $\K_1(\Var^{\leq0}_k)$ is not finitely generated.
\end{cor}
\begin{proof}
    As in \cite[page 5]{Annihilator}, the inclusion of categories $\Var_k^{\leq 0} \hookrightarrow \Var_k$ induces a morphism of $\K$-theory spectra $\K(\Var_k^{\leq 0}) \to \K(\Var_k)$. Taking $\pi_1$ and then post-composing with $\chikq$ gives us a homomorphism $\K_1(\Var_k^{\leq 0}) \to \{\pm 1\} \times k^{\times}/k^{\times 2}$. 

    If $L/k$ is a quadratic extension with $\sigma \in \Gal(L/k)$ the non-trivial element, then $[\Spec(L), \sigma] \in \K_1(\Var_k^{\leq 0})$, so the result follows by Corollary \ref{cor:sn of quadratic}. 
\end{proof}

We now turn our attention to a different structure on $\K(\Var_k)$: its $\Einfty$-ring structure.
\begin{rem}\label{rem:Einfty}
For the purpose of this paper, we largely consider $\K(\Var_k)$ and $\End_{\SH(k)}(\1_k)$ just as infinite loop spaces. However, it is worth noting that they have additional structure: both naturally are $\Einfty$-ring spectra. It is an open question, asked independently in private communication to the second author by Zakharevich and to the third author by Wickelgren, whether the map $\chi \colon  \K(\Var_k) \to \End_{\SH(k)}(\1_k)$ can be promoted to a map of $\Einfty$-ring spectra. We provide some evidence in the affirmative. 
\end{rem}
Note that $\K_1(\Var_k)$ is equipped with the structure of a $\K_0(\Var_k)$-module, since $\K(\Var_k)$ is an $\Einfty$-ring spectrum. This $\K_0(\Var_k)$-module structure admits an explicit description. 
\begin{lem}\label{lem:K0modulestructure}
If $Y$ is a variety and $[X,f]\in \K_1(\Var_k)$, then the $\K_0(\Var_k)$-module structure on $\K_1(\Var_k)$ satisfies $[Y]\cdot [X,f]=[X\times Y,f\times\mr{id}_Y]$.
\end{lem}
\begin{proof}
    This follows immediately by \cite[Theorem 1.4]{Campbell}.
\end{proof}
\begin{rem}
   The above uses the SW-category construction of $\K(\Var_k)$,we could instead work with the presentation of $\K_1(\Var_k)$ given to us by assemblers, in which case, this would recover the $\K_0(\Var_k)$-module structure on $\K_1(\Var_k)$ as determined in \cite[Proposition 7.6, Example 7.14, Theorem 8.3]{Zakharevich2022}.
\end{rem}

Similarly, $\KQ$ is an $\Einfty$-ring spectrum, which means that $\KQ_1(k)$ is equipped with the structure of a $\KQ_0(k)=\GW(k)$-module. We can also describe this action explicitly. 
\begin{lem}
    The $\GW(k)$-module structure on $\{\pm 1\} \times k^\times/k^{\times 2}$ is given by
    \[
    q \cdot (a,b) = (a^{\mathrm{rank}(q)}, \mathrm{disc}(q) b^{\mathrm{rank}(q)}).
    \]
\end{lem}
\begin{proof}
    Since $\GW(k)=\KQ_0(k)$ is generated by $1$-dimensional bilinear forms, it is enough to prove the result when $q= \langle \alpha \rangle$ for $\alpha \in k^\times$. 

    The map $\pi_{1,0}(\1_k) \to \KQ_1(k)$ is a $\GW(k)$-module homomorphism, since, as in \cite[Theorem 1.5]{PaninWalter2019}, the unit map $u\colon  \1_k \to \KQ$ is a morphism of $\Einfty$-ring spectra. The result follows by applying \cite[Theorem 7.1]{OrmsbyOstvaer}, which shows that the topological Hopf element $\eta_{top} \in \pi_{1,0}(\1_k)$ satisfies $u( \langle \alpha \rangle \eta_{top}) = (-1, \alpha)$. 
\end{proof}
The map $\chi_0\colon  \K_0(\Var_k) \to \GW(k)$ therefore allows us to endow any $\GW(k)$-module with the structure of a $\K_0(\Var_k)$-module, where $[X] \in \K_0(\Var_k)$ acts through $\chi_0$. 

\begin{thm}\label{thm:K0Module}
    Let $Y/k$ be any smooth projective variety, and let $f\colon  X \to X$ be an automorphism. Then $\chikq( [Y] \cdot [X, f]) = [Y] \cdot \chikq([X,f])$.
\end{thm}
\begin{proof}
    Using the description of the $\K_0(\Var_k)$-module structure on $\K_1(\Var_k)$ in Lemma $\ref{lem:K0modulestructure}$, we see that $ [Y] \cdot [X, f] = [ X \times_k Y, f \times \mathrm{id}]$. The K\"unneth formula for Hodge cohomology allows us to write $\Hdg(X \times_k Y) = \Hdg(X) \otimes_k \Hdg(Y)$, along with their bilinear forms. Under this identification $(f \otimes \id)^*_{\Hdg} = f^*_{\Hdg} \otimes \id$.

       Theorem \ref{maintheorem:Hermitian} implies $\chikq([Y] \cdot [X, f]) = (\mathrm{det} (f^*_{\Hdg} \otimes \id), \sn(f^*_{\Hdg} \otimes \id))$. The result follows by applying Lemma $\ref{lem:product}$ to $f^*_{\Hdg} \otimes \id$.
\end{proof}
Theorem $\ref{thm:K0Module}$ provides evidence that $\chikq$ is a $\K_0(\Var_k)$-module homomorphism and therefore that $\chi$ is indeed a morphism of $\Einfty$-ring spectra: if $\chi$ is a map of $\Einfty$-ring spectra, this automatically implies that the induced morphisms $\K_n(\Var_k) \to \pi_n(\End_{\SH(k)}(\1_k))$ are $\K_0(\Var_k)$-module homomorphisms.

\section{A birational invariant from $\chikq$}\label{section:birational}
Over the complex numbers, the study of conjugacy classes of involutions in the plane Cremona group $\mathrm{Cr}_2(k)$ was initiated by Bertini in 1877 \cite{Bertini}, and later a complete classification was given in \cite{BayleBeauville}. In \cite{DolgachevIskovskikh}, it was asked whether such conjugacy classification problems could be extended to non-algebraically closed fields. Recent work studies the classification problem over $\mb{R}$ in \cite{CMYZ}. Over other fields there are fewer tools available for the problem. This section shows that these types of classification problems have a natural home in $\K_1(\mathrm{Var})_k$, and deduces new results about involutions by using $\chikq$. We first recall the strategy for classification problems as outlined in \cite{BayleBeauville}. We use Dolgachev-Iskovskikh's theory of equivariant Sarkisov links outlined in  \cite{DolgachevIskovskikh}, to show that such links give rise to relations in $\K_1(\mr{Var}_k)$. Then we show that these relations are in the kernel of $\chikq$ under certain conditions. This gives a new and computable conjugacy invariant of birational involutions in $\mr{Cr}_2^{\mr{ab}}(k)$.
\subsection{Definitions in birational geometry}

\begin{defn}
    Let $X,Y$ be $k$ varieties. A \emph{rational map} $\phi\colon X \dashrightarrow Y$ is given by a Zariski dense open $U\subset X$ and a regular morphism $\phi \colon U\to Y$. A rational map is \emph{birational} if the image of $\phi\colon U\to Y$ is a Zariski dense open set $V\subset Y$, and $\phi\colon U \to V$ is an isomorphism.
    
    A \emph{birational morphism} is a morphism $g\colon X\to Y$ defined on all of $X$ that defines a birational map.
\end{defn}

\begin{defn}\label{defn:PlaneCremona}
    Let $Y$ be a variety over $k$. Let $\mr{Bir}(Y)$ denote the group of birational maps $\phi\colon Y\dashrightarrow Y$. When $Y=\mb{P}^n_k$, we write
    \[
    \mr{Cr}_n(k):=\mr{Bir}(\mb{P}^n_k)
    \]
    and call $\mr{Cr}_n(k)$ the \emph{Cremona group on $n$-letters} over $k$.  When $n=2$, $\mr{Cr}_2(k)$ is also called the \emph{plane Cremona group}.
\end{defn}
Our method for studying conjugacy classes of elements of $\mr{Bir}(Y)$ with a fixed order heavily relies on the fact that finite order elements of $\mr{Bir}(Y)$ are \emph{regularisable}.
\begin{defn}\label{defn:regularisable}
   We say that $f \in \mathrm{Bir}(Y)$ is \emph{regularisable} if there exists a variety $X$, a birational map $\phi \colon X\dashrightarrow Y$ and an automorphism $f_X\colon X \to X$ such that $\phi \circ f_X\circ\phi^{-1}=f$. We say that $(X,f_X)$ \emph{regularises} $f$.

    If $Y=\mathbb{P}^n$, we say $f \in \mathrm{Bir}({\mathbb{P}^n})$ is \emph{linearisable} if it is regularisable and there is a regular morphism $\mathbb{P}^n \to \mathbb{P}^n$ which regularises $f$. We say $f$ is \emph{geometrically linearisable} if $f$ is linearisable after base changing to $\kbar$.  
\end{defn}
\begin{thm}[\cite{Cheltsov} Thm.~1.7]\label{thm:regularize}
    Let $Y$ be a variety over $k$. Then finite order elements of $\mr{Bir}(Y)$ are regularisable.
\end{thm}
For the remainder of this section, fix a positive integer $n$.
\begin{defn}
    A (rational) \emph{$\Z/n\Z$-variety} is a pair $(X,f_X)$ where $X$ is a (rational) variety over $k$ and $f_X$ is a homomorphism $f_X\colon\Z/n\Z\to \Aut(X)$. Note that $f_X$ is uniquely determined by $f_X(1)$, so we will instead think of $f_X \in \Aut(X)[n]$. As in Remark $\ref{rem:equivariantrov}$, a $\Z/n\Z$-variety $(X, f_X)$ determines a class $[X, f_X] \in \K_1(\Var_k)$. A \emph{$\Z/n\Z$-morphism} of $\Z/n\Z$-varieties, $g\colon (X,f_X) \to (Y,f_Y)$ is a morphism of varieties $g\colon X \to Y$ such that $g \circ f_X = f_Y \circ g$. Similarly, a \emph{$\Z/n\Z$-rational map} is a rational map $\phi \colon X \dashrightarrow Y$ such that $\phi \circ f_X = f_Y \circ \phi$.

\end{defn}

The following theorem of Prokhorov reframes the study of conjugacy classes in $\mr{Bir}(X)$ to the study of $\mb{Z}/n\mb{Z}$-varieties rational to $X$.
\begin{thm}[Corollary of Proposition 2.1 of \cite{Prokhorov2}]
There is a bijective correspondence between conjugacy classes of order $n$ elements of $\mr{Bir}(Y)$, and $\mb{Z}/n\mb{Z}$-varieties birational to $Y$ up to $\mb{Z}/n\mb{Z}$-birational equivalence.
\end{thm}
\begin{defn}
    A $\mb{Z}/n\mb{Z}$-variety $(X,f_X)$ is \emph{$\Z/n\Z$-minimal} if, for every $(X', f_{X'})$ with a $\Z/n\Z$-birational morphism $f\colon (X, f_X) \to (X', f_{X'})$, we have that $f$ is an isomorphism.
\end{defn}

\begin{thm}\cite[Theorem 2.6]{ESS25}
Every $\mb{Z}/n\mb{Z}$-surface is $\Z/n\Z$-birational to a (not necessarily unique) $\mb{Z}/n\mb{Z}$-minimal surface.
\end{thm}
\begin{rem}\label{rem:minimalsurface}
    If $(X,f_X)$ and $(X',f_{X'})$ are rational $\mb{Z}/n\mb{Z}$-minimal surfaces (i.e. $X$ and $X'$ are birational to $\mb{P}^2$), there may still exist a $\mb{Z}/n\mb{Z}$-birational map $(X,f_X)\dashrightarrow (X,f_{X'})$ (i.e. defined only on a Zariski open) that cannot be extended to an isomorphism. For example, consider the map $\phi\colon  \mathbb{P}^1 \times \mathbb{P}^1 \dashrightarrow \mathbb{P}^2$ defined by $( [a,b], [c,d]) \mapsto [ad, bc, bd]$. Even with a trivial $\Z/n\Z$ action, both $\mathbb{P}^1 \times \mathbb{P}^1$ and $\mathbb{P}^2$ are $\Z/n\Z$-minimal surfaces, but this birational map does not extend to an isomorphism.
\end{rem}

\begin{defn}
    Let $Y$ be a $k$-variety, and $\mathrm{BirEq}_{\mb{Z}/n\mb{Z}}^{\mr{min}}(Y)$ denote the subgroup of $\K_1(\mathrm{Var}_k)$ generated by symbols $[(X, f_X)]$ where $(X,f_X)$ is a $\Z/n\Z$-minimal variety and $X$ is birational to $Y$. Write 
    \[
    \mathrm{Rat}_{\mb{Z}/n\mb{Z}}^{d,\mr{min}}:=\mathrm{BirEq}_{\mb{Z}/n\mb{Z}}^{\mr{min}}(\mb{P}^d).
    \]
    
\end{defn}
Observe that the \emph{group} $\mathrm{BirEq}_{\mb{Z}/n\mb{Z}}^{\mr{min}}(Y)$ is the refinement of the \emph{set} of $\mb{Z}/n\mb{Z}$-minimal varieties. While the difference between these two objects is significant, a natural way of producing an invariant of $\mr{BirEq}_{\mb{Z}/n\mb{Z}}(Y)$ is to construct a group homomorphism out of $\mathrm{BirEq}_{\mb{Z}/n\mb{Z}}^{\mr{min}}(Y)$ and check that the values are constant on $\mb{Z}/n\mb{Z}$-birational classes.

\subsection{Relationship to $\K_1(\Var_k)$}
In this subsection, we explain how the structure of $\K_1(\mr{Var}_k)$ provides a natural connection between the involution classification problem, i.e. the study of $\mr{Rat}^{2,\mr{min}}_{\mb{Z}/n\mb{Z}}$, and $\mr{Cr}_2^{\mr{ab}}(k)$.  We show that we can relate $\mr{Rat}^{2,\mr{min}}_{\mb{Z}/n\mb{Z}}$ and $\mr{Cr}_2^{\mr{ab}}(k)$ by a natural short exact sequence in Corollary $\ref{cor:K1Rat}$. This short exact sequence is constructed as follows.

As in \cite[Page 5]{Annihilator}, the category $\mathrm{Var}^{\leq n}_k$ of varieties of dimension $\leq n$ can be given the structure of an assembler, so has an associated $\K$-theory spectrum $\K(\mathrm{Var}^{\leq n})$. The inclusion $\mathrm{Var}^{\leq n}\hookrightarrow \mathrm{Var}^{\leq n+1}$ furnishes a long exact sequence in homotopy groups
    \begin{equation}\label{eq:lesK(Var)}
    \dots \to \K_m(\mathrm{Var}_k^{\leq n})\to \K_m(\mathrm{Var}_k^{\leq n+1})\to \pi_m(\K(\mathrm{Var}_k^{\leq {n+1}}),\K(\mathrm{Var}_k^{\leq n})) \to \dots  .
    \end{equation}
    We will write $\partial_m:\pi_m(\K(\mathrm{Var}_k^{\leq n+1}),\K(\mathrm{Var}_k^{\leq n})) \to \K_{m-1}(\Var_k^{\leq n})$ for the boundary map of the above long exact sequence. When $m=1$, Theorem A of \cite{Annihilator} gives 
    \[\pi_1(\K(\mathrm{Var}_k^{\leq n}),\K(\mathrm{Var}_k^{\leq n-1}))=\bigoplus_{[X]\in\mathrm{Bir}_n(k)}\mathrm{Aut}_k(k(X))^{\mathrm{ab}}\]
where $\mathrm{Bir_n(k)}$ is the set of birational isomorphism classes of $n$-dimensional $k$-varieties. In particular, consider the case when $n=2$. Then $\Aut_k(k(\mathbb{P}^2))^{ab} = \mathrm{Cr}_2^{\mathrm{ab}}$ is a subgroup of $\pi_1(\K(\Var_k^{\leq 2}), \K(\Var_k^{\leq 1}))$, so $\partial_1$ restricts to give a homomorphism $\mathrm{Cr}_2^{\mathrm{ab}}(k) \to \K_0(\Var_k^{\leq 1})$. This homomorphism is denoted $\tilde{c}$ in \cite{LinShinder}.

Let $X$ be  an $n$-dimensional $k$-variety. By \cite[Lemma 3.2]{Annihilator} we have a formula for $\partial_1$. For $\phi\in\mathrm{Aut}(k(X))$, pick $U$ and $V$, open subvarieties of $X$, such that $\phi$ defines an isomorphism $U\xrightarrow{\cong} V$ (such $U,V$ are not unique). Then \cite[Lemma 3.2]{Annihilator} shows
\[
\partial_1(\phi)=[X\setminus U]-[X\setminus V]\in \K_0(\mathrm{Var}^{\leq n-1}).
\]
A large class of birational automorphisms are killed by $\partial_1$.

\begin{prop}\label{prop:vanishingboundary}
    Let $X$ be an $n$-dimensional variety over $k$. Let $f \in \mathrm{Bir}(X)$ and suppose that $f$ is regularisable. Then $\partial_1(f)=0$.\begin{footnote}{More generally, $\partial_1$ vanishes on pseudo-regularisable elements (see \cite{LinShinder}), but this is not necessary for the purposes of this paper}\end{footnote}
\end{prop}
\begin{proof}
   Since $f$ is regularisable, there exists a regular automorphism $\tilde{f}\colon \tilde{X}\to \tilde{X}$ and a birational map $\phi\colon X\dashrightarrow \tilde{X}$ such that $\phi^{-1}\circ\tilde{f}\circ\phi=f$. As a result, $\tilde{f}$ is identified with $f$ in $\bigoplus_{[Y]\in\mathrm{Bir}_n(k)}\mathrm{Aut}_k(k(Y))^{\mathrm{ab}}$, and since $\tilde{f}$ is an automorphism, $\partial_1(\tilde{f})=0$.
\end{proof}

The boundary map $\partial_1$ was used in \cite{LinShinder} to answer questions of Dolgachev regarding properties of generators in $\mathrm{Cr}_n(k)$, where they produced elements of $\mathrm{Cr}_n(k)$ for $n>2$ and for various fields $k$ on which $\partial_1\neq 0$. Our work is focused in dimension $2$ and over fields of characteristic $0$. In \cite{LamySchneider}, Lamy and Schneider prove that the plane Cremona group is generated by involutions over a perfect field.  

\begin{cor}\label{cor:cremona2vanish}
    The plane Cremona group $\mathrm{Cr}_2^{\mathrm{ab}}(k)$ lies in the kernel of the boundary map $\partial_1\colon \bigoplus_{[X] \in \mathrm{Bir}_2(k)} \Aut_k(k(X))^{\mathrm{ab}} \to \K_0(\Var_k^{\leq 1})$.
\end{cor}
\begin{proof}
    By \cite[Theorem 1.1]{LamySchneider}, the plane Cremona group $\mathrm{Cr}_2(k)$ is generated by involutions, so the same is true for its abelianisation. Theorem $\ref{thm:regularize}$ implies that every involution in $\mathrm{Cr}_2(k)$ is regularisable, so Proposition $\ref{prop:vanishingboundary}$ gives the result.
\end{proof}
\begin{cor}\label{cor:K1Rat}
    The natural map $\mathrm{Rat}_{\mb{Z}/n\mb{Z}}^{2,\mr{min}} \to \mathrm{Cr}_2^{ab}(k)$  is surjective.
\end{cor}
\begin{proof}
    This follows by the proof of Corollary $\ref{cor:cremona2vanish}$, but instead picking a $\Z/2\Z$-minimal variety to regularise the generators of $\mathrm{Cr}_2(k)$, which is possible by continuously blowing down if our variety is not $\Z/2\Z$-minimal.
\end{proof}
Write $\Psi_1$ to denote the kernel of the map in Corollary $\ref{cor:K1Rat}$. By construction, $\mathrm{Rat}_{\mb{Z}/n\mb{Z}}^{2,\mr{min}}$ is a subgroup of $\K_1(\Var_k^{\leq 2})$, and so $\Psi_1 \subseteq \K_1(\Var_k^{\leq 1})$.

\begin{prop}\label{prop:relationtypes}
    The group $\Psi_1$ is generated by elements of one of the following forms.
    \begin{enumerate}
        \item Differences  $[X,f_X] - [X',f_{X'}]$ where $(X,f_X)$ is $\mb{Z}/2\mb{Z}$-birational to $(X', f_{X'})$.
        \item Sums of the form $\sum_{j=1}^m [X_j, f_j]$, where each $f_j$ is an automorphism, and if we write $\sigma_j$ to denote the image of $f_j$ in the Cremona group, then each $\sigma_1, \ldots, \sigma_m$ are all pairwise non-conjugate and $\sigma_1 \circ \ldots \circ \sigma_m = \mathrm{Id} \in \mathrm{Cr}_2(k)$.
    \end{enumerate}
\end{prop}
\begin{proof}
The fact that we can assume each $f_j$ is an involution is a result of \cite[Theorem 1.1]{LamySchneider}, so we only need to show these are the relations.

It is clear that relations $(1)$ and $(2)$ are contained in $\Psi_1$. Suppose $\sum_{j}[X_j, f_j]$ lies in $\Psi_1$ for some finite collection $[X_j,f_j]\in\mr{Rat}^{\mr{min}}_{\mb{Z}/2\mb{Z}}$. If there are two elements $[X_{j_1},f_{j_1}],[X_{j_2}, f_{j_2}]$ mapping to the same involution in $\mr{Cr}_2^{\mr{ab}}$, then $[X_{j_1},f_{j_1}]-[X_{j_2}, f_{j_2}]=0$. So assume that all of the $[X_j,f_j]$ terms appearing in the sum regularise different involutions in $\mr{Cr}_2^{\mr{ab}}$. The sum therefore ends up giving us a relation of type (2).
\end{proof}
\begin{rem}
While $\Psi_1$ is presented by differences of elements in $\K_1(\mr{Var}_k^{\leq 2})$, it is naturally a subgroup of $\K_1(\mr{Var}_k^{\leq 1})$. Consider a relation of type $(1)$, i.e. an element $[X, f_X]-[X', f_{X'}]$ where $f_X, f_{X'}$ regularise the same birational map $f\colon \mb{P}^2\dashrightarrow\mb{P}^2$. Then there exist dense open subsets $U\subset X,V\subset X'$ with $U\cong V$ and $[U,f_X|_U]=[V,f_{X'}|_{V}]$. As a result, 
\[
[X, f_X]-[X', f_{X'}]=[X \setminus U,f_{\tilde{X}\setminus U}]-[X' \setminus V, f_{X' \setminus V}]\in \K_1(\mr{Var}^{\leq 1}).
\]
\end{rem}

\begin{rem}\label{rem:CremonaInvariants}
The relations in Proposition $\ref{prop:relationtypes}$ come in two different types: type (1) relations encode birational invariants, whereas type (2) relations encode relations in $\mathrm{Cr}_2^{ab}(k)$.

    Suppose then that $\Phi$ is a homomorphism $\mathrm{Rat}_{\mb{Z}/n\mb{Z}}^{d,\mr{min}} \to G$ where $G$ is an abelian group. The morphism $\Phi$ factors through $\mathrm{Cr}_2^{ab}(k)$ if and only if $\Phi(\Psi_1) = 0$.

    Suppose that $\Phi$ only vanishes on the relations of type (1) from Proposition $\ref{prop:relationtypes}$. Then for any $f \in \mathrm{Cr}_2^{ab}(k)$ we see that if $(X,f_X)$ and $(X', f_{X'})$ are $\Z/2\Z$-minimal varieties regularising $f$, we have $\Phi([X, f_X]) = \Phi([X', f_{X'}])$. In particular, if $\Phi$ vanishes on relations of type $(1)$, it gives us a birational invariant of Cremona group elements.  
\end{rem}

\subsection{Invariants of involutions in $\mathrm{Cr}_2(k)$}
In this section, we use the previous subsection to show that the connection between the group $\mathrm{BirEq}_{\mb{Z}/n\mb{Z}}^{\mr{min}}(Y)$ and the set of conjugacy classes of order $n$ elements of $\mr{Bir}(Y)$ is already interesting in the case $n=2$, and $Y=\mb{P}^2$. The $\mb{Z}/2\mb{Z}$-minimal model programme for rational surfaces, initiated by Mori, states that any birational map between $\mb{Z}/2\mb{Z}$-minimal rational surfaces $(X,f_X)\dashrightarrow (X',f_{X'})$ factors into $\mb{Z}/2\mb{Z}$-Sarkisov links of four different types (see Definition \ref{def:linktypes} and Theorem \ref{thm:sarklink}).

In Proposition $\ref{propn:linkeqns}$, we show that $\Z/2\Z$-Sarkisov links give rise to relations between $[(X,f_X)]$ and $[(Y, f_Y)]$ in $\mathrm{Rat}_{\mb{Z}/n\mb{Z}}^{2,\mr{min}}$. An extensive list of all possible Sarkisov links between different rational surfaces can be found in \cite[Theorem 2.6]{Iskovskikh}, allowing us to enumerate all possible relations coming from $\Z/2\Z$-Sarkisov links. This allows us to prove Theorem \ref{maintheorem:birational}, which shows that $\chikq$ descends to an invariant on conjugacy classes geometrically non-linearisable involutions in $\mr{Cr}_2(k)^{\mr{ab}}$. 

Over an algebraically closed field, we have the following classification of involutions of the plane Cremona group.
\begin{thm}\cite[Theorem 1.4]{BayleBeauville}\label{thm:listofminimals}
    Let $(X,f_X)$ be a $\mb{Z}/2\mb{Z}$-minimal rational surface over an algebraically closed field. Then one of the following hold.
    \begin{enumerate}
        \item There exists a smooth $\mb{P}^1$ fibre bundle $\pi\colon X\to\mb{P}^1$ and a non-trivial involution $\tau$ of $\mb{P}^1$ such that $\pi\circ f_X=\tau\circ \pi$.
        \item There exists a fibration $\pi\colon X\to\mb{P}^1$ such that $\pi\circ f_X=\pi$; the smooth fibres of $\pi$ are rational curves on which $f_X$ induces a non-trivial involution, and any singular fibre is the union of two rational curves exchanged by $f_X$ meeting at one point.
        \item $X\cong\mb{P}^2$.
        \item $X\cong\mb{P}^1\times\mb{P}^1$ and under this identification $f_X(x,y)=(y,x)$.
        \item $X$ is a degree 2 del Pezzo surface and $f_X$ is its Geiser involution.
        \item $X$ is a degree 1 del Pezzo surface and $f_X$ is its Bertini involution.
    \end{enumerate}
\end{thm}
\begin{defn}\label{defn:Fixf}
    Following \cite[(2.2)]{BayleBeauville}, if $f$ is an involution in $\mr{Cr}_2(k)$, write $\mr{Fix}(f)$ to be the union of the normalisations of the non-rational curves fixed by $f$ in $\mathbb{P}^2$.
\end{defn}

\begin{thm}[Theorem 2.6-2.7 of \cite{BayleBeauville}]\label{thm:threegeomtypes}
    Let $\kbar$ be an algebraically closed field and let $f\in\mr{Cr}_2(\kbar)$ be an involution. Then the normalisation of $\mr{Fix}(f)$ determines the conjugacy class of $i$. Moreover $f$ is conjugate to exactly one of the following types of involutions:
    \begin{itemize}
        \item a de Jonqui\`{e}res involution,
        \item a Bertini involution,
        \item a Geiser involution.
    \end{itemize}
\end{thm}
Over $\mb{R}$, \cite[Main Corollary]{CMYZ} shows that $\mr{Fix}(f)$ does \emph{not} determine the conjugacy class of $f$. Theorem \ref{maintheorem:birational} later shows that the same is true for any characteristic $0$ field admitting a non-trivial quadratic extension.

\begin{lem}\label{rem:deg2linear}
Let $\kbar$ be an algebraically closed field and let $f \in \mathrm{Cr}_2(\kbar)[2]$. Suppose that $f$ is linearisable. Then $f$ is conjugate to a de Jonqui\`eres involution.
\end{lem}
\begin{proof}
    Note that a linearisable involution fixes a curve of degree $2$ in $\mathbb{P}^2$, so in particular $\mathrm{Fix}(f)$ is empty (since we omit rational curves from $\mathrm{Fix}(f)$). The result follows since the Bertini involution on a del Pezzo surface of degree $1$ fixes a curve of genus $4$ by \cite[Example 2.4(b)]{BayleBeauville} and the Geiser involution on a del Pezzo surface of degree $2$ fixes a curve of genus $3$ by \cite[Example 2.4(a)]{BayleBeauville}.
\end{proof}

The following theorem follows from the proof of Theorems 1.4, 2.6 and 2.7 of \cite{BayleBeauville}.
\begin{thm}\label{thm:algclosedlinear}
    Suppose $(X,f_X)$ is a $\mb{Z}/2\mb{Z}$ rational surface over an algebraically closed field. Then $(X,f_X)$ regularises a linearisable involution if and only if $K_X^2>4$. 
\end{thm}
\begin{proof}
    We exhaust the cases enumerated in Theorem \ref{thm:listofminimals}, leaving case (2) until the end. 
    \begin{itemize}
    \item If $(X,f_X)$ is in case (1). Then since all fibres of $X \to \mathbb{P}^1$ are smooth, $K_X^2=8$, and the proof of \cite[Theorem 2.6-2.7]{BayleBeauville} shows that $(X,f_X)$ is $\mb{Z}/2\mb{Z}$-birational to a linear involution on $\mb{P}^2$. 
    \item If $(X,f_X)$ is in case (3), so that $X=\mb{P}^2$. Therefore $K_{\mb{P}^2}^2=9$ and $f_X$ is linear.
    \item If $(X,f_X)$ is in case (4) then $K_X^2=8$, and the proof of \cite[Theorem 2.6-2.7]{BayleBeauville} shows that $(X,f_X)$ is $\mb{Z}/2\mb{Z}$-birational to a linear involution.
    \item Suppose that $(X,f_X)$ is in case (5) or case (6) respectively. Then $K_X^2 = 2$ or $1$ respectively by definition. Therefore $K_X^2\leq 4$ and Theorem \ref{thm:threegeomtypes} shows that $(X,f_X)$ is not linearisable.
    \item If $(X,f_X)$ is in case (2), this breaks into two sub-cases:
    \begin{itemize}
        \item If $\pi\colon X\to\mb{P}^1$ is smooth, the proof of Theorem \ref{thm:threegeomtypes} shows that $(X,f_X)$ is $\mb{Z}/2\mb{Z}$ birational to a surface in case (4), so this reduces to the earlier resolution of case (4).
        \item If $\pi$ is not smooth \cite[Complement 1.5]{BayleBeauville} shows that $\pi$ has $2g+2$ singular fibres. By \cite[Lemma 3.11]{Prokhorov}, we see $K_X^2=8-(2g+2)$. For $g=0$, the fixed curve of the associated birational involution in $\mb{P}^2$ is a smooth conic. By \cite[Proposition 2.7]{BayleBeauville}, this is a degree $2$ de Jonqui\'{e}res involution, which is linearisable by Lemma \ref{rem:deg2linear}. For $g\geq 1$, $(X,f_X)$ is not linearisable as $f_X$ fixes an irrational curve, and $K_X^2 = 8-(2g+2) \leq 4$. This proves the result.\qedhere
    \end{itemize} 
    \end{itemize} 
\end{proof}
\begin{cor}\label{cor:nonclosedlinear}
Let $f$ be an involution in $\mr{Cr}_2(k)[2]$. If $f$ is geometrically non-linearisable, then all $\mb{Z}/2\mb{Z}$-minimal surfaces $(X,f_X)$ regularising $f$ satisfy $K_X^2\leq 4$. 
\end{cor}
\begin{proof}
    Let $(X,f_X)$ be a $\mb{Z}/2\mb{Z}$-minimal surface over $k$ which regularises $f$. Let $(\overline{X},f_{\overline{X}})$ be the base change to $\kbar$, and note $K_X^2=K_{\overline{X}}^2$. If $(\overline{X},\overline{f_X})$ is $\mb{Z}/2\mb{Z}$-minimal, $f$ not being linearisable implies that $K_{\overline{X}}^2\leq 4$ by Theorem \ref{thm:algclosedlinear}. However, $(\overline{X},\overline{f_X})$ may not be $\Z/2\Z$-minimal over $\kbar$. In this case, by \cite[Lemma 1.1]{BayleBeauville}, there exists $E \subseteq \overline{X}$ where $E$ is a smooth rational $(-1)$-curve and either $\overline{f_X}(E) = E$ or $\overline{f_X}(E) \cap E = \emptyset$. In either case we can blow down the orbit of $\overline{f_X}(E)$ to get another rational $\mb{Z}/2\mb{Z}$-surface $(\overline{X_1}, \overline{f_{X_1}})$. Since $\overline{X_1}$ is a blow-down of $\overline{X}$, $K_{\overline{X}}^2<K_{\overline{X_1}}^2$. One can continue in this way until we have a $\Z/2\Z$-minimal surface $(\overline{X}_{\min}, \overline{f_{X_{min}}})$, so that there are no more smooth rational (-1)-curves $E \subseteq \overline{X}_{\min}$ satisfying $\overline{f_{X_{min}}}(E)=E$ or $\overline{f_{X_{min}}}(E)\cap E=\emptyset$.
    
    Since $\overline{f}$ is non-linearisable, $K_{\overline{X}}^2<K_{\overline{X_{\mr{min}}}}^2\leq 4$ by Theorem \ref{thm:algclosedlinear}, as required.
\end{proof}
\begin{rem}
    Corollary \ref{cor:nonclosedlinear} does not say anything whether there exist geometrically linearisable involutions over non-algebraically closed fields which are regularised by a surface $X$ with $K_{X}^2\leq 4$. It would be interesting to enumerate such surfaces to better understand the difference between geometrically linearisable and linearisable. 
\end{rem}

We now define an important class of surfaces for regularising Cremona group elements.
\begin{defn}\label{defn:MFS}
    Let $X$ be a surface and $B$ a variety of dimension $\leq 1$. A \emph{fibration} $\pi\colon  X \to B$ is a surjective morphism with connected fibres such that the anticanonical bundle $\omega_{X/k}^{\otimes -1}$ is $\pi$-relatively ample.  If $(X,f_X)$ and $(B,f_B)$ are $\Z/n\Z$-varieties, then a \emph{$\Z/n\Z$-fibration} $\pi \colon (X,f_X) \to (B,f_B)$ is a fibration where $\pi$ is a $\Z/n\Z$-morphism.

    Let $(X,f_X)$ be a surface. Let $\mr{NS}(X)$ denote the Néron Severi group of $X$. For a fibration $\pi\colon X\to B$ let $\mr{NS}(X/B)$ be the relative Néron Severi group of $\pi$. 
    
    A \emph{$\Z/n\Z$-Mori fibre space} is a $\Z/n\Z$-surface $(X,f_X)$ equipped with a $\Z/n\Z$-fibration $\pi\colon  (X,f_X) \to (B,f_B)$ such that $\mathrm{rk} ( \mathrm{NS}(X/B)^{\Z/n\Z}) = 1$. If $n=1$, we simply say $X$ is a \emph{Mori fibre space}.
\end{defn}

\begin{rem}
    If $B=\Spec(k)$, the fibration condition means that $X$ is a del Pezzo surface. If $B$ is a curve, then $X \to B$ is a conic bundle, and if $X$ is rational then $B \cong \mathbb{P}^1_k$. 
    
    Rational $\Z/n\Z$-Mori fibre spaces therefore split into two types: del Pezzo type and conic bundle type. When $X$ is of del Pezzo type, the condition on the Néron--Severi group $\mathrm{rank} ( \mathrm{NS}(X_{\kbar}/B_{\kbar})^{\Z/n\Z}) = 1$ implies that $(X,f_X)$ is $\Z/n\Z$-minimal. However, as in \cite[Remark 2.13]{ESS25}, it is not necessarily true that conic bundle type $\Z/n\Z$-Mori fibre spaces are $\Z/n\Z$-minimal. Nevertheless \cite[Theorem 2.6]{ESS25} shows that every rational $\Z/n\Z$-surface is $\Z/n\Z$-birational to a $\Z/n\Z$-Mori fibre space (note that if $X$ is a rational surface then $K_X$ is not nef by \cite[Theorem 2.2, Page 169]{Kollar}). In particular, restricting to rational $\Z/n\Z$-Mori fibre spaces allows us to study $\mathrm{Rat}^d_{\Z/n\Z}$.
\end{rem}

\begin{thm}[Theorem 5 of \cite{DolgachevIskovskikh2}]\label{thm:regularisablebymfs}
    Let $f \in \mathrm{Cr}_2(k)$ be an element of prime order $n$. Then there exists a $\Z/n\Z$-Mori fibre space $(X,f_X)$, such that $f_X$ regularises $f$.
\end{thm}

The following key result plays an important role in the proof of Theorem \ref{maintheorem:birational}.

\begin{thm}\label{thm:BirationalInvariance}
    Let $(X,f_X)$ be a rational $\Z/n\Z$-Mori fibre space, with $K_X^2 \leq 4$. Suppose that $(X, f_X) \dashrightarrow (Y, f_Y)$ is a $\Z/n\Z$-birational map of $\Z/n\Z$-Mori fibre spaces. Then \[\chikq([X, f_X]) = \chikq([Y, f_Y]).\]
\end{thm}
To prove this result, we use the $G$-minimal model programme to decompose the birational map $(X, f_X) \dashrightarrow (Y, f_X)$ into elementary steps, called \emph{Sarkisov links}. We then show that these elementary steps give rise to elements of $\K_1(\Var_k)$ in the kernel of $\chikq$, and therefore each elementary step does not change the value of $\chikq$.

For ease of notation, we sometimes do not write $f_X$, writing simply $X$ when discussing $\Z/n\Z$-varieties $(X,f_X)$.
\begin{defn}\label{def:linktypes} Following \cite[\S 2.2]{Iskovskikh}, \cite[Chapter 16, 16.23]{LamyCremonaBook}, we define \emph{$\Z/n\Z$-Sarkisov links} as follows. Let $X \xrightarrow{\pi_X} B_X, Y \xrightarrow{\pi_Y} B_Y$ be smooth projective rational $\Z/n\Z$-Mori fibre spaces with bases $B_X, B_Y$ respectively. A \emph{$\Z/n\Z$-Sarkisov link} $X \dashrightarrow Y$ is a $\Z/n\Z$-birational map $X \dashrightarrow Y$ such that one of the four following conditions holds:
    \begin{enumerate}
        \item The map $X \dashrightarrow Y$ is a \emph{type I $\Z/n\Z$-Sarkisov link} if $B_X = \Spec(k)$, $B_Y = \mathbb{P}^1_k$, and there exists $Z \subseteq X$, an irreducible $\Z/n\Z$-subvariety of dimension $0$, such that $Y = \mr{Bl}_Z(X)$, and $X \dashrightarrow Y$ is the inverse of the $\Z/n\Z$-birational morphism
        \[\begin{tikzcd}          X\arrow[d]&\mr{Bl}_Z(X)=Y\arrow[d]\arrow[l]\\
            \mr{Spec}(k)&\mathbb{P}^1. \arrow[l]
        \end{tikzcd}\]
        \item The map $X \dashrightarrow Y$ is a \emph{type II $\Z/n\Z$-Sarkisov link} if $B_X = B_Y$, and there exist $Z_X \subseteq X$, $Z_Y \subseteq Y$ irreducible $\Z/n\Z$-subvarieties of dimension $0$ and a commutative diagram of $\Z/n\Z$-varieties
        \[\begin{tikzcd}
            X\arrow[d]&\mr{Bl}_{Z_X}(X)=\mr{Bl}_{Z_Y}(Y)\arrow[l, "p_X"]\arrow[r, "p_Y"]&Y\arrow[d]\\
            B_X \arrow[rr,equal]&&B_Y
        \end{tikzcd}\]
        such that the $\Z/n\Z$-birational map $X \dashrightarrow Y$ is given by $p_Y \circ p_X^{-1}$.
        \item The map $X \dashrightarrow Y$ is a \emph{type III $\Z/n\Z$-Sarkisov link} if its inverse is a type I $\Z/n\Z$-Sarkisov link.
        \item The map $X \dashrightarrow Y$ is a \emph{type IV $\Z/n\Z$-Sarkisov link} if $Y=X$ and $X \dashrightarrow Y$ is the identity (i.e., a type IV $\Z/n\Z$-Sarkisov link is a change of conic bundle structure on $X$). This can be represented by the following commutative diagram
        \[\begin{tikzcd}            X\arrow[rr,equal]\arrow[d]&&Y\arrow[d]\\
            B_X\arrow[dr]&&B_Y\arrow[dl]\\
            &\mr{Spec}(k).&
        \end{tikzcd}\]

    \end{enumerate}
\end{defn}
We will also make use of the following special subset of type II Sarkisov links.
\begin{defn}\label{defn:BertiniGeiserType}
    Let $X$ be a $\Z/n\Z$ del Pezzo surface and let $X \dashrightarrow Y$ be a type II $\Z/n\Z$-Sarkisov link. We say that the Sarkisov link is \emph{of Bertini type} (resp. \emph{of Geiser type}) if there exists $Z_X \subseteq X$, $Z_Y \subseteq Y$ closed irreducible $\Z/n\Z$-subvarieties of dimension $0$, such that the Sarkisov link is given by the diagram
          \[\begin{tikzcd}
            X\arrow[d]&\mr{Bl}_{Z_X}(X)=\mr{Bl}_{Z_Y}(Y)\arrow[l]\arrow[r]&Y\arrow[d]\\
            \Spec(k)\arrow[rr,equal]&&\Spec(k)
        \end{tikzcd}\]
        where $\mr{Bl}_{Z_X}(X)$ is a del Pezzo surface of degree $1$ (resp. a del Pezzo surface of degree $2$) and if $\beta$ is the Bertini (resp. Geiser) involution on $\mathrm{Bl}_{Z_X}(X)$, then $\beta(E_{Z_X}) = E_{Z_Y}$.
\end{defn}

\begin{lem}\label{lem:BertiniGeiserIsomCentres}
    Let $X$ and $Y$ be two $\Z/n\Z$ del Pezzo surfaces, and let
    \[\begin{tikzcd}
            X\arrow[d]&\mr{Bl}_{Z_X}(X)=\mr{Bl}_{Z_Y}(Y)\arrow[l]\arrow[r]&Y\arrow[d]\\
            \Spec(k)\arrow[rr,equal]&&\Spec(k)
        \end{tikzcd}\]
        be a $\Z/n\Z$-Sarkisov link of Bertini or Geiser type. Then $Z_X \cong Z_Y$. 
\end{lem}
\begin{proof}
    Note that $E_{Z_i} \cong \mathbb{P}^1 \times_{\Spec(k)} Z_i$, and since $\beta(E_{Z_X}) = E_{Z_Y}$, we have that $E_{Z_X} \cong E_{Z_Y}$. This forces $Z_X \cong Z_Y$, as required. 
\end{proof}

We introduce Sarkisov links to apply the following theorem, which appears as \cite[7.2]{SZ21} and \cite[Theorem 16.28]{LamyCremonaBook}, and is an equivariant version of \cite[Theorem 2.5]{Iskovskikh}.
\begin{thm}\label{thm:sarklink}
    Let $\phi\colon X\dashrightarrow Y$ be a $\mb{Z}/n\mb{Z}$-birational map of rational $\mb{Z}/n\mb{Z}$-Mori fibre space surfaces. Then $\phi$ factors as a composition of $\Z/n\Z$-Sarkisov links.
\end{thm}
We also note that Sarkisov links have the following relation to $\K_1(\Var_k)$.
\begin{prop}\label{propn:linkeqns}
    Let $(X,f_X),(Y,f_Y)$ be $\Z/n\Z$-Mori fibre spaces. A $\Z/n\Z$-Sarkisov link $X \dashrightarrow Y$ gives rise to one of the following equations in $\K_1(\Var_k)$.
    \begin{enumerate}
        \item If $X \dashrightarrow Y$ is type I, then
        \[[X,f_X]-[Bl_Z(X),\tilde{f}_X]=[Z,f_X|_{Z}]-[E_Z,\tilde{f}_X|_{E_Z}] \in \K_1(\Var_k),\]
        where $Z \subseteq X$ is the irreducible $\Z/n\Z$-subvariety of dimension $0$ and $E_Z$ is the exceptional divisor of the blow-up.
        \item If $X \dashrightarrow Y$ is a type II $\Z/n\Z$-Sarkisov link, then
        \[[X,f_X]-[Y,f_Y]=[Z_X,f_X|_{Z_X}]-[Z_Y,f_Y|_{Z_Y}]+[E_{Z_X},\tilde{f}_X|_{E_{Z_X}}]-[E_{Z_Y},\tilde{f}_Y|_{E_{Z_Y}}] \in \K_1(\Var_k),\]
        where $E_{Z_i}$ is the exceptional divisor lying over $Z_i$ for $i=X,Y$.
        \item If $X \dashrightarrow Y$ is a type III $\Z/n\Z$-Sarkisov link, then we obtain the same formula as in part $1$ but with $X$s and $Y$s swapped.
        \item If $X \dashrightarrow Y$ is a type IV $\Z/n\Z$-Sarkisov link, we obtain the trivial relation.
    \end{enumerate}
\end{prop}
\begin{proof}
   Parts $(1)$ and $(3)$ follow by Lemma $\ref{cor:BlowUpRelation}$. For part (2), set $W = \mathrm{Bl}_{Z_1}X = \mathrm{Bl}_{Z_2}Y$. The result follows by rewriting
    \[
    [X,f_X]-[Y,f_Y] = \left([X,f_X]- [W,\tilde{f}_X]\right) + \left([W, \tilde{f}_Y] - [Y,f_Y]\right)
    \]
    and using Lemma $\ref{cor:BlowUpRelation}$. Part (4) is trivial, since type IV Sarkisov links only change the fibration, rather than the surface itself.
\end{proof}
    As in Proposition $\ref{prop:relationtypes}$ let $\Psi_1 = \mathrm{ker}\left(\mathrm{Rat}_{\mb{Z}/n\mb{Z}}^{2,\mr{min}} \to \mathrm{Cr}_2^{ab}(k)\right)$.
\begin{cor}\label{cor:reltype1}
Relations of type $(1)$ in Proposition $\ref{prop:relationtypes}$ are generated by the four relations appearing in Proposition $\ref{propn:linkeqns}$.
\end{cor}
\begin{proof}
    Relations of type (1) in Proposition $\ref{prop:relationtypes}$ are given by elements $[(X, f_X)] - [(X', f_{X'})]$ where $(X, f_X)$ and $(X', f_{X'})$ are $\Z/n\Z$-Mori fibre spaces which are $\Z/n\Z$-birational.

    Pick a $\Z/n\Z$-birational equivalence $\phi\colon  X \dashrightarrow X'$. Theorem $\ref{thm:sarklink}$ implies that $\phi$ factors
    \[
X = X_0 \dashrightarrow X_1 \dashrightarrow \ldots \dashrightarrow X_n = X'
    \]
    where each $X_j \dashrightarrow X_{j+1}$ is a Sarkisov link, so each $[(X_{j}, f_{j})] - [(X_{j+1}], f_{j+1})]$ is given by a Sarkisov link appearing in Proposition $\ref{propn:linkeqns}$. The result follows by rewriting
    \[
    [(X, f_X)] - [(X', f_{X'})] = \sum_{j=0}^{n-1} [(X_{j}, f_{j})] - [(X_{j+1}, f_{j+1})]. \qedhere \]
    
\end{proof}

Corollary $\ref{cor:reltype1}$ and Remark $\ref{rem:CremonaInvariants}$ gives us a strategy for proving Theorem $\ref{thm:BirationalInvariance}$: compute the value of $\chikq$ of the equations appearing in Proposition $\ref{propn:linkeqns}$, and show that $\chikq$ of the right hand side of these equations vanishes for the sorts of Sarkisov links that appear.

\begin{lem}\label{lem:bgtypek1}
    Let $(X,f_X)$, $(Y,f_Y)$ be $\Z/n\Z$-Mori fibre spaces. Let $(X,f_X) \dashrightarrow (Y,f_Y)$ be a type II $\Z/n\Z$-Sarkisov link of Bertini or Geiser type. Then $[(X,f_X)] = [(Y, f_Y)]$.
\end{lem}
\begin{proof}
    Lemma $\ref{lem:BertiniGeiserIsomCentres}$ implies that if $Z_X$, $Z_Y$ are the blow-up centres of the Sarkisov link, then $Z_X \cong Z_Y$. Moreover since the Bertini (or Geiser) involution is central in $\mathrm{Aut}( \mathrm{Bl}_{Z_X}(X))$, we see that $(E_{Z_X}, f|_{E_{Z_X}})$ is $\Z/2\Z$-equivariantly isomorphic to $(E_{Z_Y}, f|_{Z_Y})$, so $[(E_{Z_X}, i_{Z_X})] = [E_{Z_Y}, f_{Z_Y}] \in \K_1(\Var_k)$, which implies $[Z_X, f|_{Z_X}]=[Z_Y, f|_{Z_Y}]$. The result follows by the equation for type II Sarkisov links in Proposition $\ref{propn:linkeqns}$.
\end{proof}

\begin{lem}\label{lem:Type2ConicBundle}
    Suppose that $(X,f_X)$ and $(Y,f_Y)$ are $\Z/n\Z$-conic bundle surfaces, and let $X \dashrightarrow Y$ be a type II $\Z/n\Z$-Sarkisov link. Then $\chikq([X,f_X])=\chikq([Y,f_Y])$.
\end{lem}

\begin{proof}
    Recall that a type II $\Z/n\Z$-Sarkisov link is given by a diagram
     \[\begin{tikzcd}
            X\arrow[d]&\mr{Bl}_{Z_X}(X)=\mr{Bl}_{Z_Y}(Y)\arrow[l, "\pi_X"]\arrow[r, "\pi_Y"]&Y\arrow[d]\\
            \mathbb{P}^1 \arrow[rr,equal]&&\mathbb{P}^1,
        \end{tikzcd}\]
        where $Z_X$, $Z_Y$ are irreducible $\Z/n\Z$-subvarieties of $X,Y$ of dimension $0$. Applying Lemma $\ref{lem:BlowUpdetsn}$ to the equation in Proposition $\ref{propn:linkeqns}$, the result follows if $Z_X \cong Z_Y$ as irreducible $\Z/n\Z$-varieties. This follows by an identical argument to \cite[Proposition 3.1]{SchneiderRelations}, which shows that the map $\pi_Y\colon  \mathrm{Bl}_{Z_Y}(Y) \to Y$ is given by contracting the strict transform of the fibres containing $Z_X$. This result is stated over a general field without a $\Z/n\Z$-action. The argument extends identically from studying a $\Gal_k$-action to a $\Z/n\Z \times \Gal_k$-action.
\end{proof}

We are now able to prove Theorem $\ref{thm:BirationalInvariance}$.

\begin{proof}[Proof of Theorem $\ref{thm:BirationalInvariance}$]
    Theorem $\ref{thm:sarklink}$ tells us that it is sufficient to prove the result in the case where $(X,f_X) \dashrightarrow (Y,f_Y)$ is a $\Z/n\Z$-Sarkisov link. Since $K_X^2 \leq 4$ we may use \cite[Proposition 2.10]{ESS25} to reduce to one of the following cases:
    \begin{enumerate}
        \item $X$ is a del Pezzo surface of degree $4$, $Y$ is a conic bundle surface with $K_Y^2 = 3$ and $X \dashrightarrow Y$ is a $\Z/n\Z$-Sarkisov link of type I.
        \item $X$ is a conic bundle with $K_X^2 =3$, $Y$ is a del Pezzo surface of degree $4$ and $X \dashrightarrow Y$ is a $\Z/n\Z$-Sarkisov link of type III.
        \item The surfaces $X,Y$ are del Pezzo surfaces and $X \dashrightarrow Y$ is a $\Z/n\Z$-Sarkisov link of type II of Bertini or Geiser type.
        \item The surfaces $X,Y$ are both conic bundle surfaces and $X \dashrightarrow Y$ is a $\Z/n\Z$-Sarkisov link of type II.
    \end{enumerate}
    In case (1), $Y = \mathrm{Bl}_Z(X)$ and since $K_Y^2 = K_X^2 - 1$, we see that $\# Z(\kbar)=1$. This means that $Z$ consists of a single $\Spec(k)$ point which is fixed under the $\Z/2\Z$-action, and therefore the result holds by Lemma $\ref{cor:blowupfixedpoint}$. This also proves case (2), since this is the inverse. Case (3) follows by Lemma $\ref{lem:bgtypek1}$, and case (4) follows from Lemma $\ref{lem:Type2ConicBundle}$.
\end{proof}

\begin{cor}\label{cor:BirationalInvarianceCor}
    Let $f$ be a conjugacy class of elements in $\mathrm{Cr}_2(k)$ of order $n$. Suppose that there exists a $\Z/n\Z$-Mori fibre space $(X,f_X)$ which regularises $f$ with $K_X^2 \leq 4$. Then $\chikq([X,f_X])$ is independent of the choice of $\Z/n\Z$-Mori fibre space regularisation.
\end{cor}
\begin{proof}
    This is the result of Theorem $\ref{thm:BirationalInvariance}$ in the language of Cremona groups.
\end{proof}
Theorem $\ref{thm:BirationalInvariance}$ then allows us to prove Theorem \ref{maintheorem:birational}.

\begin{proof}[Proof of Theorem \ref{maintheorem:birational}]
    Recall that to prove Theorem \ref{maintheorem:birational}, we need to show that for $n$ a prime number if $f \in \mathrm{Cr}(k)[n]$ geometrically non-linearisable, then for two $\Z/n\Z$-Mori fibre space models of $f$ $(X,f_X)$ and $(X', f_{X'})$, we have $\chikq([X,f_X]) = \chikq([X', f_{X'}])$.

     Existence of $(X,f_X)$ follows from Theorem $\ref{thm:regularisablebymfs}$. Since $f$ is geometrically non-linearisable, $K_X^2 \leq 4$ by Corollary $\ref{cor:nonclosedlinear}$. The result therefore holds by Theorem $\ref{thm:BirationalInvariance}$. 
\end{proof}

\begin{defn}\label{defn:chikqCremona}
    In light of Theorem \ref{maintheorem:birational}, if $f \in \mathrm{Cr}_2(k)[n]$ where $n$ is prime and $f$ is geometrically non-linearisable, write $\chikq(f) := \chikq([X, f_X])$ where $(X,f_X)$ is a choice of $\Z/n\Z$-Mori fibre space regularisation of $f$.
\end{defn}

\begin{rem}
    If $f \in \mathrm{Cr}_2(k)[n]$ where $n$ is an odd prime, then while Theorem $\ref{thm:regularisablebymfs}$ guarantees the existence of a regularisation $(X,f_X)$, Corollary $\ref{cor:SquaresVanish}$ shows that we always have $\chikq(f) = (1,1)$, so Corollary $\ref{cor:BirationalInvarianceCor}$ is not interesting when $n$ is odd.

    When $n=2$, the result of Theorem \ref{maintheorem:birational} can be made very computable using Theorem $\ref{thm:threegeomtypes}$. Let $(X,f_X)$ denote a $\Z/2\Z$-minimal Mori fibre space regularisation of $f$. 
    
    If $f$ is conjugate to a Bertini involution, we can pick $X$ to be a del Pezzo surface of degree $1$ and then $\chikq(f) = \chikq([X, f_X])=(1,\mathrm{disc}(\chi_0(X)))$ by Corollary $\ref{cor:BertiniComputation}$. 
    
    If $f$ is conjugate to a Geiser involution, we may pick $X_i$ to be a del Pezzo surface of degree $2$ so $\chikq(f) = \chikq([X, f_X]) + (-1, 2 \cdot \mathrm{disc}(\chi_0(X)))$ by Lemma $\ref{lem:GeiserComputation}$.

    If $f$ is conjugate to a de Jonqui\`eres involution, then we can pick $X$ to be a conic bundle surface and $\chikq(f) = \chikq([X, f_X]) = ( (-1)^{\mathrm{rank}(\chi_0(X))}, \mathrm{disc}(\chi_0(X)))$ by Corollary $\ref{cor:conicbundlesn}$. 
\end{rem}

\begin{rem}\label{rem:NotBirationalInvariant}
    The condition that $f$ is geometrically non-linearisable is crucial in Theorem \ref{maintheorem:birational}, and relaxing it means the result fails to be true. For example, consider the involution $f_1 \colon  \mathbb{P}^1 \times \mathbb{P}^1 \to \mathbb{P}^1 \times \mathbb{P}^1$ given by $(x,y) \mapsto (y,x)$. Lemma $\ref{lem:surject on det}$ shows that $\chikq([\mathbb{P}^1 \times \mathbb{P}^1, f_1]) = (-1,-2)$. Since $\mathbb{P}^1 \times \mathbb{P}^1$ is a minimal surface, $(\mathbb{P}^1 \times \mathbb{P}^1, f_1)$ is trivially $\Z/2\Z$-minimal. Similarly, let $f_2\colon  \mathbb{P}^2 \to \mathbb{P}^2$ be defined by $[x,y,z] \mapsto [y,x,z]$. It is easy to see that $f_2$ acts by the identity on the Hodge cohomology of $\mathbb{P}^2$ and therefore $\chikq([\mathbb{P}^2, f_2]) = (1,1)$. Moreover, $\mathbb{P}^2$ is also a minimal surface. Defining the birational map $\phi\colon  \mathbb{P}^1 \times \mathbb{P}^1 \dashrightarrow \mathbb{P}^2$ as in Remark $\ref{rem:minimalsurface}$ by $\phi( [a,b], [c,d]) = [ad, bc, bd]$, we obtain a $\Z/2\Z$-birational map $(\mathbb{P}^1 \times \mathbb{P}^1, f_1) \dashrightarrow (\mathbb{P}^2, f_2)$. In $\K_1(\Var_k)$ language, the element $[(\mathbb{P}^1 \times \mathbb{P}^1, f_1)] - [(\mathbb{P}^2, f_2)]$ is a relation of type (1) in $\Psi_1$ which does not lie in $\mathrm{ker}(\chikq)$. 
\end{rem}

It is possibly surprising that $\chikq$ gives a birational invariant in some cases, however there are other results in the literature concerning birational invariance of related invariants. For example, in \cite[Theorem 1.1]{Batyrev}, Batyrev famously shows that if $X$ and $Y$ two smooth projective birational Calabi--Yau varieties, then they have the same compactly supported Euler characteristic. Using \cite[Remark 1.3.1]{LevineAspects} to compare the $\mathbb{A}^1$-Euler characteristic with the classical Euler characteristic gives $\mathrm{rank}(\chi_0(X)) = \mathrm{rank}(\chi_0(Y))$. In \cite[Theorem 7.3]{YauZaslow}, the same proof shows that if $X$ and $Y$ are two birational smooth projective Calabi--Yau varieties then $\mathrm{disc}(\chi_0(X)) = \mathrm{disc}(\chi_0(Y))$. It is an open question whether in this case $\chi_0(X) = \chi_0(Y)$. 
\section{Applications to the structure of the plane Cremona group}\label{sec:application}
We end the paper by giving two applications of Theorem \ref{maintheorem:birational}, inspired by the work in \cite[\S7]{CMYZ} and \cite[\S10]{CMYZ} respectively. The first application, Theorem  $\ref{thm:deJonquieres}$, uses $\chikq$ to construct an explicit family of de Jonqui\`eres involutions which are non-conjugate but are geometrically conjugate. The second application, Theorem \ref{maintheorem:fixedcurve},  gives a partial converse to \cite[Proposition 2.7]{BayleBeauville} over many non-algebraically closed fields. 

\subsection{Non-conjugate de Jonqui\`eres involutions}\label{subsection:example}
In this subsection we apply Theorem \ref{maintheorem:birational} to  produce a family of de Jonqui\`eres involutions which are non-conjugate over $k$ but become conjugate over $\kbar$, by adapting an example from \cite[Lemma 7.3]{CMYZ}.

 Following \cite[Lemma 6.2]{CMYZ}, let $Y_d$ be the surface in weighted projective space $\mathbb{P}(3,3,1,1)$ defined by the equation
 \[
 Y_d: xy = f(z,t) := (z^2-dt^2)(z^2-4dt^2)(z^2-9dt^2).
 \]
The surface $Y_d$ is equipped with a $\Z/2\Z$-action: namely, $f_d\colon [x,y,z,t] \mapsto [y,x,z,t]$. 
    
 By an identical argument to \cite[Lemma 6.2]{CMYZ}, by a slight abuse of notation, if we let $(X_d, f_d)$ denote the $\Z/2\Z$-minimal desingularisation of $(Y_d, f_d)$, then $X_d$ is a conic bundle, and $(X_d, f_d)$ is a $\Z/2\Z$-Mori fibre space by \cite[Lemma 6.2]{CMYZ}.

\begin{lem}\label{lem:deJonquieresXd}
        The surface $(X_d, f_d)$ is a $\Z/2\Z$-Mori fibre space model of de Jonqui\`eres involution of $\mathbb{P}^2$ defined by the rational map
\[
                [u,v,z] \mapsto \left[v, \frac{(u^2-d)(u^2-4d)(u^2-9d)}{v}, z\right]
\]
\end{lem}
\begin{proof}
    Applying the birational map $ Y_d \dashrightarrow \mathbb{P}^2$ given by $[x,y,z,t] \mapsto [xt^2, x, t^3]$, and using that $(X_d, f_d) \to (Y_d, f_d)$ is a $\Z/2\Z$-birational morphism gives the result.
\end{proof}
\begin{defn}
    In light of Lemma $\ref{lem:deJonquieresXd}$, let $\phi_d \in \Cr_2(k)$ be the de Jonqui\`eres involution
    \[
       [u,v,z] \mapsto \left[v, \frac{(u^2-d)(u^2-4d)(u^2-9d)}{v}, z\right].
    \]
\end{defn}

\begin{lem}\label{lem:Hirzebruch}
    Let $\mathbb{F}_3$ denote the Hirzebruch ruled surface. Then $X$ is isomorphic to $\mathrm{Bl}_Z \mathbb{F}_3$, where $Z$ is isomorphic to the disjoint union of three copies of $\Spec(k(\sqrt{d}))$.
\end{lem}
\begin{proof}
    This follows by an identical argument to \cite[Lemma 4.7]{RZ18}.
\end{proof}
\begin{cor}\label{cor:geometricallylinearisable}
    The involution $\phi_d$ is not geometrically linearisable.
\end{cor}
\begin{proof}
    Lemma $\ref{lem:Hirzebruch}$ implies $K_X^2 = 8-6 = 2$, which gives the result by Corollary $\ref{cor:nonclosedlinear}$.
\end{proof}
\begin{lem}\label{lem:deJonqutation}
    We have that $\chikq([X_d, f_d]) = (-1, d)$.
\end{lem}
\begin{proof}
    Since $\mathbb{F}_3$ is a Zariski locally trivial $\mathbb{P}^1$-bundle over $\mathbb{P}^1$, we see that 
\[\chi_0(\mathbb{F}_3) = \chi_0(\mathbb{P}^1) \cdot \chi_0(\mathbb{P}^1) = 2\mathbb{H}.\]
Lemma $\ref{lem:Hirzebruch}$ and \cite[Proposition 1.4(5)]{LevineAspects} show that  $\chi_0(X_d) = 2\mathbb{H} + 3 \langle -2 \rangle + 3 \langle -2d\rangle$, so applying Corollary $\ref{cor:conicbundlesn}$ gives the result.
\end{proof}
\begin{lem}\label{lem:conj}
    Over $k(\sqrt{d})$, the surface $(X_d, f_d)$ is $\Z/2\Z$-equivariantly isomorphic to $(X_1, f_1)$. In particular, $\phi_d$ is conjugate to $\phi_1$ when $d$ is a square in $k$.
\end{lem}
\begin{proof}
    Suppose $\sqrt{d} \in k$. Then there is an isomorphism  $Y_d \cong Y_1$ given on co-ordinates by $[x,y,z,t] \mapsto [x,y,z,\frac{1}{\sqrt{d}}t]$. Resolving singularities gives the result for $X_d$.
\end{proof}

\begin{thm}\label{thm:deJonquieres}
    In $\mathrm{Cr}_2(k)$, the de Jonqui\`eres involutions $\phi_d$ and $\phi_{d'}$ are conjugate if and only if $d = d' \in k^\times/k^{\times 2}$.
\end{thm}
\begin{proof}
The backwards direction holds by Lemma $\ref{lem:conj}$. For the forwards direction, note that $\phi_d$ is not geometrically linearisable by Corollary $\ref{cor:geometricallylinearisable}$. Therefore Theorem \ref{maintheorem:birational} tells us that $\chikq(\phi_d)$ is well defined and invariant on conjugacy classes. Lemma $\ref{lem:deJonqutation}$ computes $\chikq(\phi_d) = \chikq([X_d, f_d]) = (-1, d)$, so $\chikq(\phi_{d}) = \chikq(\phi_{d'})$ if and only if $d$ and $d'$ are equivalent modulo squares. The result follows since by Theorem \ref{maintheorem:birational}, $\chikq(\phi_d)$ is constant on conjugacy classes in $\Cr_2(k)$. 
\end{proof}

    Theorem $\ref{thm:deJonquieres}$ shows that $\chikq$ gives a genuine arithmetic invariant that can distinguish between conjugacy classes of elements of $\Cr_2(k)$ which become conjugate over $\kbar$.

\subsection{Non-conjugate involutions with the same fixed curve}\label{subsection:final}
Over $\mathbb{C}$, \cite[Proposition 2.7]{BayleBeauville} implies that an involution in $\mathrm{Cr}_2(\mathbb{C})$ is completely determined by the isomorphism class of the non-rational components of the fixed curve of the involution. Over $\mathbb{R}$, \cite[Main Corollary]{CMYZ} shows that this is not true, but over other fields this question has been open. In this subsection, we prove that over a perfect characteristic $\neq 2$ field $k$ with $k \neq \mb{F}_3, \mb{F}_5$, there are non-conjugate involutions with isomorphic non-rational fixed curve, which is precisely Theorem \ref{maintheorem:fixedcurve}. We adapt the construction of \cite[\S10]{CMYZ} to construct a family of involutions of conic bundle surfaces which all have the same fixed curve. We show that these conic bundle surfaces are $\Z/2\Z$-rational Mori fibre spaces under mild hypotheses (Lemma $\ref{lem:Xdrationality}$), which we can always arrange to satisfy over any $\neq 2$ field with a quadratic extension (Corollary $\ref{cor:anyfield}$), unless $k=\mb{F}_3$ or $\mb{F}_5$. We then compute $\chikq$ for such involutions (Corollary $\ref{cor:chikqXabd}$). Involutions on rational surfaces in turn gives us involutions in the plane Cremona group. We then apply Theorem \ref{maintheorem:birational} to the constructed involutions, which proves Theorem \ref{maintheorem:fixedcurve}. 

Fix $k$ to be such that $k^{\times}/k^{\times 2}$ is non-trivial and let $\lambda \in k^{\times}$ be such that $-\lambda$ is not a square. Let $Y_{a,b,d}$ be the surface in $\mathbb{P}^2 \times \mathbb{A}^1$ given by the equation 
\[
Y_{a,b,d}\colon x^2 + \lambda (t^2-a)(t^2-b)y^2 + (t^2-d)z^2=0\]
where $([x,y,z],t)$ are the co-ordinates on $\mathbb{P}^2 \times \mathbb{A}^1$. Let $f_{a,b,d} \colon Y_{a,b,d} \to Y_{a,b,d}$ be the involution given by $([x,y,z],t) \mapsto ([x,y,-z],t)$. We see that the projection $Y_{a,b,d} \to \mathbb{A}^1$ given by $([x,y,z],t) \mapsto t$ exhibits $Y_{a,b,d }$ as a conic bundle over $\mathbb{A}^1$. As in \cite[\S10]{CMYZ}, we may embed $Y_{a,b,d}$ into $X_{a,b,d}$ where $X_{a,b,d}$ is a smooth projective surface with a conic bundle structure $X_{a,b,d} \to \mathbb{P}^1$ where the fibre over $t=\infty$ is smooth. The involution $f_{a,b,d}$ extends to an involution on $X_{a,b,d}$, which we will abuse notation and also call $f_{a,b,d}$.

\begin{lem}\label{lem:Xdrationality}
 Suppose that $d-a$ and $d-b$ are both in $k^{\times 2}$ and that $X_{a,b,d}(k) \neq \emptyset$. Then $X_{a,b,d}$ is a rational $\Z/2\Z$-Mori fibre space.
\end{lem}
\begin{proof}
     Let $r^2 = d-a$ and let $s^2 = d-b$. The conic bundle $X_{a,b,d} \to \mathbb{P}^1$ has $6$ singular fibres, so $K_{X_{a,b,d}}^2 = 8-6 = 2$.

    Over the point $t^2=a$, the singular fibre is the conic defined by the equation $x^2 - (d-a)z^2 = 0$, which splits as the union of the $(-1)$-curves $(x-rz)=0$ and $(x+rz)=0$. These $(-1)$-curves are permuted by the $\Z/2\Z$-action, so does not contribute to the relative Néron--Severi group. Similarly over $t^2=b$ the singular fibre splits as the union of the $(-1)$-curves $(x-sz)=0$ and $(x+sz=0)$ which are permuted by the $\Z/2\Z$-action.

    Over the quadratic point $t^2=d$, the singular fibre is the conic defined by the equation $x^2 - (-\lambda) (d-a)(d-b)y^2$. Since $-\lambda$ is not a square by construction, this is a non-split singular fibre, so does not contribute to the relative Néron--Severi group. Putting this together implies that $\mathrm{rank}(\mathrm{NS}(X_{a,b,d}/\mathbb{P}^1)^{\Z/2\Z}) = 1$, so $X_{a,b,d}$ is a $\Z/2\Z$-Mori fibre space.

    It remains to show that $X_{a,b,d}$ is rational. Consider the the four $(-1)$-curves on $X_{a,b,d}$ base changed to $\kbar$, defined by the equations
    \begin{align}
        t = \sqrt{a}, & &  x-rz=0 \\
        t=-\sqrt{a}, & & x-rz=0\\
        t=\sqrt{b}, & & x-sz=0\\
        t=-\sqrt{b}, & & x-sz=0.
    \end{align}
    This is a $\Gal_k$-orbit of four disjoint $(-1)$-curves, so there exists $\tilde{X}_{a,b,d}$ which is a contraction of these $(-1)$-curves. By construction, $K_{\tilde{X}_{a,b,d}}^2 = 2+4=6$. The assumption that $X_{a,b,d}(k)$ is non-empty implies that $\tilde{X}_{a,b,d}(k) \neq \emptyset$, and therefore $\tilde{X}_{a,b,d}(k)$ is a del Pezzo surface of degree $6$ with a $k$-point, so is rational by \cite[Proposition 9.4.20]{Poonen}.
\end{proof}

\begin{lem}\label{lem:fixedcurve}
    The fixed curve of $f_{a,b,d}$ in $X_{a,b,d}$ is isomorphic to the hyperelliptic curve defined by the equation $Y^2 = \lambda(X^2 - a)(X^2-b)$. In particular, $\mathrm{Fix}(f_{a,b,d})$ does not depend on $d$ and has positive genus.
\end{lem}
\begin{proof}
    The fixed curve is the completion of the curve $z=0$, and $x^2 - \lambda(t^2-a)(t^2-b)y^2 = 0$.
   This curve lies in the affine patch $y=1$, and normalising gives the result.
\end{proof}

\begin{lem}\label{lem:discXd}
    The quadratic Euler characteristic $\chi_0(X_{a,b,d})$ is given by
    \[
    \chi_0(X_{a,b,d}) = 2\mathbb{H} +  2\langle 2\rangle + \langle 2a\rangle + \langle 2b\rangle + \langle 2\lambda \rangle + \langle 2\lambda d \rangle
    \]

\end{lem}
\begin{proof}
  By the working in Lemma $\ref{lem:Xdrationality}$, we see that $X_{a,b,d}$ is a conic bundle surface with singular fibres given by split conics over $\Spec(k(\sqrt{a})$ and $\Spec(k(\sqrt{b}))$ and a singular fibre given by a non-split conic over $\Spec(k(\sqrt{d}))$ which is split by adjoining $\sqrt{\lambda}$.

  We proceed by applying \cite[Corollary 11.8]{LevineAspects}. Let $\pi\colon X_{a,b,d} \to \mathbb{P}^1$ be the conic bundle morphism. Then define 
  \[
  \chi_0(X_{a,b,d}) = \sum_{y} \mathrm{Tr}_{k(y)/k}(\langle \delta_y \rangle) - D(\pi)\cdot \mathbb{H},
  \]
  where $y$ runs over all of the nodes in the singular fibres of $\pi$, $\delta_y \in k^{\times}/k^{\times2}$ is the non-square class such that the tangents to the node at $y$ are defined over $k(y)(\sqrt{\delta_y})$ and $D(\pi)$ is an integer which can be computed explicitly by a formula, which we can avoid doing. 

  Lemma $\ref{lem:Xdrationality}$ shows that we have three singular nodes $y_1,y_2,y_3$, with \begin{align*}
  k(y_1) = k(\sqrt{a}) \text{ and }\delta_{y_1}=1 \\
  k(y_2) = k(\sqrt{b}) \text{ and }\delta_{y_2} = 1 \\
  k(y_3) = k(\sqrt{d}) \text{ and }\delta_{y_3} = \lambda. 
  \end{align*}
  Substituting this into the equation for $\chi_0(X_{a,b,d})$ gives
  \[
  \chi_0(X_{a,b,d}) = \langle 2 \rangle + \langle 2a \rangle + \langle 2 \rangle + \langle 2b \rangle+ \langle 2\lambda \rangle+ \langle 2\lambda d \rangle  - D(\pi) \cdot \mathbb{H}.
  \]
  The result follows since $X_{a,b,d}$ is a del Pezzo surface of degree $2$ so $\mathrm{rank}(\chi_0(X_{a,b,d})) = 10$, which gives $D(\pi) = -2$.
\end{proof}
\begin{cor}\label{cor:chikqXabd}
    We can compute $\chikq([X_{a,b,d},f_{a,b,d}]) = (1, abd\lambda) \in \{\pm1\} \times k^\times/k^{\times2}$.
\end{cor}
\begin{proof}
    Follows immediately from Corollary $\ref{cor:conicbundlesn}$ and Lemma $\ref{lem:discXd}$.
\end{proof}

In light of Lemma $\ref{lem:fixedcurve}$ and Corollary $\ref{cor:chikqXabd}$, to produce elements in $\mathrm{Cr}_2(k)$ with the same fixed curve which are not conjugate, it is sufficient to find some a pair $(a,b)$ and two elements $d, d'$ with $d \neq d' \pmod{k^\times}$ such that $X_{a,b,d}$ and $X_{a,b,d'}$ are both rational. We can find these relatively easily thanks to the following two results.
\begin{lem}\label{lem:squares}
    Let $k$ be an infinite field and let $d,d' \in k^{\times}$. Then there exist infinitely many $b \in k^\times$ such that $d-b$ and $d'-b$ are both non-zero squares.
\end{lem}
\begin{proof}
    Let $\mathbb{A}^2_k$ have co-ordinates $u,v$ and consider the conic $C\colon u^2-v^2 = d'-d$. Note that if $(\alpha,\alpha') \in C(k)$, then setting $a = d-\alpha^2=d'-(\alpha')^2$ gives us that $d-a, d'-a$ are both squares. The result therefore follows if $C(k)$ is infinite. Let $\overline{C}$ be the projectivisation of $C$ in $\mathbb{P}^2$, so that $\overline{C}$ is a projective conic and note that $[1,1,0] \in \overline{C}(k)$. Therefore, $\overline{C}$ has a $k$-point, so has infinitely many $k$-points, and therefore $C(k)$ is infinite, as required.
\end{proof}

\begin{lem}\label{lem:atleast5elts}
Suppose $k \neq \mb{F}_3, \mb{F}_5$. Then there exists $s \in k^{\times}$ such that $s,s+1$ are both non-zero squares.
\end{lem}
\begin{proof}
    Recall that $\mathrm{char}(k) \neq 2$, so $k \neq \mb{F}_3, \mb{F}_5$ implies that $|k|>5$. Therefore we can pick $u\in k^{\times}$ such that $u^4 \neq 1$. Define $s$ to be the element $\frac12 (u - u^{-1})^2$. Then $s$ is a square, and $s+1 = \frac12 (u+u^{-1})^2$. Moreover, $s,s+1$ are non-zero by assumption.
\end{proof}

\begin{cor}\label{cor:anyfield}
     Suppose that $k \neq \mb{F}_3, \mb{F}_5$, and pick $s$ as in Lemma $\ref{lem:atleast5elts}$. Then there exist $b,d \in k^{\times}$ such that $d$ is not a square and $X_{-1,b,d}$ and $X_{-1,b,s}$ are both rational.

\end{cor}
\begin{proof}
    First note that $X_{-1,b,s} \neq \emptyset$ for all $b$: it has a point at $([\sqrt{s},0,1],-1) \in \mathbb{P}^2 \times \mathbb{A}^1$. Moreover, note that $s + 1$ is a square by assumption. By assumption $k^{\times}/k^{\times 2}$ is non-trivial. We claim that there exists $u \in k^{\times}$ with $u^2-1$ non-square. Let $\alpha \in k^{\times}$ be a non-square. Set $u = \frac{\alpha+1}{\alpha-1}$, so then $u^2-1 = \alpha \cdot \left( \frac{2}{\alpha-1} \right)^2$. This means $u^2-1 \equiv \alpha \pmod{k^{\times2}}$, so is a non-square. Fix $u \in k^{\times}$ such that $u^2-1$ is not a square and let $d = u^2-1$. By construction $d+1$ is a square. If $k$ is infinite, Lemma $\ref{lem:squares}$ allows us to find $b$ such that $d-b$ and $s-b$ are squares, and if $k$ is finite, we can find this by partitioning $k^{\times}$ into squares and non-squares and using the pigeonhole principle. Finally $([u,0,1],1) \in X_{-1, b,d}(k)$, so $X_{-1,b,d}$ and $X_{-1, b, s}$ are rational by Lemma $\ref{lem:Xdrationality}$. 
\end{proof}

If $X_{a,b,d}$ is rational, let $\phi_{a,b,d} \in \mathrm{Cr}_2(k)$ be the element associated to the involution $f_{a,b,d}$.

\begin{cor}\label{cor:nonconjugatesamefixedcurve}
   Let $k \neq \mb{F}_3, \mb{F}_5$, and fix $a,b \in k^{\times}$ and let $d,d' \in k^{\times}$ be such that $X_{a,b,d}$ and $X_{a,b,d'}$ are both rational and $d \neq d' \pmod{k^{\times 2}}$. Then the involutions $\phi_{a,b,d}, \phi_{a,b,d'} \in \mathrm{Cr}_2(k)$ have isomorphic fixed curves but are non-conjugate in $\mathrm{Cr}_2(k)$.
\end{cor}
\begin{proof}
    The fact that $\phi_{a,b,d}$ and $\phi_{a,b,d'}$ have isomorphic fixed curves is Lemma $\ref{lem:fixedcurve}$. Since $(X_{a,b,d}, f_{a,b,d})$ is a $\Z/2\Z$-Mori fibre space with $K_X^2=2<4$, Theorem \ref{maintheorem:birational} tells us that $\chikq(\phi_{a,b,d}) = \chikq([X_{a,b,d}, f_{a,b,d}])$ is a conjugacy class invariant of $\phi_{a,b,d}$. Corollary $\ref{cor:chikqXabd}$ computes $\chikq(\phi_{a,b,d}) = (1,abd\lambda)$. In particular, if $d \neq d' \pmod{k^{\times2}}$, then we see that $\chikq(\phi_{a,b,d}) \neq \chikq(\phi_{a,b,d'})$ and so $\phi_{a,b,d}$ is not conjugate to $\phi_{a,b,d'}$ by Theorem \ref{maintheorem:birational}.
\end{proof}

This allows us to show the final main theorem of the paper, Theorem \ref{maintheorem:fixedcurve}.

\begin{proof}[Proof of Theorem \ref{maintheorem:fixedcurve}]
Recall that we need to show that for any perfect, characteristic not $2$ field $k$ with $k \neq \mb{F}_3, \mb{F}_5$ such that $k$ admits a quadratic extension, there exist $\phi, \phi' \in \mathrm{Cr}_2(k)[2]$ which have isomorphic fixed curves but are non-conjugate in $\mathrm{Cr}_2(k)$. Since $k \neq \mb{F}_3, \mb{F}_5$, and since $k^{\times}/k^{\times 2} \neq \{1\}$ Corollary $\ref{cor:anyfield}$ shows that there exists $d \in k^{\times}$ which is not a square, $s \in k^\times$ which is a square and $b \in k$ such that $X_{-1, b,d}$ and $X_{-1, b, s}$ are both rational, and therefore $\phi_{-1,b,d}, \phi_{-1,b,s} \in \mathrm{Cr}_2(k)[2]$. Since $s$ is a square we see $d \neq s \in k^{\times}/k^{\times2}$ by assumption. The result holds by Corollary $\ref{cor:nonconjugatesamefixedcurve}$.
\end{proof}

\begin{rem}
    In Theorem \ref{maintheorem:fixedcurve}, we leverage the fact that $X_{-1,b,d}$ and $X_{-1,b,s}$ are (non-equivariantly) birational, yet $\chi_0(X_{-1,b,d}) \neq \chi_0(X_{-1,b,s})$. In particular, the proof of Theorem \ref{maintheorem:fixedcurve} leverages the fact that $\chi_0$ is not a birational invariant for these surfaces, despite the fact that Theorem \ref{maintheorem:birational} shows that $\chikq$ is an equivariant birational invariant.
\end{rem}

\bibliography{k1}{}

@article {Zak17-k1,
    AUTHOR = {Zakharevich, Inna},
     TITLE = {On {$K_1$} of an assembler},
   JOURNAL = {J. Pure Appl. Algebra},
  FJOURNAL = {Journal of Pure and Applied Algebra},
    VOLUME = {221},
      YEAR = {2017},
    NUMBER = {7},
     PAGES = {1867--1898},
      ISSN = {0022-4049,1873-1376},
   MRCLASS = {19B99 (14C35 18F30)},
  MRNUMBER = {3614978},
MRREVIEWER = {L.\ N.\ Vaserstein},
       DOI = {10.1016/j.jpaa.2016.12.028},
       URL = {https://doi.org/10.1016/j.jpaa.2016.12.028},
}

@article{Bertini,
  title={Ricerche sulle trasformazioni univoche involutorie nel piano},
  author={Bertini, Eugenio},
  journal={Annali di Matematica Pura ed Applicata (1867-1897)},
  volume={8},
  number={1},
  pages={244--286},
  year={1877},
  publisher={Springer}
}

@article {Campbell,
    AUTHOR = {Campbell, Jonathan A.},
     TITLE = {The {$K$}-theory spectrum of varieties},
   JOURNAL = {Trans. Amer. Math. Soc.},
  FJOURNAL = {Transactions of the American Mathematical Society},
    VOLUME = {371},
      YEAR = {2019},
    NUMBER = {11},
     PAGES = {7845--7884},
      ISSN = {0002-9947,1088-6850},
   MRCLASS = {19D99},
  MRNUMBER = {3955537},
MRREVIEWER = {Jens\ Hornbostel},
       DOI = {10.1090/tran/7648},
       URL = {https://doi.org/10.1090/tran/7648},
}

@article {Cheltsov,
    AUTHOR = {Cheltsov, I. A.},
     TITLE = {Regularization of birational automorphisms},
   JOURNAL = {Mat. Zametki},
  FJOURNAL = {Matematicheskie Zametki},
    VOLUME = {76},
      YEAR = {2004},
    NUMBER = {2},
     PAGES = {286--299},
      ISSN = {0025-567X,2305-2880},
   MRCLASS = {14E07},
  MRNUMBER = {2098998},
MRREVIEWER = {M.\ Kh.\ Gizatullin},
       DOI = {10.1023/B:MATN.0000036763.78939.9f},
       URL = {https://doi.org/10.1023/B:MATN.0000036763.78939.9f},
}

@incollection {DolgachevIskovskikh,
    AUTHOR = {Dolgachev, Igor V. and Iskovskikh, Vasily A.},
     TITLE = {Finite subgroups of the plane {C}remona group},
 BOOKTITLE = {Algebra, arithmetic, and geometry: in honor of {Y}u. {I}.
              {M}anin. {V}ol. {I}},
    SERIES = {Progr. Math.},
    VOLUME = {269},
     PAGES = {443--548},
 PUBLISHER = {Birkh\"auser Boston, Boston, MA},
      YEAR = {2009},
      ISBN = {978-0-8176-4744-5},
   MRCLASS = {14E07 (14J26)},
  MRNUMBER = {2641179},
MRREVIEWER = {Julie\ D\'eserti},
       DOI = {10.1007/978-0-8176-4745-2\_11},
       URL = {https://doi.org/10.1007/978-0-8176-4745-2_11},
}

@article {CMYZ,
    AUTHOR = {Cheltsov, Ivan and Mangolte, Fr\'ed\'eric and Yasinsky, Egor
              and Zimmermann, Susanna},
     TITLE = {Birational involutions of the real projective plane},
   JOURNAL = {J. Eur. Math. Soc. (JEMS)},
  FJOURNAL = {Journal of the European Mathematical Society (JEMS)},
    VOLUME = {28},
      YEAR = {2026},
    NUMBER = {8},
     PAGES = {3653--3711},
      ISSN = {1435-9855,1435-9863},
   MRCLASS = {14E07 (14E08 14L30)},
  MRNUMBER = {5062490},
       DOI = {10.4171/jems/1537},
       URL = {https://doi.org/10.4171/jems/1537},
}

@article {Iskovskikh,
    AUTHOR = {Iskovskikh, V. A.},
     TITLE = {Factorization of birational mappings of rational surfaces from
              the point of view of {M}ori theory},
   JOURNAL = {Uspekhi Mat. Nauk},
  FJOURNAL = {Uspekhi Matematicheskikh Nauk},
    VOLUME = {51},
      YEAR = {1996},
    NUMBER = {4(310)},
     PAGES = {3--72},
      ISSN = {0042-1316,2305-2872},
   MRCLASS = {14E05 (14E07 14M20)},
  MRNUMBER = {1422227},
MRREVIEWER = {M.\ Kh.\ Gizatullin},
       DOI = {10.1070/RM1996v051n04ABEH002962},
       URL = {https://doi.org/10.1070/RM1996v051n04ABEH002962},
}

@article {LinShinder,
    AUTHOR = {Lin, Hsueh-Yung and Shinder, Evgeny},
     TITLE = {Motivic invariants of birational maps},
   JOURNAL = {Ann. of Math. (2)},
  FJOURNAL = {Annals of Mathematics. Second Series},
    VOLUME = {199},
      YEAR = {2024},
    NUMBER = {1},
     PAGES = {445--478},
      ISSN = {0003-486X,1939-8980},
   MRCLASS = {14E07 (14C35 19A99)},
  MRNUMBER = {4681148},
MRREVIEWER = {Kenneth\ A.\ Brown},
       DOI = {10.4007/annals.2024.199.1.6},
       URL = {https://doi.org/10.4007/annals.2024.199.1.6},
}

@article {Riouappendix,
    AUTHOR = {Levine, Marc and Yang, Yaping and Zhao, Gufang and Riou,
              Jo\"el},
     TITLE = {Algebraic elliptic cohomology theory and flops {I}},
   JOURNAL = {Math. Ann.},
  FJOURNAL = {Mathematische Annalen},
    VOLUME = {375},
      YEAR = {2019},
    NUMBER = {3-4},
     PAGES = {1823--1855},
      ISSN = {0025-5831,1432-1807},
   MRCLASS = {14E30 (55N22 55N34)},
  MRNUMBER = {4023393},
MRREVIEWER = {Lennart\ Meier},
       DOI = {10.1007/s00208-019-01880-x},
       URL = {https://doi.org/10.1007/s00208-019-01880-x},
}

@book{Dolgachev,
  author    = {Dolgachev, Igor V.},
  title     = {Classical Algebraic Geometry: A Modern View},
  publisher = {Cambridge University Press},
  address   = {Cambridge},
  year      = {2012},
  isbn      = {978-1-107-01765-8},
  doi       = {10.1017/CBO9781139084437},
  url       = {https://doi.org/10.1017/CBO9781139084437}
}

@article{Prokhorov,
  title={The rationality problem for conic bundles},
  author={Prokhorov, Yu G},
  journal={Russian Mathematical Surveys},
  volume={73},
  number={3},
  pages={375--456},
  year={2018},
  publisher={London Mathematical Society, Turpion Ltd and the Russian Academy of Sciences}
}

@article {LamySchneider,
    AUTHOR = {Lamy, St\'ephane and Schneider, Julia},
     TITLE = {Generating the plane {C}remona groups by involutions},
   JOURNAL = {Algebr. Geom.},
  FJOURNAL = {Algebraic Geometry},
    VOLUME = {11},
      YEAR = {2024},
    NUMBER = {1},
     PAGES = {111--162},
      ISSN = {2313-1691,2214-2584},
   MRCLASS = {14E07},
  MRNUMBER = {4680015},
MRREVIEWER = {Alexandr\ V.\ Pukhlikov},
       DOI = {10.14231/ag-2024-004},
       URL = {https://doi.org/10.14231/ag-2024-004},
}

@article {ZakhCampbell,
    AUTHOR = {Campbell, Jonathan A. and Zakharevich, Inna},
     TITLE = {D\'evissage and localization for the {G}rothendieck spectrum
              of varieties},
   JOURNAL = {Adv. Math.},
  FJOURNAL = {Advances in Mathematics},
    VOLUME = {411},
      YEAR = {2022},
     PAGES = {Paper No. 108710, 80},
      ISSN = {0001-8708,1090-2082},
   MRCLASS = {19E08 (14C35 18N25)},
  MRNUMBER = {4512395},
MRREVIEWER = {Vivek\ Sadhu},
       DOI = {10.1016/j.aim.2022.108710},
       URL = {https://doi.org/10.1016/j.aim.2022.108710},
}

@article {Annihilator,
    AUTHOR = {Zakharevich, Inna},
     TITLE = {The annihilator of the {L}efschetz motive},
   JOURNAL = {Duke Math. J.},
  FJOURNAL = {Duke Mathematical Journal},
    VOLUME = {166},
      YEAR = {2017},
    NUMBER = {11},
     PAGES = {1989--2022},
      ISSN = {0012-7094,1547-7398},
   MRCLASS = {14C35 (18F30 19E99)},
  MRNUMBER = {3694563},
MRREVIEWER = {Paolo\ Aluffi},
       DOI = {10.1215/00127094-0000016X},
       URL = {https://doi.org/10.1215/00127094-0000016X},
}

@article {Zak17-assemblers,
    AUTHOR = {Zakharevich, Inna},
     TITLE = {The {$K$}-theory of assemblers},
   JOURNAL = {Adv. Math.},
  FJOURNAL = {Advances in Mathematics},
    VOLUME = {304},
      YEAR = {2017},
     PAGES = {1176--1218},
      ISSN = {0001-8708,1090-2082},
   MRCLASS = {19M05 (14C35 18F25 18F30 55N15)},
  MRNUMBER = {3558230},
MRREVIEWER = {Kevin\ Hutchinson},
       DOI = {10.1016/j.aim.2016.08.045},
       URL = {https://doi.org/10.1016/j.aim.2016.08.045},
}

@incollection{Morel2004,
  author    = {Morel, Fabien},
  title     = {An introduction to {$\mathbb{A}^1$}-homotopy theory},
  booktitle = {Contemporary Developments in Algebraic {$K$}-Theory},
  editor    = {Karoubi, Max and Kuku, Aderemi O. and Pedrini, Claudio},
  series    = {ICTP Lecture Notes},
  volume    = {15},
  pages     = {357--441},
  publisher = {Abdus Salam International Centre for Theoretical Physics},
  address   = {Trieste},
  year      = {2004}
}

@incollection {MT08,
    AUTHOR = {Muro, Fernando and Tonks, Andrew},
     TITLE = {On {$K_1$} of a {W}aldhausen category},
 BOOKTITLE = {{$K$}-theory and noncommutative geometry},
    SERIES = {EMS Ser. Congr. Rep.},
     PAGES = {91--115},
 PUBLISHER = {Eur. Math. Soc., Z\"urich},
      YEAR = {2008},
      ISBN = {978-3-03719-060-9},
   MRCLASS = {19B99 (19D10 55U35)},
  MRNUMBER = {2513334},
MRREVIEWER = {Charles\ Weibel},
       DOI = {10.4171/060-1/4},
       URL = {https://doi.org/10.4171/060-1/4},
}

@misc{Nan24,
      title={Weight Filtrations and Derived Motivic Measures}, 
      author={Anubhav Nanavaty},
      year={2024},
      howpublished={arXiv:2401.06879},
      archivePrefix={arXiv},
      primaryClass={math.AG},
      url={https://arxiv.org/abs/2401.06879}, 
}

@article {LPLS24,
    AUTHOR = {Levine, Marc and Pepin Lehalleur, Simon and Srinivas,
              Vasudevan},
     TITLE = {Euler characteristics of homogeneous and weighted-homogeneous
              hypersurfaces},
   JOURNAL = {Adv. Math.},
  FJOURNAL = {Advances in Mathematics},
    VOLUME = {441},
      YEAR = {2024},
     PAGES = {Paper No. 109556, 86},
      ISSN = {0001-8708,1090-2082},
   MRCLASS = {14F42 (14B05 14C15 14C30)},
  MRNUMBER = {4708158},
       DOI = {10.1016/j.aim.2024.109556},
       URL = {https://doi.org/10.1016/j.aim.2024.109556},
}

@article {AMBOWZ22,
    AUTHOR = {Arcila-Maya, Niny and Bethea, Candace and Opie, Morgan and
              Wickelgren, Kirsten and Zakharevich, Inna},
     TITLE = {Compactly supported {$\mathbb{A}^1$}-{E}uler characteristic and
              the {H}ochschild complex},
   JOURNAL = {Topology Appl.},
  FJOURNAL = {Topology and its Applications},
    VOLUME = {316},
      YEAR = {2022},
     PAGES = {Paper No. 108108, 24},
      ISSN = {0166-8641,1879-3207},
   MRCLASS = {14F42 (13D03 19E15 55M05)},
  MRNUMBER = {4438953},
MRREVIEWER = {Daniel\ Dugger},
       DOI = {10.1016/j.topol.2022.108108},
       URL = {https://doi.org/10.1016/j.topol.2022.108108},
}

@article{Ananyevskiy2016,
  author  = {Ananyevskiy, Alexey},
  title   = {On the relation of special linear algebraic cobordism to {W}itt groups},
  journal = {Homology, Homotopy and Applications},
  year    = {2016},
  volume  = {18},
  number  = {1},
  pages   = {205--230},
  doi     = {10.4310/HHA.2016.v18.n1.a11}
}

@article{Ng,
  title={$ \mathrm{K}_1(\mathrm{Var}) $ is presented by stratified birational equivalences},
  author={Ng, Ming},
  journal={arXiv preprint arXiv:2510.20433},
  year={2025}
}

@ARTICLE{Ron25,
  title     = {The {G}rothendieck ring of varieties and algebraic {K}-theory of
               spaces},
  author    = {R{\"o}ndigs, Oliver},
  journal   = {Open Book Ser.},
  publisher = {Mathematical Sciences Publishers},
  volume    = {6},
  number    = {1},
  pages     = {165--196},
  month     = {Mar},
  year      = {2025},
}

@article {RSO19,
    AUTHOR = {R\"ondigs, Oliver and Spitzweck, Markus and {\O}stv{\ae}r, Paul
              Arne},
     TITLE = {The first stable homotopy groups of motivic spheres},
   JOURNAL = {Ann. of Math. (2)},
  FJOURNAL = {Annals of Mathematics. Second Series},
    VOLUME = {189},
      YEAR = {2019},
    NUMBER = {1},
     PAGES = {1--74},
      ISSN = {0003-486X,1939-8980},
   MRCLASS = {14F42 (55Q45)},
  MRNUMBER = {3898173},
MRREVIEWER = {Markus\ Szymik},
       DOI = {10.4007/annals.2019.189.1.1},
       URL = {https://doi.org/10.4007/annals.2019.189.1.1},
}

@article {HSS17,
    AUTHOR = {Hoyois, Marc and Scherotzke, Sarah and Sibilla, Nicol\`o},
     TITLE = {Higher traces, noncommutative motives, and the categorified
              {C}hern character},
   JOURNAL = {Adv. Math.},
  FJOURNAL = {Advances in Mathematics},
    VOLUME = {309},
      YEAR = {2017},
     PAGES = {97--154},
      ISSN = {0001-8708,1090-2082},
   MRCLASS = {14F05 (14F42 18D05 19D55)},
  MRNUMBER = {3607274},
MRREVIEWER = {Satoshi\ Mochizuki},
       DOI = {10.1016/j.aim.2017.01.008},
       URL = {https://doi.org/10.1016/j.aim.2017.01.008},
}

@misc{stacks,
    shorthand    = {Stacks},
    author       = {The {Stacks Project Authors}},
    title        = {\textit{Stacks Project}},
    howpublished = {\url{https://stacks.math.columbia.edu}},
    year         = {2018},
  }

@article{KLMMS24,
      title={Scissors automorphism groups I: Homological stability and K-theory}, 
      author={Alexander Kupers and Ezekiel Lemann and Cary Malkiewich and Jeremy Miller and Robin J. Sroka},
      year={(to appear)},
      journal={J. Reine Angew. Math.},
}

@article {Bit04,
    AUTHOR = {Bittner, Franziska},
     TITLE = {The universal {E}uler characteristic for varieties of
              characteristic zero},
   JOURNAL = {Compos. Math.},
  FJOURNAL = {Compositio Mathematica},
    VOLUME = {140},
      YEAR = {2004},
    NUMBER = {4},
     PAGES = {1011--1032},
      ISSN = {0010-437X,1570-5846},
   MRCLASS = {14F42 (14C35 14E05)},
  MRNUMBER = {2059227},
MRREVIEWER = {Thomas\ Geisser},
       DOI = {10.1112/S0010437X03000617},
       URL = {https://doi.org/10.1112/S0010437X03000617},
}

@article {Hoy14,
    AUTHOR = {Hoyois, Marc},
     TITLE = {A quadratic refinement of the
              {G}rothendieck-{L}efschetz-{V}erdier trace formula},
   JOURNAL = {Algebr. Geom. Topol.},
  FJOURNAL = {Algebraic \& Geometric Topology},
    VOLUME = {14},
      YEAR = {2014},
    NUMBER = {6},
     PAGES = {3603--3658},
      ISSN = {1472-2747,1472-2739},
   MRCLASS = {14F42 (11E81 55M20)},
  MRNUMBER = {3302973},
MRREVIEWER = {Daniel\ C.\ Isaksen},
       DOI = {10.2140/agt.2014.14.3603},
       URL = {https://doi.org/10.2140/agt.2014.14.3603},
}

@article {LR20,
    AUTHOR = {Levine, Marc and Raksit, Arpon},
     TITLE = {Motivic {G}auss-{B}onnet formulas},
   JOURNAL = {Algebra Number Theory},
  FJOURNAL = {Algebra \& Number Theory},
    VOLUME = {14},
      YEAR = {2020},
    NUMBER = {7},
     PAGES = {1801--1851},
      ISSN = {1937-0652,1944-7833},
   MRCLASS = {14F42 (55N20 55N35)},
  MRNUMBER = {4150251},
MRREVIEWER = {Matthias\ Wendt},
       DOI = {10.2140/ant.2020.14.1801},
       URL = {https://doi.org/10.2140/ant.2020.14.1801},
}

@book {Sch85,
    AUTHOR = {Scharlau, Winfried},
     TITLE = {Quadratic and {H}ermitian forms},
    SERIES = {Grundlehren der mathematischen Wissenschaften [Fundamental
              Principles of Mathematical Sciences]},
    VOLUME = {270},
 PUBLISHER = {Springer-Verlag, Berlin},
      YEAR = {1985},
     PAGES = {x+421},
      ISBN = {3-540-13724-6},
   MRCLASS = {11Exx (11-02 12D15 15A63 16A16 16A28)},
  MRNUMBER = {770063},
MRREVIEWER = {R.\ Ware},
       DOI = {10.1007/978-3-642-69971-9},
       URL = {https://doi.org/10.1007/978-3-642-69971-9},
}

@book {Lam05,
    AUTHOR = {Lam, T. Y.},
     TITLE = {Introduction to quadratic forms over fields},
    SERIES = {Graduate Studies in Mathematics},
    VOLUME = {67},
 PUBLISHER = {American Mathematical Society, Providence, RI},
      YEAR = {2005},
     PAGES = {xxii+550},
      ISBN = {0-8218-1095-2},
   MRCLASS = {11Exx},
  MRNUMBER = {2104929},
MRREVIEWER = {K.\ Szymiczek},
       DOI = {10.1090/gsm/067},
       URL = {https://doi.org/10.1090/gsm/067},
}

@misc{Lieb21,
	author        = {Joshua Lieber},
	title         = {From Six Functors Formalisms to Derived Motivic Measures},
	eprint        = {arxiv:2109.07065},
	archivePrefix = {arXiv},
	primaryClass  = {math.AG},
	year          = {2021}
}

@article {Kuijper,
    AUTHOR = {Kuijper, Josefien},
     TITLE = {A descent principle for compactly supported extensions of
              functors},
   JOURNAL = {Ann. K-Theory},
  FJOURNAL = {Annals of K-Theory},
    VOLUME = {8},
      YEAR = {2023},
    NUMBER = {3},
     PAGES = {489--529},
      ISSN = {2379-1683,2379-1691},
   MRCLASS = {14F42 (18F10 55U99)},
  MRNUMBER = {4632374},
MRREVIEWER = {Daniel\ Isaksen},
       DOI = {10.2140/akt.2023.8.489},
       URL = {https://doi.org/10.2140/akt.2023.8.489},
}

@article{BD17,
	author  = {Bondarko, Mikhail V. and D{\'e}glise, Fr{\'e}d{\'e}ric},
	title   = {Dimensional Homotopy $t$-Structures in Motivic Homotopy Theory},
	journal = {Advances in Mathematics},
	volume  = {311},
	pages   = {91--189},
	year    = {2017},
	doi     = {10.1016/j.aim.2017.02.003}
}

@article{CWZ19,
	author  = {Campbell, Jonathan and Wolfson, Jesse and Zakharevich, Inna},
	title   = {Derived {$\ell$}-adic Zeta Functions},
	journal = {Advances in Mathematics},
	volume  = {354},
	pages   = {106760},
	year    = {2019},
	doi     = {10.1016/j.aim.2019.106760}
}

@article {BGT,
    AUTHOR = {Blumberg, Andrew J. and Gepner, David and Tabuada, Gon\c calo},
     TITLE = {A universal characterization of higher algebraic {$K$}-theory},
   JOURNAL = {Geom. Topol.},
  FJOURNAL = {Geometry \& Topology},
    VOLUME = {17},
      YEAR = {2013},
    NUMBER = {2},
     PAGES = {733--838},
      ISSN = {1465-3060,1364-0380},
   MRCLASS = {19D10 (18D20 19D25 19D55 55N15 55U40)},
  MRNUMBER = {3070515},
MRREVIEWER = {Ross\ Staffeldt},
       DOI = {10.2140/gt.2013.17.733},
       URL = {https://doi.org/10.2140/gt.2013.17.733},
}

@article {MT1type,
    AUTHOR = {Muro, Fernando and Tonks, Andrew},
     TITLE = {The 1-type of a {W}aldhausen {$K$}-theory spectrum},
   JOURNAL = {Adv. Math.},
  FJOURNAL = {Advances in Mathematics},
    VOLUME = {216},
      YEAR = {2007},
    NUMBER = {1},
     PAGES = {178--211},
      ISSN = {0001-8708,1090-2082},
   MRCLASS = {19B99 (18G50 18G55)},
  MRNUMBER = {2353254},
MRREVIEWER = {Luca\ Barbieri Viale},
       DOI = {10.1016/j.aim.2007.05.008},
       URL = {https://doi.org/10.1016/j.aim.2007.05.008},
}

@misc{Zakharevich2022,
      title={$E_\infty$-Ring structures on the $K$-theory of assemblers and point counting}, 
      author={Inna Zakharevich},
      year={2022},
      eprint={2209.08640},
      archivePrefix={arXiv},
      primaryClass={math.KT},
      url={https://arxiv.org/abs/2209.08640}, 
}

@article{Hoy17,
	author  = {Hoyois, Marc},
	title   = {The six operations in equivariant motivic homotopy theory},
	journal = {Advances in Mathematics},
	volume  = {305},
	year    = {2017},
	pages   = {197--279},
	doi     = {10.1016/j.aim.2016.09.031}
}

@article {Bass,
    AUTHOR = {Bass, Hyman},
     TITLE = {Clifford algebras and spinor norms over a commutative ring},
   JOURNAL = {Amer. J. Math.},
  FJOURNAL = {American Journal of Mathematics},
    VOLUME = {96},
      YEAR = {1974},
     PAGES = {156--206},
      ISSN = {0002-9327,1080-6377},
   MRCLASS = {15A66 (10C30 15A63 20G35)},
  MRNUMBER = {360645},
MRREVIEWER = {W.\ Scharlau},
       DOI = {10.2307/2373586},
       URL = {https://doi.org/10.2307/2373586},
}

@book{Hartshorne,
  author    = {Hartshorne, Robin},
  title     = {Algebraic Geometry},
  series    = {Graduate Texts in Mathematics},
  volume    = {52},
  publisher = {Springer-Verlag},
  address   = {New York},
  year      = {1977},
  doi       = {10.1007/978-1-4757-3849-0}
}

@article{BayleBeauville,
  author  = {Bayle, Lionel and Beauville, Arnaud},
  title   = {Birational involutions of {$\mathbb{P}^2$}},
  journal = {Asian Journal of Mathematics},
  year    = {2000},
  volume  = {4},
  number  = {1},
  pages   = {11--18},
  doi     = {10.4310/AJM.2000.v4.n1.a2}
}

@article{SZ21,
  author  = {Schneider, Julia and Zimmermann, Susanna},
  title   = {Algebraic subgroups of the plane {Cremona} group
             over a perfect field},
  journal = {{\'E}pijournal de G{\'e}om{\'e}trie Alg{\'e}brique},
  volume  = {5},
  year    = {2021},
  note    = {Article 14},
  doi     = {10.46298/epiga.2021.6715},
  url     = {https://epiga.episciences.org/8709}
}

@book{Kollar,
  author    = {Koll{\'a}r, J{\'a}nos},
  title     = {Rational Curves on Algebraic Varieties},
  series    = {Ergebnisse der Mathematik und ihrer Grenzgebiete. 3. Folge},
  volume    = {32},
  publisher = {Springer-Verlag},
  address   = {Berlin},
  year      = {1996},
  isbn      = {978-3-540-60168-5},
  doi       = {10.1007/978-3-662-03276-3},
  url       = {https://doi.org/10.1007/978-3-662-03276-3}
}

@misc{ESS25,
  author        = {Elagin, Alexey and Schneider, Julia and Shinder, Evgeny},
  title         = {Atomic decompositions for derived categories of {$G$}-surfaces},
  year          = {2025},
  eprint        = {2512.05064},
  archivePrefix = {arXiv},
  primaryClass  = {math.AG},
  doi           = {10.48550/arXiv.2512.05064},
  url           = {https://arxiv.org/abs/2512.05064}
}

@article{YauZaslow,
  author  = {Pajwani, Jesse and P{\'a}l, Ambrus},
  title   = {An arithmetic {Yau--Zaslow} formula},
  journal = {Compositio Mathematica},
  year    = {2026},
  volume  = {162},
  number  = {4},
  pages   = {792--851},
  doi     = {10.1017/S0010437X26102887}
}

@article{LevineAspects,
  author  = {Levine, Marc},
  title   = {Aspects of enumerative geometry with quadratic forms},
  journal = {Documenta Mathematica},
  year    = {2020},
  volume  = {25},
  pages   = {2179--2239},
  doi     = {10.4171/DM/797}
}

@unpublished{LamyCremonaBook,
  author = {Lamy, St{\'e}phane},
  title  = {The {Cremona} Group},
  year   = {2026},
  note   = {Submitted manuscript, Cremona Book Project},
  url    = {https://www.math.univ-toulouse.fr/~slamy/blog/cremona.html}
}

@article{OrmsbyOstvaer,
title = {Stable motivic $\pi_1$ of low-dimensional fields},
journal = {Advances in Mathematics},
volume = {265},
pages = {97-131},
year = {2014},
issn = {0001-8708},
doi = {https://doi.org/10.1016/j.aim.2014.07.024},
url = {https://www.sciencedirect.com/science/article/pii/S0001870814002576},
author = {Kyle M. Ormsby and Paul Arne {\O}stv{\ae}r}
}

@incollection{Batyrev,
  author    = {Batyrev, Victor V.},
  title     = {Birational {Calabi--Yau} {$n$}-folds have equal {Betti} numbers},
  booktitle = {New Trends in Algebraic Geometry},
  editor    = {Hulek, Klaus and Reid, Miles and Peters, Chris and Catanese, Fabrizio},
  series    = {London Mathematical Society Lecture Note Series},
  volume    = {264},
  pages     = {1--12},
  publisher = {Cambridge University Press},
  address   = {Cambridge},
  year      = {1999},
  doi       = {10.1017/CBO9780511721540.002}
}

@misc{BG21,
  author        = {Braunling, Oliver and Groechenig, Michael},
  title         = {Classes in {Zakharevich} {$K$}-groups constructed from {Quillen} {$K$}-theory},
  year          = {2021},
  eprint        = {2109.01136},
  archiveprefix = {arXiv},
  primaryclass  = {math.KT},
  doi           = {10.48550/arXiv.2109.01136},
  url           = {https://arxiv.org/abs/2109.01136}
}

@article{PaninWalter2019,
  author  = {Panin, Ivan and Walter, Charles},
  title   = {On the motivic commutative ring spectrum {$\mathbf{BO}$}},
  journal = {St. Petersburg Mathematical Journal},
  year    = {2019},
  volume  = {30},
  number  = {6},
  pages   = {933--972},
  doi     = {10.1090/spmj/1578}
}

@article{DolgachevIskovskikh2,
  author  = {Dolgachev, Igor V. and Iskovskikh, Vasily A.},
  title   = {On elements of prime order in the plane {Cremona}
             group over a perfect field},
  journal = {International Mathematics Research Notices},
  volume  = {2009},
  number  = {18},
  pages   = {3467--3485},
  year    = {2009},
  doi     = {10.1093/imrp/rnp061}
}

@incollection{Prokhorov2,
  author    = {Prokhorov, Yuri},
  title     = {Finite groups of birational transformations},
  booktitle = {European Congress of Mathematics},
  publisher = {EMS Press},
  year      = {2023},
  pages     = {413--437},
  doi       = {10.4171/8ECM/20}
}

@article{SchneiderRelations,
  author  = {Schneider, Julia},
  title   = {Relations in the {Cremona} group over a perfect field},
  journal = {Annales de l'Institut Fourier},
  volume  = {72},
  number  = {1},
  pages   = {1--42},
  year    = {2022},
  doi     = {10.5802/aif.3463}
}

@inproceedings{Voevodsky,
  author    = {Voevodsky, Vladimir},
  title     = {{$\mathbb{A}^1$}-homotopy theory},
  booktitle = {Proceedings of the International Congress of Mathematicians,
               Vol. I (Berlin, 1998)},
  series    = {Documenta Mathematica},
  number    = {Extra Vol. I},
  year      = {1998},
  pages     = {579--604},
  doi       = {10.4171/DMS/1-1/21},
  url       = {https://doi.org/10.4171/DMS/1-1/21}
}

@misc{Broring,
  author        = {Br{\"o}ring, Louisa F.},
  title         = {Quadratic {Euler} Characteristic of Geometrically Cyclic Branched Coverings},
  year          = {2026},
  eprint        = {2605.13425},
  archiveprefix = {arXiv},
  primaryclass  = {math.AG},
  doi           = {10.48550/arXiv.2605.13425},
  url           = {https://arxiv.org/abs/2605.13425}
}

@article{Viergever,
  author  = {Viergever, Anna M.},
  title   = {The quadratic {Euler} characteristic of a smooth projective same-degree complete intersection},
  journal = {Annals of K-Theory},
  year    = {2025},
  volume  = {10},
  number  = {4},
  pages   = {597--652},
  doi     = {10.2140/akt.2025.10.597}
}

@book{Shafarevich,
  author     = {Shafarevich, Igor R.},
  title      = {Basic Algebraic Geometry 1: Varieties in Projective Space},
  edition    = {3},
  translator = {Reid, Miles},
  publisher  = {Springer},
  address    = {Heidelberg},
  year       = {2013},
  doi        = {10.1007/978-3-642-37956-7},
  isbn       = {978-3-642-37955-0}
}

@article{Schlichting,
  author        = {Schlichting, Marco},
  title         = {The {Mayer--Vietoris} principle for
                   {Grothendieck--Witt} groups of schemes},
  journal       = {Inventiones Mathematicae},
  volume        = {179},
  number        = {2},
  pages         = {349--433},
  year          = {2010},
  doi           = {10.1007/s00222-009-0219-1},
  eprint        = {0811.4632},
  archivePrefix = {arXiv},
  primaryClass  = {math.KT}
}

@incollection{DoldPuppe1999,
  author    = {Dold, Albrecht and Puppe, Dieter},
  title     = {Duality, trace and transfer},
  booktitle = {Proceedings of the International Conference on
               Geometric Topology (Warsaw, 1978)},
  publisher = {PWN},
  address   = {Warsaw},
  year      = {1980},
  pages     = {81--102}
}

@article{RZ18,
  author  = {Robayo, Maria Fernanda and Zimmermann, Susanna},
  title   = {Infinite algebraic subgroups of the real {Cremona} group},
  journal = {Osaka Journal of Mathematics},
  year    = {2018},
  volume  = {55},
  number  = {4},
  pages   = {681--712},
  doi     = {10.18910/70816}
}

@article{Zimmermann2018,
  author  = {Zimmermann, Susanna},
  title   = {The Abelianization of the real {Cremona} group},
  journal = {Duke Mathematical Journal},
  volume  = {167},
  number  = {2},
  pages   = {211--267},
  year    = {2018},
  doi     = {10.1215/00127094-2017-0028}
}

@book{Poonen,
  author    = {Poonen, Bjorn},
  title     = {Rational Points on Varieties},
  series    = {Graduate Studies in Mathematics},
  volume    = {186},
  publisher = {American Mathematical Society},
  address   = {Providence, RI},
  year      = {2017},
  isbn      = {978-1-4704-3773-2},
  doi       = {10.1090/gsm/186}
}

@article{BHSVW,
  author  = {Bilu, Margaret and Ho, Wei and Srinivasan, Padmavathi
             and Vogt, Isabel and Wickelgren, Kirsten},
  title   = {Quadratic enrichment of the logarithmic derivative
             of the zeta function},
  journal = {Transactions of the American Mathematical Society, Series B},
  volume  = {11},
  pages   = {1183--1225},
  year    = {2024},
  doi     = {10.1090/btran/201}
}

@article{Azo,
 author = {Azouri, Ran},
 title = {Motivic {Euler} characteristic of nearby cycles and a generalised quadratic conductor formula},
 journal = {Journal of Algebra},
 volume = {667},
 pages = {109--164},
 year = {2025},
 doi = {10.1016/j.jalgebra.2024.12.010},
}
\bibliographystyle{alpha}
\end{document}